\documentclass[a4paper,final]{elsarticle}

\xdef\bibDir{Bibliography}
\xdef\gfxpath{Figures}

\xdef\LDDTPRfigSource{\gfxpath/Tikzsource}
\xdef\LDDTPTPfigSource{\gfxpath/Tikzsource/LDD-TwoPhase}

\xdef\LDDTPRcaptionPrefix{LDD-TPR: }

\newif\ifdevelopmentVersion

\newif\ifshowlinenumbers

\usepackage{ifluatex}
\ifluatex
    \usepackage{lualatex-math}
    \usepackage{shellesc} 
\else
    \PassOptionsToPackage{utf8}{inputenx} 
    \usepackage[utf8]{inputenx}
    \input{ix-utf8enc.dfu}
\fi
\PassOptionsToPackage{T1}{fontenc} 
\usepackage{fontenc}
\PassOptionsToPackage{ngerman,american,british}{babel} 
\usepackage{babel}

\usepackage[full]{textcomp} 
\usepackage{eucal}  
\usepackage{mathrsfs} 
\PassOptionsToPackage{fleqn}{amsmath}       
\usepackage{amsmath}
\usepackage{amsfonts}
\usepackage{stmaryrd}
\usepackage{exscale}
\usepackage{amstext}

\usepackage{amsthm}
\usepackage{thm-patch,aliasctr,parseargs,keyval}
\usepackage{thmtools}

\usepackage{bm} 
\usepackage{graphicx}
\graphicspath{{\gfxpath/}} 
\usepackage[a4paper,left = 1.5 cm, right = 1.5 cm,top = 2 cm, bottom = 2.5cm]{geometry}
\usepackage{scrhack} 
\usepackage{xspace}  
\usepackage{mparhack} 
\PassOptionsToPackage{printonlyused,smaller}{acronym}
  \usepackage{acronym} 

\usepackage[shortlabels]{enumitem} 

\usepackage{tabularx} 
\usepackage{makecell} 
\usepackage{booktabs} 
\usepackage{multirow}
\usepackage[format=plain, indention=2ex]{caption} 
\usepackage{wrapfig}
\PassOptionsToPackage{dvipsnames*,svgnames,table}{xcolor}
\usepackage{xcolor}
\usepackage{subcaption}

\usepackage{hhline} 

\definecolor{mediumblue}{RGB}{0,0,205}
\definecolor{navyblue}{RGB}{0,0,128}
\definecolor{midnightblue}{RGB}{25,25,112}
\definecolor{royalblue4}{RGB}{39,64,139}
\definecolor{blue3}{RGB}{0,0,205}
\definecolor{steelblue3}{RGB}{79,148,205}
\definecolor{steelblue4}{RGB}{54,100,139}
\definecolor{uniStuttgartHellblau}{RGB}{0,190,255}
\definecolor{uniStuttgartMittelblau}{RGB}{0,81,158}

\definecolor{brown}{RGB}{165,42,42} 
\definecolor{brown3}{RGB}{205,51,51} 
\definecolor{brown4}{RGB}{139,35,35} 
\definecolor{red3}{RGB}{205,0,0} 
\definecolor{tomato}{RGB}{205,79,57} 
\definecolor{firebrick3}{RGB}{205,38,38}
\definecolor{firebrick4}{RGB}{139,26,26}

\definecolor{gold}{RGB}{255,215,0}
\definecolor{gold3}{RGB}{238,201,0} 
\definecolor{darkgoldenrod1}{RGB}{255,185,15}
\definecolor{goldenrod1}{RGB}{255,193,37}
\definecolor{goldenrod}{RGB}{218,165,32}

\definecolor{orange}{rgb}{.9,.6,.1}
\definecolor{dunkelorange}{rgb}{.9,.5,.0}
\definecolor{orange2}{RGB}{238,154,0}
\definecolor{orange3}{RGB}{205,133,0}

\definecolor{darkgreen}{RGB}{0,100,0}
\definecolor{green3}{RGB}{0,205,0} 
\definecolor{olivedrab}{RGB}{107,142,35}
\definecolor{olivedrab2}{RGB}{179,238,58}
\definecolor{olivedrab3}{RGB}{154,205,50}
\definecolor{forestgreen}{RGB}{34,139,34}
\definecolor{darkolivegreen}{RGB}{85,107,47}
\definecolor{darkolivegreen4}{RGB}{110,139,61}
\definecolor{khaki3}{RGB}{205,198,115}

\definecolor{grey}{rgb}{0.5,0.5,0.5}
\definecolor{dimgrey}{RGB}{105,105,105}
\definecolor{dimgrey2}{RGB}{153,153,153}
\definecolor{dimgrey3}{RGB}{181,181,181}
\definecolor{lightgrey}{RGB}{211,211,211}
\definecolor{lightergrey}{RGB}{201,201,201}
\definecolor{verylightgrey}{RGB}{222,222,222}
\definecolor{uniStuttgartAnthrazit}{RGB}{62,68,76}

\definecolor{braun}{rgb}{.6,.5,.1}
\definecolor{braun2}{rgb}{.6,.4,.1}
\definecolor{tan4}{RGB}{149,90,43}

\definecolor{wheat3}{RGB}{205,186,150}
\definecolor{wheat4}{RGB}{139,126,102}

\makeatletter
  \def\thmt@headstyle@margincolored{%
    \makebox[0pt][r]{\color{black!40}\NUMBER\ }\NAME\NOTE
  }
\makeatother

\newcommand{\definitionstylecolor}{black}
\newcommand{\exampleandremarkstylecolor}{black}
\newcommand{\lammacorollarystylecolor}{black}
\newcommand{\notestylecolor}{black}
\newcommand{\customtheoremstylecolor}{black}
\declaretheoremstyle[
    spaceabove=6pt, 
    spacebelow=6pt,
    headfont={\color{\definitionstylecolor}\bfseries },
    notefont=\normalfont, notebraces={(}{)},
    bodyfont=\itshape,
    postheadspace=1em,
]{customdefintion}
\declaretheorem[name=Definition,numberwithin=section,style=customdefintion,
refname={Definition,Definitions}]{definition}

\declaretheorem[name=Notation,sibling=definition,style=customdefintion,
refname={Notation,Notations}]{notation}
\declaretheorem[name=Notations,numbered=no,style=customdefintion,
refname={Notations,Notations}]{notations*}
\declaretheorem[name=Assumption,sibling=definition,style=customdefintion,
refname={Assumption,Assumptions}]{assumption}

\declaretheoremstyle[
spaceabove=6pt, spacebelow=6pt,
headfont={\color{\exampleandremarkstylecolor}\bfseries},
notefont=\normalfont, notebraces={(}{)},
bodyfont=\normalfont,
postheadspace=1em,
]{custombspbem} 

\declaretheorem[name=Remark,sibling=definition,style=custombspbem,
refname={Remark,Remarks}]{remark}

\declaretheoremstyle[
spaceabove=6pt, spacebelow=6pt,
headfont={\color{\lammacorollarystylecolor}\bfseries},
notefont=\normalfont, notebraces={(}{)},
bodyfont=\itshape,
postheadspace=1em,
]{customlemma}
\declaretheorem[name=Lemma,sibling=definition,style=customlemma,
refname={Lemma,Lemmata}]{lemma}

\declaretheoremstyle[
spaceabove=6pt, spacebelow=6pt,
headfont={\color{\notestylecolor}\bfseries},
notefont=\normalfont, notebraces={(}{)},
bodyfont=\itshape,
postheadspace=1em,
]{customnotiz}

\declaretheorem[name=Problem,sibling=definition,style=customnotiz,
refname={Problem,Problems}]{problem}
\declaretheorem[name=Problem,numbered=no,style=customnotiz,
refname={Problem,Problems}]{problem*}
\declaretheorem[name=Note,sibling=definition,style=customnotiz,
refname={Note,Notes}]{note}

\declaretheoremstyle[
spaceabove=6pt, spacebelow=6pt,
headfont={\color{darkgreen}\bfseries},
notefont=\normalfont, notebraces={(}{)},
bodyfont=\normalfont,
postheadspace=1em,
]{customberechnung}

\declaretheoremstyle[
spaceabove=6pt, spacebelow=6pt,
headfont={\color{\customtheoremstylecolor}\bfseries},
notefont=\normalfont, notebraces={(}{)},
bodyfont=\itshape,
postheadspace=1em,
]{customsatz}

\declaretheorem[name=Theorem,sibling=definition,style=customsatz,
refname={Theorem,Theorems}]{theorem}
\declaretheorem[name=Theorem,numbered=no,style=customsatz,
refname={Theorem,Theorems}]{theorem*}

\declaretheoremstyle[
spaceabove=6pt, spacebelow=6pt,
headfont=\bfseries,
notefont=\normalfont, notebraces={(}{)},
bodyfont=\normalfont,
postheadspace=1em,
]{custombew}

\ifdevelopmentVersion
    \usepackage{shellesc}
    \usepackage{tikz}
    \usepackage{pgfplots}
    \usepackage{pgfplotstable}
    \pgfplotscreateplotcyclelist{thesisSeus}{
      uniStuttgartMittelblau,densely dashdotdotted,every mark/.append style={solid,fill=uniStuttgartMittelblau!70},mark=halfsquare right*\\
      uniStuttgartMittelblau,densely dashdotted,every mark/.append style={solid,fill=uniStuttgartMittelblau!70},mark=triangle*\\
      uniStuttgartAnthrazit,every mark/.append style={fill=uniStuttgartAnthrazit!70},mark=square*\\
      uniStuttgartAnthrazit,densely dashdotdotted,every mark/.append style={solid,fill=uniStuttgartAnthrazit!70},mark=pentagon*\\
      uniStuttgartHellblau!100!black,every mark/.append style={fill=uniStuttgartHellblau},mark=halfsquare*\\
      uniStuttgartHellblau!100!black,densely dashed,every mark/.append style={solid,fill=uniStuttgartHellblau!70},mark=oplus*\\
      uniStuttgartMittelblau,mark=star\\
      uniStuttgartAnthrazit,every mark/.append style={fill=uniStuttgartAnthrazit!70},mark=diamond*\\
      uniStuttgartHellblau!100!black,densely dashdotted,every mark/.append style={solid,fill=uniStuttgartHellblau!70},mark=*\\
      uniStuttgartMittelblau,densely dashed,every mark/.append style={solid,fill=uniStuttgartMittelblau!70},mark=square*\\
      uniStuttgartAnthrazit,densely dashed,every mark/.append style={solid,fill=uniStuttgartAnthrazit!70},mark=triangle*\\
      uniStuttgartHellblau!100!black,densely dashed,every mark/.append style={solid},mark=star\\
      uniStuttgartMittelblau,densely dashed,every mark/.append style={solid,fill=uniStuttgartMittelblau!70},mark=diamond*\\
      uniStuttgartAnthrazit,densely dotted,thick,every mark/.append style={solid,fill=uniStuttgartAnthrazit!70},mark=triangle\\
      uniStuttgartHellblau!100!black,densely dashed,every mark/.append style={solid},mark=halfsquare left*\\
      uniStuttgartMittelblau,densely dashdotdotted,every mark/.append style={solid,fill=uniStuttgartMittelblau!70},mark=+\\
    }

    \usepackage{xfp}
    \pgfplotsset{
      compat=newest,%
      cycle list name={thesisSeus},%
    }
    \usetikzlibrary{arrows,%
                    arrows.meta,%
                    petri,%
                    topaths,%
                    fit,%
                    positioning,%
                    decorations.pathmorphing,%
                    backgrounds,
                    calc,
                    patterns%
    }%

    \usetikzlibrary{external}
    \usepgfplotslibrary{external}
    \tikzset{%
      external/optimize=true,
      }%
    \tikzset{external/shell escape={-shell-escape}}
      \tikzset{
	png export/.style={
	  external/system call/.add={}{; convert -density 300 -transparent white "\image.pdf" "\image.png"},
	  /pgf/images/external info,
	  /pgf/images/include external/.code={%
	    \includegraphics[width=\pgfexternalwidth,height=\pgfexternalheight]{##1.png}%
	  },
	}
      }

    \tikzifexternalizing{%
    }{%
    \usepackage{pdfpages}
    }%
\fi

\usepackage{listings}
\usepackage{siunitx}

\usepackage{microtype}

\usepackage{hyperref}
\hypersetup{%
  unicode=true,
  colorlinks=true, linktocpage=true, pdfstartpage=3, pdfstartview=FitV,%
  breaklinks=true, pageanchor=true,%
  pdfpagemode=UseNone, %
  plainpages=false, bookmarksnumbered, bookmarksopen=true, bookmarksopenlevel=1,%
  hypertexnames=true, pdfhighlight=/O,
  urlcolor=CTurl, linkcolor=CTlink, citecolor=CTcitation, 
  pdfsubject={},%
  pdfkeywords={},%
  pdfcreator={pdfLaTeX},%
  pdfproducer={LaTeX with hyperref}%
}

\makeatletter
\@ifpackageloaded{babel}%
  {%
    \addto\extrasamerican{%
    }%
    \addto\extrasngerman{%
    }%
    }{\relax}
\makeatother

\ifshowlinenumbers
  \usepackage{lineno}
  \modulolinenumbers[5]
\fi

\allowdisplaybreaks[1] 

\providecommand{\mLyX}{L\kern-.1667em\lower.25em\hbox{Y}\kern-.125emX\@}

\newcommand{\spl}{\bigl\langle}
\newcommand{\spr}{\bigr\rangle}
\newcommand{\tspl}{\langle}
\newcommand{\tspr}{\rangle}
\newcommand{\bspl}{\Bigl\langle}
\newcommand{\bspr}{\Bigr\rangle}

\newcommand{\restrictTo}[2]{{#1}{}{|_{#2}}}

\newcommand{\Richards}{Richards }
\newcommand{\richards}{Richards}
\newcommand{\TPRfull}{two-phase–\Richards }

\newcommand{\TPR}{TP–R }
\newcommand{\TPTP}{TP–TP }
\newcommand{\tpr}{TP–R}
 
\newcommand{\RR}{\mathbb{R}}
\newcommand{\RRd}{\mathbb{R}^d}
\newcommand{\NN}{\mathbb{N}}

\newcommand{\krond}[2][\alpha]{\delta_{#1 #2}}
\providecommand{\abs}[1]{\lvert#1\rvert}
\providecommand{\norm}[1]{\lVert#1\rVert}

\newcommand{\tr}[1]{\gamma_{#1}}

\newcommand{\completion}[2]{\overline{#1}{}^{#2}}

\DeclareMathOperator{\dv}{\bm{\nabla}\cdot}

\DeclareMathOperator{\Div}{div}

\newcommand{\nw}{n\!w}

\newcommand{\interface}{\Gamma}
\newcommand{\dom}[1]{\Omega_{#1}}
\newcommand{\domT}[1]{\Omega_{#1,T}}
\newcommand{\closDom}[1]{\overline{\dom{#1}}}
\newcommand{\spaceTimeDom}[1]{\dom{#1} \times (0,T)}
\newcommand{\spaceTimeInterface}{\interface \times (0,T)}
\newcommand{\union}[2]{#1 \cup \, #2}
\newcommand{\subdomNum}{{W}}
\newcommand{\vt}[1]{\bm{#1}}
\newcommand{\nv}{\bm{n}}
\newcommand{\ind}{\mathcal{I}}
\newcommand{\indR}{\ind^{R}}
\newcommand{\indTP}{\ind^{T\!P}}
\newcommand{\indTPR}{\ind^{T\!P}_R}
\newcommand{\indTPRl}[1]{\ind^{T\!P}_{#1,R}}
\newcommand{\indRTP}{\ind_{T\!P}^R}
\newcommand{\indRTPl}[1]{\ind_{#1,T\!P}^{R}}

\newcommand{\addnullcolor}{black}

\newcommand{\oB}[1]{\partial\Omega_#1 \cap \partial\Omega}

\newcommand{\porosity}[1]{\Phi_{#1}}
\newcommand{\intf}[1]{\Gamma_{#1}}

\newcommand{\paln}[2][l]{p_{\alpha,#1}^{#2}}
\newcommand{\palnh}[2][l]{p_{\alpha,#1}^{h,#2}}
\newcommand{\pwln}[2][l]{p_{w,#1}^{#2}}

\newcommand{\pnwln}[2][l]{p_{\nw,#1}^{#2}}

\newcommand{\palni}[2][l]{p_{\alpha,#1}^{n,#2}}
\newcommand{\palnih}[2][l]{p_{\alpha,#1}^{h,#2}}
\newcommand{\pwlni}[2][l]{p_{w,#1}^{n,#2}}

\newcommand{\pnwlni}[2][l]{p_{\nw,#1}^{n,#2}}

\newcommand{\Sln}[2][l]{\porosity{#1}S_{#1}^{#2}}

\newcommand{\Slni}[2][l]{\porosity{#1}S_{#1}^{n,#2}}

\newcommand{\dtSR}[1]{\porosity{#1}\partial_t S_{#1}(p_a,\pwln[#1]{})}
\newcommand{\dtS}[1]{\porosity{#1}\partial_t  S_{#1}(\pnwln[#1]{},\pwln[#1]{})}
\newcommand{\gradPalnPlusGravity}[2][l]{\bm{\nabla} \bigl( \paln[#1]{#2} + z_{\alpha} \bigr)}
\newcommand{\gradPnwlnPlusGravity}[2][l]{\bm{\nabla} \bigl( \pnwln[#1]{#2} + z_{\nw} \bigr)}
\newcommand{\gradPwlnPlusGravity}[2][l]{\bm{\nabla} \bigl( \pwln[#1]{#2} + z_w \bigr)}
\newcommand{\gradPalniPlusGravity}[2][l]{\bm{\nabla} \bigl( \palni[#1]{#2} + z_{\alpha} \bigr)}

\newcommand{\flux}[2][l]{\vt{F_{#2,#1}}}
\newcommand{\fluxn}[2][l]{\vt{F_{#2,#1}^{n}}}
\newcommand{\fluxnMinusOne}[2][l]{\vt{F_{#2,#1}^{n-1}}}
\newcommand{\fluxb}[2][]{\vt{F_{#2,l}^{#1}}}
\newcommand{\fluxai}[2][l]{\vt{F_{\alpha,#1}^{n,#2}}}
\newcommand{\fluxaih}[2][l]{\vt{F_{\alpha,#1}^{h,#2}}}
\newcommand{\fluxwi}[2][l]{\vt{F_{w,#1}^{n,#2}}}
\newcommand{\fluxgi}[2][l]{\vt{F_{\nw,#1}^{n,#2}}}
\newcommand{\fluxnwi}[2][l]{\vt{F_{\nw,#1}^{n,#2}}}

\newcommand{\kalni}[2][l]{k_{\alpha,#1}^{n,#2}}
\newcommand{\kwlni}[2][l]{k_{w,#1}^{n,#2}}

\newcommand{\kaln}[2][l]{k_{\alpha,#1}^{#2}}
\newcommand{\kwln}[2][l]{k_{w,#1}^{#2}}

\newcommand{\knwln}[2][l]{k_{\nw,#1}^{#2}}

\newcommand{\gali}[2][l]{g_{\alpha,#1}^{#2}}
\newcommand{\galih}[2][l]{g_{\alpha,#1}^{h,#2}}
\newcommand{\gwli}[2][l]{g_{w,#1}^{#2}}

\newcommand{\gnwli}[2][l]{g_{\nw,#1}^{#2}}
\newcommand{\gwl}[1][l]{g_{w,#1}}
\newcommand{\gnwl}[1][l]{g_{\nw,#1}}

\newcommand{\gal}[1][l]{g_{\alpha,#1}}

\newcommand{\epali}[2][l]{e_{p,#1}^{\alpha,#2}}
\newcommand{\epwli}[2][l]{e_{p,#1}^{w,#2}}
\newcommand{\epnwli}[2][l]{e_{p,#1}^{\nw,#2}}
\newcommand{\gradepali}[2][l]{\bm{\nabla}e_{p,#1}^{\alpha,#2}}

\newcommand{\egali}[2][l]{e_{g,#1}^{\alpha,#2}}
\newcommand{\egwli}[2][l]{e_{g,#1}^{w,#2}}
\newcommand{\egnwli}[2][l]{e_{g,#1}^{\nw,#2}}
\newcommand{\epaliOnGamma}[2][l]{e_{p,#1}^{\alpha,#2}{}}
\newcommand{\epwliOnGamma}[2][l]{e_{p,#1}^{w,#2}{}}
\newcommand{\epnwliOnGamma}[2][l]{e_{p,#1}^{\nw,#2}{}}
\newcommand{\egaliOnGamma}[2][l]{e_{g,#1}^{\alpha,#2}{}}
\newcommand{\egwliOnGamma}[2][l]{e_{g,#1}^{w,#2}{}}
\newcommand{\egnwliOnGamma}[2][l]{e_{g,#1}^{\nw,#2}{}}

\newcommand{\testfunc}[2][l]{\varphi_{#2,#1}}

\newcommand{\Fs}{\mathcal{V}}
\newcommand{\dFs}{\Fs_1 \times \Fs_2}
\newcommand{\prodFs}[2][]{\widehat{\Fs}_{#1}^{#2}}
\newcommand{\Tracespace}{H^{1/2}_{00}(\Gamma)}

\newcommand{\Ts}[2][00]{H^{1/2}_{#1}(\Gamma_{#2})}
\newcommand{\dualTs}[2][00]{H^{1/2}_{#1}(\Gamma_{#2})'}
\newcommand{\tGamma}{{\Gamma}}
\newcommand{\onGamma}[1]{{#1}{{\mid}_\Gamma}}

\newcommand{\Fenics}{{\scshape Fenics }}
\newcommand{\fenics}{{\scshape Fenics}}
\newcommand{\Dolfin}{{\scshape Dolfin }}
\newcommand{\dolfin}{{\scshape Dolfin}}

\newcommand{\mshr}{{\scshape Mshr}}
\newcommand{\Code}[1]{\lstinline!#1! }
\newcommand{\code}[1]{\lstinline!#1!}

\newcommand{\meshrez}{\lstinline!mesh_resolution!}

\newcommand{\Python}{{\scshape Python}}

\newcommand{\LDD}{{\scshape LDD }}

\newcommand{\ldd}{{\scshape LDD}}

\newcommand{\LDDTPTP}{{\scshape \ldd-\TPTP}}
\newcommand{\LDDTPR}{{\scshape \ldd-\TPR}}

\newcommand{\mesh}{\mathcal{T}}
\newcommand{\submesh}[1]{\mathcal{T}_{#1}}
\newcommand{\FEMp}{\mathcal{P}_1\varLambda^0}
\newcommand{\FEMgl}{\mathcal{P}_1\varLambda^2}

\providecommand{\DIFdel}[1]{} 

\journal{Journal of \LaTeX\ Templates}

\begin{document}
\begin{frontmatter}
\title{Towards Hybrid Two-Phase Modelling Using Linear Domain Decomposition}


\author[davidaddress]{David Seus\corref{correspondingAuthor}}
\address[davidaddress]{
  Institute of Applied Analysis and Numerical Simulation, University of Stuttgart, 
  Pfaffenwaldring 57, 
  70569 Stuttgart, 
  Germany
}

\author[florinaddress]{Florin A. Radu}
\address[florinaddress]{
  Department of Mathematics, 
  University of Bergen, 
  P.O. Box 7800, 
  N-5020 Bergen, 
  Norway
}
\author[davidaddress]{Christian Rohde}


\begin{abstract} 
  The viscous  flow  of two  immiscible fluids  in a porous medium on the Darcy scale is governed by a system of nonlinear parabolic equations. If infinite mobility of one phase can be assumed (e.g.~in soil layers in contact with the atmosphere)
the system can be substituted by the scalar Richards model. Thus, the domain of the porous medium  may be  partitioned into  disjoint subdomains with 
either the full two-phase or the simplified Richards model dynamics.  
Extending the one-model approach  from  \cite{SeusMitra2018, SeusEnumath2019} we suggest coupling conditions  for  this hybrid model
approach.  Based on an Euler implicit discretisation, a linear iterative (-type) domain decomposition scheme is proposed, and proven to be convergent.  
The theoretical findings are verified by  a comparative numerical study  that  in particular confirms the efficiency of the hybrid ansatz as compared to full two-phase model computations. 
\end{abstract}


\end{frontmatter}


\section{Introduction} 
Multiphase flow  through porous media occurs   for a wide variety of natural and technical processes.
Examples in  soil-related environmental sciences comprise   enhanced oil recovery,  remediation of contaminated soils, $CO_2$ storage  or   evaporation processes in the vadose zone. In the technological realm we mention the  design of filters, fuel cells or damping materials.
Mathematical modelling and numerical simulation are essential tools for the understanding of multiphase flow processes. However,  due  to  varying   material properties or changing 
flow regimes the  governing equations can become  strongly  heterogeneous leading to severe mathematical and in particular computational problems. To  meet these challenges  domain 
decomposition methods are an established approach (see e.g.~\cite{Lions1988}). The basic 
idea is to split the  domain in subdomains such that 
each of these subdomains can be equipped with its own model 
and numerical solver. Following an iterative scheme and by 
construction of  analytically and numerically appropriate 
coupling conditions an approximate solution on the original mono-domain can then be recovered.

In contrast to existing approaches  for homogeneous  two-phase flow  modellings, the purpose of  the present contribution is the development and analysis  of a non-overlapping domain decomposition method for \textit{hybrid} two-phase flow modellings. We consider for a porous domain on the Darcy scale the dynamics of two  incompressible and viscous fluids, that are assumed to be immiscible. Let the fluids be denoted  as the wetting ($w$) and the nonwetting ($nw$) phase, respectively. The domain is decomposed into subdomains with the flow either governed 
by the full two-phase (TP) model or by the simpler  Richards (R) model. The latter  applies e.g.\ for high mobilities of the nonwetting phase. The partition might  come  along with  changes in the relative permeability functions, fluid viscosities and densities, as well as in porosities and intrinsic permeabilities.  The major advantage of the hybrid approach is the possible gain of computing time that 
can be obtained when substituting the  full two-phase model system by the  approximative scalar Richards equation on parts of the domain.

First, in Section  \ref{LDD-TPR:section:problemFormulation}, we  present  coupling conditions for the hybrid  \TPR  model across the interfaces of subdomains. In fact, 
the coupling condition for  the nonwetting flux in  the two-phase model  is not at all obvious, given that on the Richards model domain there is no equation for the nonwetting phase. This leads to an unmatched number
of unknowns on the different subdomains. We therefore introduce  two different  coupling conditions depending on the  (non)occurrence of gravitational forces. 
Extending   our approach  for  homogeneous two-phase flow models in \cite{SeusMitra2018,SeusEnumath2019}, we proceed then  in Section \ref{LDD-TPR:section:problemFormulation} with the time-discrete problem and introduce   a  domain decomposition solver 
based  on  simultaneous L-scheme linearisation, see \cite{Pop2004, List2016}. The resulting scheme is called \LDDTPR solver.   We  provide a consistency result that ensures in the case of convergence of the \LDDTPR solver that the mono-domain  solution is recovered (Lemma \ref{lemma:LDD-TPR:iterationinterpretationlemma}).  Section \ref{LDD-TPR:section:ConvergenceOfScheme}  contains  the  core analytical result of the paper, that is the   
 convergence of the \LDDTPR solver in Theorem \ref{LDD-TPR:convergence:theorem}. 
 The idea of the proof is based on bounding  the  series of iteration errors which implies that the sequence of iteration errors must vanish.
 A key ingredient to achieve this is to detect matching  interface Robin-type terms such that telescopic sums are obtained. In fact, the latter is only possible if the pressure traces that are part of the Robin-type coupling condition on interfaces act as functionals via the $\Tracespace$-scalar product and not in the classical way via the dual pairing $H^{-1/2}(\interface)\times \Tracespace$,
  cf.\ Remark \ref{rem_represent} as well as \cite{Agranovich2015}.
 The convergence is guaranteed under a restriction on the time-step size which reduces to the restrictions obtained in \cite{SeusMitra2018,SeusEnumath2019} for the respective single-model cases.
 
 To limit  the notational overhead  and  to keep the focus, 
 Sections \ref{LDD-TPR:section:problemFormulation},
 \ref{LDD-TPR:section:ConvergenceOfScheme} are 
 restricted to a two-domain partition. In  Section \ref{LDD-TPR:section:problemFormulation:multiDomain}
we generalize  the  \LDDTPR solver  to  
a multi-domain situation.  
Finally, Section \ref{LDD-TPR:section:NumericalSection} provides  the validation  of the performance of the  \LDDTPR solver,  and displays  simulations on two- and multi-domain partitions for realistic soil parameters. The experiments confirm the convergence statement from Theorem \ref{LDD-TPR:convergence:theorem}  revealing linear rates. We then analyse  the influence of numerical and solver  parameters
(mesh size,  time step, Robin parameters, L-scheme parameters). For the 
multi-domain case we  focus on gravity effects. Most importantly we 
show the  advantage    of the hybrid model approach in terms of computational efficiency,  
as compared to the use of  the full two-phase flow model  on the entire mono-domain.
The paper ends with an outlook how the \LDDTPR solver
can be utilised for an error-controlled model-adaptive approach.

We conclude this introduction with a short  overview on the  literature for  related domain-decomposition methods and  solvers  for multiphase flow in porous media.
Independently of the underlying numerical approaches, domain decomposition methods allow 
to reduce the computational complexity of the problem, and to  follow parallel solver techniques. We refer to  \cite{QuarteroniValli2005,Dolean2015} for  general descriptions of the field. Optimising the parameters in the 
transmission conditions is an important issue in all
domain decomposition methods, see  e.g.\ \cite{Benne2016} and references therein.  What concerns porous media flow on the Darcy scale, we 
refer to  \cite{Skogestad2013} for an overview of different overlapping domain decomposition strategies.
Turning to two-phase flow,  a combined non-overlapping domain decomposition method and multigrid solver approach for the Richards equation has been put forward in \cite{Berninger2015}.
In  \cite{DIPIETRO2014163,Ahmed2019} algorithms for
multiphase porous media flow are introduced, including a-posteriori estimates to optimise the parameters and the number of iterations. A time-adaptive domain decomposition concept is pursued  in \cite{KURAZ20142}.
Convergence of a Schwarz waveform relaxation method is established in \cite{gander2021nonoverlapping} for 
the transport equation in the fractional flow formulation of two-phase flow.
Lunowa et al. in \cite{Lunowa2020} apply ideas from \cite{SeusMitra2018,SeusEnumath2019}   for  a dynamic capillary pressure model with hysteresis on a two-domain substructuring. The work \cite{AHMED2020113294} is concerned with two-phase flow with discontinuous capillary pressures. None of these works address
the case of a hybrid model ansatz.

We combine the domain-decomposition method for each time step with an L-scheme (see \cite{Pop2004, List2016})
to linearize the complete system. This linearisation approach, which is a stabilised Picard method  has been used for a variety of applications, e.g. nonlinear poromechanics \cite{BorregalesReveron2021} or fully coupled flow and transport  \cite{Illiano2021}. The L-scheme has the advantage of not involving the computation of derivatives in contrast to the Newton or the modified Picard method. Moreover, its implementation is very easy, it is globally convergent and the linear problems that need to be solved within each iteration are much better conditioned as the ones steming from e.g.\ the Newton method, see \cite{List2016}. Nevertheless, a drawback of L-schemes is their slower (linear) convergence in comparison to Newton's  scheme. Albeit faster converging L-schemes
have been suggested in \cite{List2016,MitraPop2019},  this article adheres to the standard L-scheme, focussing  on  an \LDD scheme for a flexible, subdomain-wise combination of the \Richards equation and the full two-phase flow model.

\section{Two-phase flow models and the  \LDDTPR solver for the two-domain case} 
\label{LDD-TPR:section:problemFormulation}
\subsection{Coupling the full two-phase flow model 
 with the \Richards model: the \TPR model}
\begin{figure}[t]
    \centering
    \ifdevelopmentVersion
      \input{\LDDTPRfigSource/DD-illustration.tikz}
    \else
      \includegraphics[width=.35\textwidth]{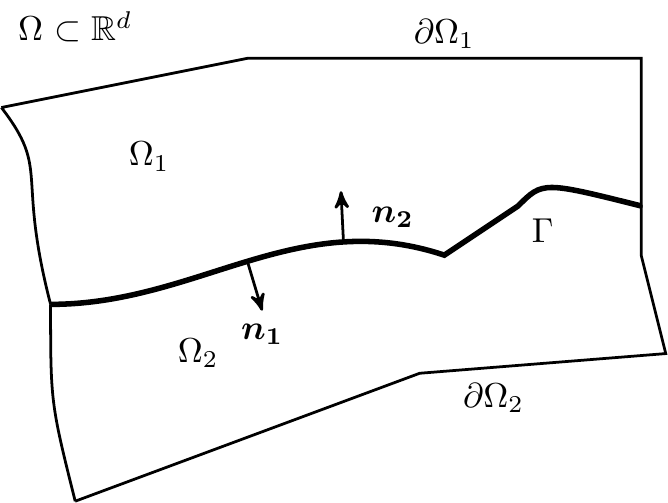}
    \fi
    \caption[Illustration of a two-domain layered soil]{Illustration of a layered soil domain $\dom{} = \dom{1} \cup \dom{2} \cup \interface$ $\subset \RRd$ with fixed interface
    $\interface$. Also shown are the normal vectors along the interface. \label{DD-two-domain-illustration} }
    \end{figure}

Let a Lipschitz domain $\dom{} \subset \RR^d$, $d\in\{2,3\}$, be decomposed into two non-overlapping Lipschitz subdomains 
$\dom{1}, \dom{2} \subset \RR^d$  such that $\dom{} = \dom{1} \cup \interface \cup \dom{2} $,  with  
 $\interface := (\closDom{1} \cap \closDom{2})\setminus \partial\Omega$ being the interface. The latter is assumed to be a $(d-1)$-dimensional Lipschitz manifold. By $\bf{n}_1,\bf{n}_2$ we denote the outer normals on the intersection of $\interface$
 and the boundaries of $\dom{1}, \dom{2}$. We refer to  Figure \ref{DD-two-domain-illustration} for a sketch of the described situation.   The entire  domain $\dom{}$ is filled by a  porous medium which is assumed to be isotropic on each subdomain. 
We consider the dynamics of two immiscible, incompressible and viscous fluids, denoted as
a wetting one ($w$) and a non-wetting one ($nw$).
Considering a hybrid ansatz we suppose the   full two-phase model to be valid in domain $\dom{2}$, cf.\ \cite{Helmig97,Bear2007},  whereas 
we assume that on $\dom{1}$  the simplified 
Richards model, cf.\ \cite{Richards1931,Richardson1922}, is justified. 
A typical situation in which this occurs is the flow of water and air 
through a porous medium that is so permeable that the  air phase can be considered  to be ``infinitely'' 
mobile,  resulting in a constant pressure field equal to the atmospheric pressure. 
In view of the model hierarchy discussed e.g.\ in \cite{HenryHilhorstEymard2012},
the \Richards model can be viewed as the limit 
of a two-phase flow regime if the  ratio
of the nonwetting and the wetting viscosity tends to zero (and hence the mobility to infinity). 
With this interpretation, other situations than water and air are conceivable for a hybrid model ansatz.

Precisely, we consider  the following coupling of 
the \Richards equation with the 
two-phase flow  model in pressure-pressure formulation. 
\begin{problem}[\TPR problem]%
  \label{problem:TPR:continuous}
  For $l\in \{1,2\}$, let $\domT{l}:= \spaceTimeDom{l}$ and 
  $\interface_T :=  \spaceTimeInterface$. 
  The coupled \emph{\TPRfull} (\tpr)  problem consists of finding phase pressures $\pwln[1]{}$, $\pwln[2]{}$ and $\pnwln[2]{}$ solving
  \begin{align}
    \dtSR{1}%
      - \dv\biggl(\frac{k_{\mathsf{i},1}}{\mu_w} k_{w,1}\bigl( S_1 \bigr) \gradPwlnPlusGravity[1]{} \biggr)
      &=  f_{w,1}
      \label{eqn:LDD-TPR:Continuous:Richards:Strong}
    \end{align}
    in  $\domT{1}$ together with
    \begin{align}
    \dtS{2}%
      - \dv\biggl(\frac{k_{\mathsf{i},2}}{\mu_w} k_{w,2}\bigl( S_2 \bigr) \gradPwlnPlusGravity[2]{}\biggr)
      &=  f_{w,2}
    \label{eqn:LDD-TPR:Continuous:TwoPhase:Wetting:Strong},\\[1ex]
    \hspace{-2ex}- \dtS{2}%
      - \dv\biggl(\frac{k_{\mathsf{i},2}}{\mu_{\nw}} k_{\nw,2} \bigl( 1 - S_2\bigr)  \gradPnwlnPlusGravity[2]{} \biggr)
      &=  f_{\nw,2}
      \label{eqn:LDD-TPR:Continuous:TwoPhase:Nonwetting:Strong}
  \end{align}
  in $\domT{2}$.  Equations (\ref{eqn:LDD-TPR:Continuous:Richards:Strong})–(\ref{eqn:LDD-TPR:Continuous:TwoPhase:Nonwetting:Strong}) are coupled via
  \begin{align}
        \begin{aligned}
	\pwln[1]{} &= \pwln[2]{}, && \hspace{1.5mm}\flux[1]{w}\cdot\vt{n_1} = \flux[2]{w}\cdot\vt{n_1}, \\
	\pnwln[2]{} &= p_a,   && \flux[2]{\nw}\cdot\vt{n_2} = \flux[1]{\nw}\cdot\vt{n_2}
	\end{aligned}& && 
        \label{eqn:LDD-TPR:ContinuousCouplingConditions:strong}
  \end{align}
  \mbox{on} $\interface_T$. The problem is closed by suitable initial as well as boundary conditions.
\end{problem}
We  supplement  the  notations  used in Problem \ref{problem:TPR:continuous} including the specification of the  fluxes in \eqref{eqn:LDD-TPR:ContinuousCouplingConditions:strong}.
For subdomain index $l\in\{1,2\}$, 
            our primary unknowns are the wetting pressure $\pwln{}$ and the nonwetting pressure $\pnwln{}$ on $\domT{l}$, respectively.
{%
            The given constant atmospheric pressure is denoted by $p_a$ and 
            on $\domT{1}$, we have $\pnwln[1]{}=p_a$, by assumption.
            The functions
            $S_l = S_l(p_{\nw,l},p_{w,l}) \in [0,1]$ denote the wetting saturations and are assumed to be functions of the phase pressures via the capillary pressure saturation relationships $p_{c,l}(S_l) = \pnwln{}- \pwln{}$, see e.g.\ \cite{Helmig97}, i.e., it is assumed that the functions $p_{c,l}$ are invertible, cf.\ Assumption \ref{env:LDD-TPR:Assumptions:Data}. 
            Since we model two-phase flow, we use the assumption that on all subdomains $\dom{l}$ only the two phases are present, i.e., 
            the nonwetting saturations $S_l^{\nw}$ can be expressed by the relations $S_l^{\nw} = 1-S_l$.
            This is already built into the equation (\ref{eqn:LDD-TPR:Continuous:TwoPhase:Nonwetting:Strong}).\\
            The porosities $\porosity{l}\in (0,1)$ on each subdomain $\dom{l}$ are assumed to be constant and furthermore, 
            we denote by $\rho_\alpha > 0$ the density and by $\mu_\alpha > 0$ the viscosity of phase $\alpha\in \{w,\nw\}$. 
            For simplicity, we assume that the intrinsic permeabilities $\vt{K_l}$ are isotropic on every subdomain, 
            i.e., $\vt{K_l} = k_{\mathsf{i},l}\vt{E_d}$.  %
            Finally, for $\alpha \in \{w, \nw\}$, $k_{\alpha,l}$ denotes the given relative permeability, $f_{\alpha,l}$ a source term  and 
            $z_\alpha = \rho_\alpha g x_{d}$  is the gravitational force term ($g$ being the gravitational
            acceleration).\\
            The fluxes in \eqref{eqn:LDD-TPR:ContinuousCouplingConditions:strong} determine the mass flow coupling between the domains. For 
            $ (\alpha,l) \in \{ (w,1), (w,2),(\nw,2) \}
            $ they are given by 
            \begin{align}
                    \flux{\alpha} = - \frac{k_{\mathsf{i},l}}{\mu_{\alpha}}  \bm{\nabla} \big(k_{\alpha,l}( S_l) 
                    \big( p_{\alpha,l} +z_\alpha\big)
                 .
                \label{eqn:LDD-TPR:FluxContinuityContinuous}
            \end{align}
It remains do determine the flux   $\flux[1]{\nw}$.
  When coupling the \Richards model with the two-phase flow equations, it is not   clear which conditions
  should be imposed in (\ref{eqn:LDD-TPR:ContinuousCouplingConditions:strong}), 
  because the nonwetting phase is considered to be present, yet remains unmodelled.
  Since on $\dom{1}$ the nonwetting pressure is assumed constant, $\pnwln{} = p_a $, the part of the nonwetting Neumann flux containing the gradient of the pressure (in a two-phase flow model) would have to vanish. 
  However, this is not the case for the gravitational  part. 
  Thus, there are  two possible ways to account for the gravitational force of the nonwetting phase on $\dom{1}$ at the interface.
  In view of the fact, that the \Richards model is the mathematical limit of the two-phase model, cf.\ \cite{HenryHilhorstEymard2012}, one choice is
  \begin{align}
    \bm{F_{\nw,1}} = \frac{k_{\mathsf{i},1}}{\mu_{\nw}} k_{\nw,1}\bigl( 1 - S_1\bigl(p_a,\pwln[1]{}\bigr) \bigr)  \vt{\bm{\nabla}} z_{\nw}. %
                                            \label{eqn:LDD-TPR:ContinuousCouplingConditions:strong:alternative}
  \end{align}
  On the other hand, one could ignore  the effect entirely, i.e.\ %
 \begin{equation}
    \bm{F_{\nw,1}} = \mathbf{0}.                        \label{eqn:LDD-TPR:CouplingConditions:continuous:zero_flux}
  \end{equation}

  The two couplings (\ref{eqn:LDD-TPR:ContinuousCouplingConditions:strong:alternative}) 
  and (\ref{eqn:LDD-TPR:CouplingConditions:continuous:zero_flux}) are suggested in an adhoc manner.
  A  rigorous derivation of coupling conditions via e.g.\ homogenisation techniques is out of the scope of the present paper.  However,  we point out that the formulation of the \LDD scheme and its proof of convergence work for  both cases.

\begin{remark}[Extended coupling conditions]
    \label{remark:LDD-TPR:discussionCouplingConditions}
  The coupling conditions in \autoref{problem:TPR:continuous} 
    are the generic domain decomposition coupling conditions providing equivalence of the substructured problem to a monodomain formulation. 
  While natural in this sense, they exhibit certain limitations from a modelling perspective. 
  Indeed, since we 
  prescribe the continuity of the phase pressures, the capillary pressures  are  
  continuous as well. 
  However, in general, capillary trapping phenomena can occur for heterogeneous soils, where a phase might not enter into another soil layer due to a pressure barrier. This translates to a pressure jump over the interface.
  Non-matching capillary pressure 
  curves that in addition are extended to multivalued functions for vanishing (wetting or nonwetting) saturations, need to be considered in this case, cf.\ \cite{Cances2012}.
  This approach reflects pressure discontinuities over the interface
  by imposing the continuity of the capillary pressures together with the
  continuity of the wetting pressure in a generalised, multivalued sense.
  However, the analytical treatment (proof of existence of solutions) of this generalised formulation  consists in {approximating nonmatching capillary 
  pressure curves by a family of matching curves} with continuous  phase pressures.\\
  From the numerical point of view it is therefore  important to 
  investigate the applicability of the \LDD solver 
  to the case of continuous pressures not only as a first step, but 
  notably so as an approximation of the more realistic discontinuous case. We refer to \cite{AHMED2020113294} for a recent contribution in this direction.

\end{remark}

\subsection{Function spaces\label{sec:spaces}}
Before we  give  the weak formulation for 
\autoref{problem:TPR:continuous} we introduce  some notions for
function spaces  on  Lipschitz domains  and their boundaries, the latter being essential 
for the  analysis of the transmission conditions in the domain decomposition method. 
In this section, $\dom{}\subset \RR^d$ denotes a generic Lipschitz domain. In particular, all notations apply to all domains  $\dom{},\dom{l}$, $l=1,2$ introduced in the previous sections. 
  
  {\scshape Spaces on $\dom{}$.} 
    $C_0^\infty(\dom{})$ denotes  the space of smooth functions with compact support in
    $\dom{}$.
    $L^2(\dom{})$ is the space of   square-integrable  functions equipped with the scalar product  $\langle\cdot,\cdot\rangle$.
    For $s\in \RR$, the space $H^s(\RR^d)$ 
    denotes the standard Sobolev-Slobodeckij space with norm   $\norm{u}_{H^s(\dom{})} $. 
    We will need
    $H_0^1(\dom{}) = \completion{C_0^\infty(\dom{})}{H^1}$, and  for vector-valued functions $v: \dom{} \rightarrow \RR^d$, the space 
    \[
     H(\Div,\dom{}):= \Big\{v \in \left[L^2(\dom{})\right]^d\, \Bigr| \, \Div v \in L^2(\dom{}) \Big\},
    \]
    together with the norm $\norm{v}_{H(\Div,\dom{})}^2:= \norm{v}_{L^2}^2 + \norm{\Div v}_{L^2}^2$, $\Div v$ being understood via the integration by parts 
    formula.
  
  {\scshape Spaces on $\interface\subset\partial\dom{}$.} \quad
    The spaces $H^s(\interface)$ for $|s|\leq 1$ are  defined by understanding  that functions on $\partial \dom{}$   
    in local coordinates
    belong to $H^s(\RR^{d-1})$.
    When the Lipschitz surface $\partial \dom{}$  is divided into  two surfaces 
    $\interface_1$ and $\interface_2$, 
    $\partial \dom{} = \interface_1 \cup \partial \interface{}_j \cup \interface_2 $,
    with their common boundaries $\partial\interface_j$ of dimension $d-2$ in turn being Lipschitz, 
    the spaces $H^s(\interface_j)$ for $|s|\leq 1$ can be  introduced in the same way. 
    For a function $u \in H^{1/2}(\interface_j)$ the extension $\widetilde{\cdot}$ by 
    zero on $\partial \dom{}\setminus \interface_j$ does not imply $u \in H^{1/2}(\partial \dom{})$, see \cite[Theorem 3.4.4]{Agranovich2015} and the discussion thereafter. 
    In order to define Neumann traces 
    in a generalised sense via the Green formula on parts of the boundary, we need to define the subspace of those functions 
    in $H^{1/2}(\interface)$ for which the 
    extension by zero belongs to $H^{1/2}(\partial \dom{})$,
    that is
    \[
      \Tracespace := \left\{ u\in H^{1/2}(\interface)\, \left|\right. 
                           \, \widetilde{u} \in H^{1/2}(\partial \dom{})\right\}, \qquad 
      \norm{u}_{\Tracespace}:= \norm{\widetilde{u}}_{H^{1/2}(\partial \dom{})}.
    \]
    With the scalar product  inherited from $H^{1/2}(\partial \dom{})$  the space $\Tracespace$
    becomes a Hilbert space. 
   With these definitions, %
   the trace 
    operator $\tr{\interface} : H^1(\dom{}) \rightarrow H^{1/2}(\interface)$ can be defined %
    as extension of the restriction on smooth functions, acting as a bounded, surjective linear 
    operator on these spaces with bounded right inverse $\mathcal{R}_{\interface} : H^{1/2}(\interface) \rightarrow H^1(\dom{}) $, cf.
    \cite[Theorem A.2.3 and p. 132 ff]{Berninger2009}, \cite[Theorem 9.2.1, p. 118]{Agranovich2015} or \cite{Brezzi1991,MacLean2000}.
    To ease the notation, we will denote the trace by $\restrictTo{u}{\interface}$ instead of 
    $\tr{\interface} u$. 
  Moreover,  there  is a unique linear continuous operator $\gamma^{\nv}_{\partial\dom{}}: H(\Div,\dom{}) \rightarrow H^{-1/2}(\interface)$ such that $\gamma^{\nv} u = \restrictTo{(u\cdot{\nv})}{\interface}$ for 
  $u\in H(\Div,\dom{}) \cap \bigl[C(\overline{\dom{}})\bigr]^d$. It is in this generalised sense that 
  we will understand Neumann fluxes.

   {\scshape Dual spaces.} \label{notation.duality}
    Denoting the dual spaces $\mathcal{L}(H^s,\RR)$ equipped with the standard norm by ${H^s}'$, 
    functionals $F\in {H^s}'$   
    can  be identified via the Riesz representation theorem as an element of $v_F \in H^s$ itself, i.e. 
    $F(u) = \langle v_F, u\rangle_{H^s}$ for $u\in H^s$ and $\langle \cdot, \cdot\rangle_{H^s}$ 
    denoting the scalar product. 
    However, extending the form 
    $(u,v)_{\RR^d} := \langle u, v\rangle_{H^0(\RR^d)}$ to $H^{-s}(\RR^n)\times H^{s}(\RR^d)$ renders the spaces 
    $H^{-s}(\RR^d)$ and $H^{s}(\RR^d)$ mutually dual, providing an alternative representation of functionals on $H^s$, cf. \cite[p.9 ff]{Agranovich2015}.
    A similar duality holds for $H^s$-spaces on $\dom{}$ and $\interface \subset \partial \dom{}$, see \cite[Theorem 5.1.12, p. 61]{Agranovich2015}. For the case $s=\frac12$ which we need here, we have $\Tracespace' = H^{-1/2}(\interface)$ for 
    $\interface \subset \partial \dom{}$.
    We will use the symbol $\langle\cdot, \cdot\rangle_{\interface}$ for 
    the evaluation $F(\varphi)$ of a functional $F\in H^{-1/2}(\interface)$ with a function $\varphi \in \Tracespace$ also referred to as dual pairing.

\begin{remark}
  \label{rem_represent}  
 Note that  the considerations on duality from above show
that   each functional $u$ in $\Tracespace'$ has two representations.
  Namely, there is a function $\widehat{u}\in H^{-1/2}(\interface)$ and another function $\overline{u}\in \Tracespace$ such that 
  \begin{align} \bigl\langle u,\varphi \bigr\rangle_{\interface} =
        \bigl\langle \widehat{u},\varphi \bigr\rangle_{H^0(\interface)}
          \text{ and } 
         \bigl\langle u,\varphi \bigr\rangle_{\interface} 
        = \bigl\langle \overline{u},\varphi \bigr\rangle_{\Tracespace}.
        \label{LDD-RR:eqn:representations_of_functional}
    \end{align}
    hold. The choice of representation will be important in the formulation 
    of the domain decomposition  scheme below.
\end{remark}

\subsection{The LDD-TP–R solver for the \TPR problem}
In this section we introduce a time-discrete weak formulation of \autoref{problem:TPR:continuous} and formulate an \LDD solver for this setting.  Based on 
Section \ref{sec:spaces} we  define spaces associated to the subdomain partition. For $l\in \{1,2\}$, we define
  \[
    \Fs_l:= %
    \left\{ u \in H^1(\dom{l})\, \bigl|\bigr. \, \restrictTo{u}{\oB{l}}
    \equiv 0 \right\}  \text{ and } %
    \Fs   := \left\{ (u_1,u_2) \in \Fs_1\times \Fs_2 \, \bigl|\bigr. \, {u_1}_{|_\interface} \equiv {u_2}_{|_\interface}\right\},%
  \]
  where the norms in the spaces $\Fs_l$ are the standard $H^1(\dom{l})$-norms, and on $\Fs$ the norm 
  ${\norm{\cdot}}_{\Fs}^2 = \sum_{l=1,2}{\norm{\cdot}}_{\Fs_l}^2$ is used. 
  $\Fs_l'$ denotes again the dual space of $\Fs_l$ and is equipped with the usual norm for functionals
  ${\norm{F}}_{\Fs_l'} = \sup_{\varphi_l \in \Fs_l} \frac{\norm{F\varphi_l}_{\Fs_l}}{\norm{\varphi_l}_{\Fs_l}}$.

\begin{remark}  
    Henceforth, we assume that the atmospheric pressure $p_a$ vanishes. This can be done without loss of generality: let $\tilde{p}$ be a physical pressure and $\tilde{p}_a$ the atmospheric  pressure.  
    By introducing $p := \tilde{p} - \tilde{p}_a$, 
    the desired 
    normalisation $p_a = 0$ is achieved and  
    (\ref{eqn:LDD-TPR:Continuous:Richards:Strong}),
    (\ref{eqn:LDD-TPR:Continuous:TwoPhase:Wetting:Strong}),
    (\ref{eqn:LDD-TPR:Continuous:TwoPhase:Nonwetting:Strong}) stay the same, 
    since $Dp = D\tilde{p}$ for all derivatives and 
    $p_c^{-1}(\tilde{p}_{\nw} - \tilde{p}_w) = S = p_c^{-1}\bigl(\tilde{p}_{\nw} -\tilde{p}_a  -( \tilde{p}_w -\tilde{p}_a)\bigr) = p_c^{-1}(p_{\nw} - p_w)$.\\
    As a consequence the  nonwetting pressure 
    unknown  in the two-phase model domain $\Omega_2$ at a discrete time step will be in 
    $H^1_0(\Omega_2)$ and not in the  space  ${\mathcal V}_2$.
\end{remark}

As  the first step towards  the \LDD solver, we formulate a time discrete version of \autoref{problem:TPR:continuous}.  
For $N\in \NN$, the introduction of the time step size 
$\tau := \tfrac{T}{N}$ partitions the interval $[0,T]$ into the $N+1$ time steps $t^n:= n\cdot \tau$, $n = 0,\dots, N$. 

  The functions $\pwln[1]{n}:\dom{1} \to \RR$ and $\pwln[2]{n}, \pnwln[2]{n}:\dom{2}\to\RR$ denote the unknown time-discrete pressures at time step $t^n$. In addition, we set $S_l^n = S_l(\pnwln[l]{n},\pwln[l]{n})$ (anticipating $ \pnwln[1]{n} = p_a=0$) and  abbreviate
    \begin{equation}
    \begin{aligned}
        \kwln{n} &:= \frac{k_{\mathsf{i},l}}{\mu_w}   k_{w,l}(S^n_l),%
        && \knwln{n} :=   \frac{k_{\mathsf{i},l}}{\mu_{\nw}} %
        k_{\nw,l}(1-S^n_l).
    \end{aligned}\label{LDD-TPTP:notation1}
    \end{equation}
Consequently the fluxes at $t^n$ write as 
    \[
        \fluxb[n]{w} := - %
      k^n_{w,l}
       \gradPwlnPlusGravity{n}, \, \, l\in \{1,2\}\text{ and }
        {\bm{F^n_{\nw,2}}} := - k^n_{\nw,2} %
        \bm{\nabla}(p^n_{\nw,2} + z_{\nw}).
    \]
Depending on the choices in  (\ref{eqn:LDD-TPR:ContinuousCouplingConditions:strong:alternative}) 
  and (\ref{eqn:LDD-TPR:CouplingConditions:continuous:zero_flux})    the time-discrete flux ${\bm{F^n_{\nw,1}}}$   is defined in the same way.

With a backward Euler discretisation in time, the time-discrete coupled  \TPR problem in weak form then reads as follows.

\begin{problem}[Time-discrete \TPR  problem]
\label{env:LDD-TPR:problemFormulation:semiDiscrete:weak} For some $n\in \NN$, let 
$(\pwln[1]{n-1},\pwln[2]{n-1}) \in \Fs$ and $\pnwln[2]{n-1} \in H_0^1(\dom{2})$
Then, the time-discrete \TPR problem consists of finding
$\bigl((\pwln[1]{n},\pwln[2]{n}), \pnwln[2]{n}\bigr) \in \Fs\times H_0^1(\dom{2})$,
such that $\fluxn[l]{\alpha}\cdot \vt{n_l} \in \dualTs{}$ holds for $l=1,2$, $\alpha \in \{n, \nw\}$ and such that the equations
  \begin{align}
     &\spl \Sln[1]{n} - \Sln[1]{n-1},\testfunc[1]{w}\spr  %
      - \tau\spl \fluxn[1]{w},\bm{\nabla}\testfunc[1]{w} \spr
      + \tau \spl \fluxn[2]{w}\cdot \vt{n_1},\testfunc[1]{w}\spr_\interface
      =  \tau\spl  f_{w,1}^n ,\testfunc[1]{w} \spr, \label{LDD-TPR:limitsystemeRichards} \\[1.1ex]
     &\spl \Sln[2]{n} - \Sln[2]{n-1},\testfunc[2]{w}\spr  %
      - \tau\spl \fluxn[2]{w},\bm{\nabla}\testfunc[2]{w} \spr
      + \tau \spl \fluxn[1]{w}\cdot \vt{n_2},\testfunc[2]{w}\spr_\interface
      =  \tau\spl  f_{w,2}^n ,\testfunc[2]{w} \spr, \label{LDD-TPR:limitsystemeqn1w} \\[1.1ex] %
    -&\spl \Sln[2]{n} - \Sln[2]{n-1},\testfunc[2]{\nw}\spr  %
      - \tau\spl \fluxn[2]{\nw},\bm{\nabla}\testfunc[2]{\nw} \spr %
      + \tau \spl \fluxn[1]{\nw}\cdot \vt{n_2},\testfunc[2]{\nw}\spr_\interface
      =  \tau\spl  f_{\nw,2}^n ,\testfunc[2]{\nw} \spr \label{LDD-TPR:limitsystemeqn1nw}
  \end{align}
  are satisfied for all $(\testfunc[1]{w},\testfunc[2]{w}) \in \Fs_1\times\Fs_2$ and $\testfunc[2]{\nw} \in \Fs_2$.
\end{problem}%

\begin{remark} 
    \label{remark:LDD-TPR:semi-discreteRemark}
    \begin{enumerate}[label=\roman*)]
        \item
        In what follows, we will assume, that there is a unique solution $\bigl((\pwln[1]{n},\pwln[2]{n}), \pnwln[2]{n}\bigr) \in \Fs\times H_0^1(\dom{2})$ 
        for the  nonlinear time-discrete  \autoref{env:LDD-TPR:problemFormulation:semiDiscrete:weak}.
        We are not aware of any  results concerning well-posedness but we expect that standard methods for nonlinear parabolic equations can be applied.
        \item 
         Note that  traces of functions are  implicitly 
         taken in \eqref{LDD-TPR:limitsystemeRichards}-\eqref{LDD-TPR:limitsystemeqn1nw}.  They are needed for the 
         dual pairings of functionals on $\interface$ or likewise the scalar product of spaces on $\interface$, i.e. 
    $\langle u, \varphi\rangle_{\interface} = \langle \onGamma{u}, \onGamma{\varphi}\rangle_{\interface}$ if $u, \varphi \in \Fs_l$.\\
        For each solution $\bigl( (\pwln[1]{n}, \pwln[2]{n}), \pnwln[2]{n} \bigr)$ of 
        \autoref{env:LDD-TPR:problemFormulation:semiDiscrete:weak} 
        the coupling conditions (\ref{eqn:LDD-TPR:ContinuousCouplingConditions:strong}) are implicitly fulfilled in a weaker form at each time step $t_n$. 
        Namely, we have $\restrictTo{\pwln[1]{n}}{\interface} = \restrictTo{\pwln[2]{n}}{\interface}$ and 
        $\restrictTo{\pnwln[2]{n}}{\interface}=0$
        in the sense of traces by the definition of the spaces $\mathcal{V}$, and  $H_0^1(\dom{2})$.
        The continuity of the fluxes, $\fluxb[l]{\alpha}\cdot\bm{n_l} = \fluxn[3-l]{\alpha}\cdot\bm{n_l}$, $\alpha \in\{w,\nw\}$,
        is given as equality of functionals in $\Tracespace'$. 
        This is true regardless of the different choices for $\fluxn[1]{\nw}$.
    \end{enumerate}
\end{remark}

Next, based on the time-discrete  \TPR problem in weak formulation, we 
define the iterative  domain decomposition ansatz with iteration number $i\in \NN_0$. Extending the notation once more,
 the functions $\pwln[1]{n,i}:\dom{1} \to \RR$ and $\pwln[2]{n.i}, \pnwln[2]{n,i}:\dom{2}\to\RR$ denote the unknown $i$th pressure iterate at time step $t^n$.
 We set  again $S_l^{n,i} = S_l(\pnwln[l]{n,i},\pwln[l]{n,i})$ 
 using  $ \pnwln[1]{n,i} = p_a=0$. With the notations 
    \begin{equation}
    \begin{aligned}
        \kwln{n,i} &:= \frac{k_{\mathsf{i},l}}{\mu_w}   k_{w,l}(S^{n,i}_l),%
        && \knwln{n,i} :=   \frac{k_{\mathsf{i},l}}{\mu_{\nw}} %
        k_{\nw,l}(1-S^{n,i}_l),
    \end{aligned} %
    \end{equation}
the flux iterates  at $t^n$ are given by 
    \[
        \fluxb[n,i]{w} := - %
      k^{n,i-1}_{w,l}
       \gradPwlnPlusGravity{n,i}, \, \, l\in \{1,2\}\text{ and }
        \fluxnwi[2]{i} := - k^{n,i-1}_{\nw,2} %
        \bm{\nabla}(p^{n,i}_{\nw,2} + z_{\nw}).
    \]
The  flux iterate  ${\bm{F^{n,i}_{\nw,1}}}$   is defined  case by case in the same way.

\LDD schemes are designed to deal at the same time with two difficulties  of 
\autoref{env:LDD-TPR:problemFormulation:semiDiscrete:weak}. 
Firstly, 
each equation in \autoref{env:LDD-TPR:problemFormulation:semiDiscrete:weak} is doubly nonlinear, nonlinearities being present in the discretised time derivative, 
as well as in the fluxes. 
Secondly, the system of equations (\ref{LDD-TPR:limitsystemeRichards}) to (\ref{LDD-TPR:limitsystemeqn1nw}) is nonlinearly coupled
and the coupling conditions contain nonlinearities themselves.
The \LDD method tackles both of these problems by linearising and decoupling the equations in one single fixed point iteration.

We assume
that $\bigl( (\pwln[1]{n-1}, \pwln[2]{n-1}), \pnwln[2]{n-1} \bigr) \in \Fs \times H_0^1(\dom{2})$ is given %
 and set as initial iterates in the $n$th time step
\begin{equation}
        \pwlni{0}:= p_{w,l}^{n-1} \quad \mbox{on } \dom{l}, \, l=1,2 \text{, and }
        \pnwlni[2]{0}:= p_{\nw,2}^{n-1}, \quad \mbox{on } \dom{2},
   \label{eqn:LDD-TPR:initialPressures}
\end{equation}
Let the  numbers $\lambda_{\alpha} \in(0,\infty)$, $\alpha \in \{n,nw\}$ be given. They are bound to control the ratio between
Dirichlet and flux-type transmission conditions. Following Lions, \cite{Lions1988,Lions1990III}, we   introduce the Robin-type interface terms  
\begin{align} 
    \gali{0} &:= \fluxb[n-1]{\alpha} \cdot \vt{n_l} - \lambda_{\alpha} \onGamma{\paln[l]{n-1} }
    \label{eqn:LDD-TPR:gl0term:strong}
\end{align}
as functionals in $\Tracespace'$ for both phases and both subdomains.
  Since on $\dom{1}$ the nonwetting pressure is constant, $\pnwln[1]{} = p_a = 0 $, 
  we define either 
  \begin{align}
    \gnwli[1]{0} := 
    \fluxnMinusOne[1]{\nw}\cdot\vt{n_1} =
   k^{n-1}_{\nw,1}
   \vt{\bm{\nabla}} z_{\nw} \cdot \vt{n_1}
    \label{eqn:LDD-TPR:gl0withGravity}
  \end{align}
 if gravity effects are {included}, corresponding to the right hand side of  (\ref{eqn:LDD-TPR:ContinuousCouplingConditions:strong:alternative}), or
 \begin{align}
    \gnwli[1]{0} := \fluxnMinusOne[1]{\nw}\cdot\vt{n_1} = 0
    \label{eqn:LDD-TPR:gl0withoutGravity}
 \end{align}
 instead, in case gravity effects are {excluded}, as is expressed through 
 the right hand side of  (\ref{eqn:LDD-TPR:CouplingConditions:continuous:zero_flux}).
 
 \begin{note}[Pressure functionals]
 \label{note:LDD-TPR:pressureFunctionalsDefiniton}
 For $p\in \Fs_l$, we define the pressure functionals 
 $\langle p, \cdot \rangle_{\interface}\in \Tracespace'$ on the interface such as the ones appearing in (\ref{eqn:LDD-TPR:gl0term:strong}) according to
 \begin{align}
    \langle p, \cdot \rangle_{\interface}:= \langle p, \cdot \rangle_{1/2}. \label{LDD-RR:eqn:pressureFunctionalDef2}
 \end{align}
 In this way, the pressure traces $\onGamma{p}$ are representing the functionals $\langle p, \cdot \rangle_{\interface}\in \Tracespace'$ w.r.t the 
 $\Tracespace$-scalar product, i.e., $\onGamma{p} = \overline{p}$, cf.\ Remark \ref{rem_represent}. This will be important in the proof of Theorem \ref{LDD-TPR:convergence:theorem}.
  \end{note}
Now, the \ldd-TPR scheme  approximates the solution   to the time-discrete Problem 
\ref{env:LDD-TPR:problemFormulation:semiDiscrete:weak} at time $t_n$ by solving subsequently the following  problem (\ldd-TPR solver), together with the initial iterates given in (\ref{eqn:LDD-TPR:initialPressures}) and (\ref{eqn:LDD-TPR:gl0term:strong}).

\begin{problem}[\ldd-\TPR solver]%
\label{env:LDD-TPR:LDD-scheme:weak}
Let $L_{w,2},L_{\nw,2}  ,L_{w,1} > 0$ and some previously known iterates
$\palni{i-1} \in \Fs_l$,  $\gali{i-1}\in \Tracespace'$ be given for  $ i\in \NN, \, n \geq 1$.\\  
Find $\bigl(\pwlni[1]{i}, \pwlni[2]{i}, \pnwlni[2]{i}\bigr) \in \dFs\times \Fs_2$,
such that
\begin{align}
  L_{w,l}&\spl \pwlni{i},\testfunc{w} \spr
    - \tau\spl \fluxwi{i},\bm{\nabla}\testfunc{w} \spr
    + \tau \spl \lambda_w\pwlni{i}
    + \gwli{i} ,\testfunc{w}\spr_\interface  \nonumber \\
  &= L_{w,l}\spl \pwlni{i-1},\testfunc{w}\spr
      - \spl \Slni{i-1} - \Sln{n-1},\testfunc{w}\spr
    + \tau\spl  f_{w,l}^n ,\testfunc{w} \spr
      \label{eqn:LDD-TPR:LDD-scheme:wetting}
\end{align}
is fulfilled for $l\in \{1,2\}$ with
\begin{align}
  \spl \gwli{i},\testfunc{w}\spr_{\interface}  :=\spl -2\lambda_{w}
  \pwlni[3-l]{i-1}-   \gwli[3-l]{i-1},\testfunc{w}\spr_{\interface}, 	 \label{eqn:LDD-TPR:LDD-scheme:gliupdate:wetting:weak}
\end{align}
as well as 
\begin{align}
  L_{\nw,2}&\spl \pnwlni[2]{i},\testfunc[2]{\nw} \spr
    - \tau\spl \fluxgi[2]{i},\bm{\nabla}\testfunc[2]{\nw} \spr
    + \tau \spl \lambda_{\nw}\pnwlni[2]{i}
    + \gnwli[2]{i} ,\testfunc[2]{\nw}\spr_\interface \nonumber \\
  &= L_{\nw,2}\spl \pnwlni[2]{i-1},\testfunc[2]{\nw}\spr
      +  \spl \Slni[2]{i-1} - \Sln[2]{n-1},\testfunc[2]{\nw}\spr
    + \tau\spl  f_{\nw,2}^n ,\testfunc[2]{\nw} \spr \label{eqn:LDD-TPR:LDD-scheme:nonwetting}
\end{align}
with
\begin{align}
    \spl \gnwli[1]{i},\testfunc[2]{\nw}\spr_{\interface} &:=\spl -2\lambda_{\nw}
      \pnwlni[2]{i-1}-\gnwli[2]{i-1},\testfunc[2]{\nw}\spr_{\interface}, \label{eqn:LDD-TPR:LDD-scheme:gliupdate:nonwetting:weak:1} \\
      \spl \gnwli[2]{i},\testfunc[2]{\nw}\spr_{\interface} 
       &:=\spl -\gnwli[1]{i-1},\testfunc[2]{\nw}\spr_{\interface}, \label{eqn:LDD-TPR:LDD-scheme:gliupdate:nonwetting:weak:2} 
\end{align}
for all $\testfunc[1]{w}\in \Fs_1$, $\testfunc[2]{w}\in \Fs_2$ and $\testfunc[2]{\nw} \in \Fs_{2}$.
\end{problem}

 Notice that in a two-domain situation the index of the adjacent domain can be denoted by $3-l$,  for any given $l\in \{1,2\}$, since $3-l=2$ for $l=1$ and $3-l=1$ for $l=2$. This type of notation has been used in \autoref{env:LDD-TPR:LDD-scheme:weak}   and will be used  henceforth.
\begin{remark}
    \label{remark:LDD:TPR:gnwliPecularity}
  It may look peculiar to introduce an update for the term $\gnwli[1]{i}$ in (\ref{eqn:LDD-TPR:LDD-scheme:gliupdate:nonwetting:weak:1}) and
  (\ref{eqn:LDD-TPR:LDD-scheme:gliupdate:nonwetting:weak:2})
  as we do not have any equation for the nonwetting phase on $\dom{1}$. However,
  it is precisely this way of updating the $\gali{i}$ terms,  that liberates the 
  nonwetting pressure iterates $\pnwlni[2]{i}$ of the requirement to being elements of $H_0^1(\dom{2})$, i.e.\ to fulfil continuity to the atmospheric pressure \emph{in each iteration}. Instead, it 
  allows to merely require that $\pnwlni[2]{i}$ is an  element of $\Fs_2$. 
  This is less restrictive. In the present formulation, 
  the \LDD solver enforces $\onGamma{\pnwlni[2]{i}} \rightarrow 0$ 
 in the limit $i\rightarrow \infty$ all by itself.
 Moreover, it enables us to formulate a scheme that treats both model assumptions 
 (\ref{eqn:LDD-TPR:ContinuousCouplingConditions:strong:alternative}) and (\ref{eqn:LDD-TPR:CouplingConditions:continuous:zero_flux}) in a unified manner.
\end{remark}

The  assertions of \autoref{remark:LDD:TPR:gnwliPecularity} will 
be verified 
once \autoref{env:LDD-TPR:problemFormulation:semiDiscrete:weak} is reformulated such that the reformulation can be recognised as the formal limit system of the solver and the  convergence  of the \ldd-\TPR solver to this reformulation is proven. 
The reformulation of \autoref{env:LDD-TPR:problemFormulation:semiDiscrete:weak} is given in the next section.

\subsection{Consistency of the  \LDDTPR solver with the  time-discrete \TPR Problem \ref{env:LDD-TPR:problemFormulation:semiDiscrete:weak}}
\label{section:LDD-TPR:modeOfAction}
Recall that for a solution $\bigl( (\pwln[1]{n}, \pwln[2]{n}), \pnwln[2]{n} \bigr) \in \Fs\times H^1_0(\dom{2})$ of 
\autoref{env:LDD-TPR:problemFormulation:semiDiscrete:weak}
the nonwetting Neumann flux $\fluxn[1]{\nw}\cdot \vt{n_1}$ is defined by the right hand side of either 
(\ref{eqn:LDD-TPR:ContinuousCouplingConditions:strong:alternative})
or (\ref{eqn:LDD-TPR:CouplingConditions:continuous:zero_flux}), cf.\ also 
(\ref{eqn:LDD-TPR:gl0withoutGravity}) and (\ref{eqn:LDD-TPR:gl0withGravity}).
Thus, the functionals 
\begin{align}
  \gwli{} &:= -\lambda_{w}\restrictTo{\pwln{n}}{\interface} + \fluxn{w}\cdot \vt{n_l} \quad (l=1,2),  \label{LDD-TPR:eqn:gwln:wetting} \\
  \gnwli[1]{} &:= -\lambda_{\nw}\underbrace{\onGamma{\pnwln[1]{n}}}_{=p_a=0} + \fluxn[1]{\nw}\cdot \vt{n_1} \text{ and }   
  \gnwli[2]{} := -\lambda_{\nw}\underbrace{\restrictTo{\pnwln[2]{n}}{\interface}}_{=0} + \fluxn[2]{\nw}\cdot \vt{n_2}, %
  \label{LDD-TPR:eqn:gwln:nonwetting2}
\end{align}
in $\Tracespace'$ 
fulfil the relations
\begin{align*}
  \gwli{} %
      &= -2\lambda_{w}\restrictTo{\pwln[3-l]{n}}{\interface} + \lambda_{w}\restrictTo{\pwln[3-l]{n}}{\interface} - \fluxn[3-l]{w}\cdot \vt{n_{3-l}}%
      = -2\lambda_{w}\restrictTo{\pwln[3-l]{n}}{\interface} - \gwli[3-l]{} \quad  (l=1,2),\\
  \gnwli[1]{} &= -2\lambda_{\nw}\restrictTo{\pnwln[2]{n}}{\interface} - \gnwli[2]{} \, \text{ and }
  \gnwli[2]{} = -\gnwli[1]{}. 
\end{align*}
Note that  $\restrictTo{\pwln[1]{n}}{\interface} = \restrictTo{\pwln[2]{n}}{\interface}$, 
$\restrictTo{\pnwln[2]{n}}{\interface}=0$
and  $\flux[1]{\alpha}\cdot\bm{n_1} = -\flux[2]{\alpha}\cdot\bm{n_2}$ for $\alpha \in \{w,\nw\}$. 
\autoref{env:LDD-TPR:problemFormulation:semiDiscrete:weak} can therefore 
be written as
\begin{align}
   \spl \Sln{n} - \Sln{n-1},\testfunc{w}\spr  %
    &- \tau\spl \fluxb[n]{w},\bm{\nabla}\testfunc{w} \spr %
    + \tau \spl \lambda_w \pwln{n}
    + \gwl,\testfunc{w}\spr_\interface
    =  \tau\spl  f_{w,l}^n ,\testfunc{w} \spr \tag{\ref{eqn:LDD-TPR:LDD-scheme:wetting}'} 
    \label{eqn:LDD-TPR:limitSystem:wetting} 
\end{align}
with 
\begin{align}
  \spl \gwl,\testfunc{w}\spr_{\interface}  =\spl -2\lambda_{w}
  \pwln[3-l]{n} -   \gwl[3-l]{},\testfunc{w}\spr_{\interface}, 
  \tag{\ref{eqn:LDD-TPR:LDD-scheme:gliupdate:wetting:weak}'}\label{eqn:LDD-TPR:limitSystem:glupdate}
\end{align}
for $l\in \{1, 2\}$ as well as 
\begin{align}
  -\spl \Sln[2]{n} - \Sln[2]{n-1},\testfunc[2]{\nw}\spr  %
    &- \tau\spl \fluxn[2]{\nw},\bm{\nabla}\testfunc[2]{\nw} \spr %
    + \tau \spl \lambda_{\nw}\pnwln[2]{n}
    + \gnwl[2]{} ,\testfunc[2]{\nw}\spr_\interface
    =  \tau\spl  f_{\nw,2}^n ,\testfunc[2]{\nw} \spr 
    \tag{\ref{eqn:LDD-TPR:LDD-scheme:nonwetting}'}\label{eqn:LDD-TPR:limitSystem:nonwetting} 
    \end{align}
    together with
    \begin{align}
    \spl \gnwl[1]{},\testfunc[2]{\nw}\spr_{\interface} &=\spl -2\lambda_{\nw}
      \pnwln[2]{n}-\gnwl[2]{},\testfunc[2]{\nw}\spr_{\interface}, 
      \tag{\ref{eqn:LDD-TPR:LDD-scheme:gliupdate:nonwetting:weak:1}'} \label{eqn:LDD-TPR:limitSystem:gli:nonwetting:1}\\
      \spl \gnwl[2]{},\testfunc[2]{\nw}\spr_{\interface} 
       &=\spl -\gnwli[1]{},\testfunc[2]{\nw}\spr_{\interface}. 
       \tag{\ref{eqn:LDD-TPR:LDD-scheme:gliupdate:nonwetting:weak:2}'} \label{eqn:LDD-TPR:limitSystem:gli:nonwetting:2}
\end{align}
We know $\restrictTo{\pnwln[2]{n}}{\interface} = 0$ since $\pnwln[2]{n} \in H^1_0(\dom{2})$ and thus the pressure functionals in (\ref{eqn:LDD-TPR:limitSystem:nonwetting})–(\ref{eqn:LDD-TPR:limitSystem:gli:nonwetting:1}) actually disappear. 
They are written out here to emphasise the structure. 

Conversely, any tuple of functions $\bigl(\pwln[1]{n}, \pwln[2]{n}, \pnwln[2]{n}\bigr) \in \dFs\times\Fs_2$ 
fulfilling (\ref{eqn:LDD-TPR:limitSystem:wetting}) – (\ref{eqn:LDD-TPR:limitSystem:gli:nonwetting:2}) 
together with  
\begin{align}
    \gnwli[1]{} =  \fluxn[1]{\nw}\cdot \vt{n_1},
    \label{eqn:LDD-TPR:gnwlFlux}
\end{align}
where $\fluxn[1]{\nw}\cdot \vt{n_1}$ is defined by the right hand side of either 
(\ref{eqn:LDD-TPR:ContinuousCouplingConditions:strong:alternative})
or (\ref{eqn:LDD-TPR:CouplingConditions:continuous:zero_flux}),
is a solution of \autoref{env:LDD-TPR:problemFormulation:semiDiscrete:weak}.

The  argument supporting this claim for the wetting phase has been given in the proof of \cite[Lemma 2]{SeusMitra2018} or \cite[Lemma 2.3.12]{SeusThesis2021}
and it carries over to the 
situation here.

Regarding the nonwetting phase on $\dom{2}$, notice that $\restrictTo{\pnwln[2]{n}}{\interface} = 0$ is contained in
(\ref{eqn:LDD-TPR:limitSystem:gli:nonwetting:2}) and (\ref{eqn:LDD-TPR:limitSystem:gli:nonwetting:1}) as 
$
    \langle\pnwln[2]{n}, \testfunc[2]{\nw} \rangle_{\interface} = 0
$
follows for all $\testfunc[2]{\nw} \in \Fs_2$ by plugging (\ref{eqn:LDD-TPR:limitSystem:gli:nonwetting:2}) into (\ref{eqn:LDD-TPR:limitSystem:gli:nonwetting:1}).
By our definition of the pressure functionals, cf.\ \autoref{note:LDD-TPR:pressureFunctionalsDefiniton}, as well as by virtue of the surjectivity of the trace operator, this means 
$\langle\restrictTo{\pnwln[2]{n}}{\interface}, \eta \rangle_{1/2}
=\langle\pnwln[2]{n}, \eta \rangle_{\interface}=0$ for all $\eta \in \Tracespace$ and thus $\restrictTo{\pnwln[2]{n}}{\interface}= 0$. 

Using this and integrating (\ref{eqn:LDD-TPR:limitSystem:nonwetting}) by parts yields 
$\gnwli[2]{} =  \fluxn[2]{\nw}\cdot \vt{n_2}$. 
Since 
$\gnwli[2]{} = -\gnwli[1]{} =  \fluxn[1]{\nw}\cdot \vt{n_2}$ by 
(\ref{eqn:LDD-TPR:limitSystem:gli:nonwetting:2}) and (\ref{eqn:LDD-TPR:gnwlFlux}), the continuity of the fluxes follows. 

Consequently, we have proven the following 

\begin{lemma}[Limit of the \LDDTPR solver]
\label{lemma:LDD-TPR:iterationinterpretationlemma}
  Let $n \in \NN, \, n \ge 1$, be fixed, the tuple
  $\bigl( (\pwln[1]{n-1}, \pwln[2]{n-1}), \pnwln[2]{n-1} \bigr) \in \Fs \times H_0^1(\dom{2})$ be given
  and assume that 
   functions $\bigl(\pwln[1]{n}$, $\pwln[2]{n}$, $\pnwln[2]{n}\bigr) \in \dFs\times\Fs_2$ 
  and $\gali{}\in \Tracespace'$ exist for $\alpha\in\{n, \nw\}$ and $l\in\{1,2\}$, such that $\gnwli[1]{}$ is given by \eqref{eqn:LDD-TPR:gnwlFlux}, %
  and such that these functions fulfil the system of equations
  \eqref{eqn:LDD-TPR:limitSystem:wetting}- \eqref{eqn:LDD-TPR:limitSystem:gli:nonwetting:2}
for all $\testfunc[1]{w} \in \Fs_1$ and $\testfunc[2]{w}, \testfunc[2]{\nw} \in \Fs_2$.\\
Then, the interface conditions
\begin{align}
\onGamma{ \pwln[1]{n} } 		&= \onGamma{ \pwln[2]{n} } \quad \mbox{ and} \quad \onGamma{ \pnwln[2]{n}}  = 0, \\
\fluxn{\alpha} \cdot \vt{n_l}%
    &= \fluxn[3-l]{\alpha} \cdot \vt{n_l}, \quad \alpha\in \{w, \nw\} \,\, (l\in\{1,2\}) 
\label{lemma:LDD-TPR:FluxContinuity}
\end{align}
are satisfied in $ \Tracespace'$  and $\bigl((\pwln[1]{n}, \pwln[2]{n}), \pnwln[2]{n}\bigr)$ solves  \autoref{env:LDD-TPR:problemFormulation:semiDiscrete:weak}.
Moreover,
\begin{align}
    \gwli{} &= -\lambda_{w}\restrictTo{\pwln{n}}{\interface} + \fluxn{w}\cdot \vt{n_l} \quad  (l=1,2) \text{ and }  {\label{LDD-TPR:eqn:gwln:wetting:new}} \\
    \gnwli[2]{} &= -\lambda_{\nw}%
    {\restrictTo{\pnwln[2]{n}}{\interface}} %
    + \fluxn[2]{\nw}\cdot \vt{n_2}  
    {\label{LDD-TPR:eqn:gwln:nonwetting2:new}}
\end{align}
in $\Tracespace'$.

Conversely, if $\bigl((\pwln[1]{n}, \pwln[2]{n}), \pnwln[2]{n}\bigr) \in \Fs\times H^1_0(\dom{2})$ is a solution of \autoref{env:LDD-TPR:problemFormulation:semiDiscrete:weak} and 
$ \gwli{}$ is defined according to (\ref{LDD-TPR:eqn:gwln:wetting:new}),  $\gnwli[2]{}$ 
according to (\ref{LDD-TPR:eqn:gwln:nonwetting2:new}) and $\gnwli[1]{}$ given by (\ref{eqn:LDD-TPR:gnwlFlux}), then $\bigl((\pwln[1]{n}, \pwln[2]{n}), \pnwln[2]{n}\bigr)$ and $\gali{}$ solve the system
(\ref{eqn:LDD-TPR:limitSystem:wetting}) – (\ref{eqn:LDD-TPR:limitSystem:gli:nonwetting:2}).
\end{lemma}

\begin{remark}
\autoref{LDD-TPR:convergence:theorem} below shows that the family $\bigl\{\bigl(\pwlni[1]{i}, \pwlni[2]{i}, \pnwlni[2]{i}\bigr)\bigr\}_{i\in\NN}$ of subsequent solutions 
to \autoref{env:LDD-TPR:LDD-scheme:weak},
  together with the 
iterates $\bigl\{\gali{i}\bigr\}_{i\in\NN}$ converge to a solution of 
(\ref{eqn:LDD-TPR:limitSystem:wetting})–(\ref{eqn:LDD-TPR:limitSystem:gli:nonwetting:2}). By the just proven lemma, this means 
solving \autoref{env:LDD-TPR:problemFormulation:semiDiscrete:weak}.
Thus, it is  justified to refer to equations (\ref{eqn:LDD-TPR:limitSystem:wetting})–(\ref{eqn:LDD-TPR:limitSystem:gli:nonwetting:2})
as the limit system to \autoref{env:LDD-TPR:LDD-scheme:weak}. 
\end{remark}

%
\section{Convergence of the \LDDTPR solver}
\label{LDD-TPR:section:ConvergenceOfScheme}
In this  core section we analyse the convergence of the \LDDTPR solver.  Before doing so we state the general assumptions needed (see also the setting in \cite{Pop2004,Radu2017}).
\begin{assumption} 
\label{env:LDD-TPR:Assumptions:Data} Let $l\in \{1,2\}$.%
 \begin{enumerate}[itemsep=0.5ex,label=\alph*)]
   \item  The intrinsic permeabilities $k_{\mathsf{i},l}$ belong to $ L^\infty(\dom{},\RR_{+})\cap C^{0,1}(\dom{},\RR_{+})$.
  \item The relative permeabilities of the wetting phases $k_{w,l}:[0,1] \rightarrow [0,1]$ are strictly mo\-notonically \emph{increasing} and
  Lipschitz continuous functions with Lipschitz constants $L_{k_{w,l}}$.
  The relative permeabilities of the nonwetting phases $k_{\nw,l}:[0,1] \rightarrow [0,1]$ on both domains are strictly mo\-notonically \emph{decreasing} (as functions of the wetting saturation) and
  Lipschitz continuous functions with Lipschitz constants $L_{k_{\nw,l}}$.
  \item \label{LDD-TPR:env:assumptions:coercivity} There are numbers   $m_1,m_2 >0$  such that we have  $ \frac{k_{\mathsf{i},l} k_{w,1} }{\mu_w} \geq  m_1 $ and 
    \begin{align}
     m_2=\min\left\{\min_{s\in[0,1]}\frac{k_{\mathsf{i},2}}{\mu_w} k_{w,2}(s), \min_{s\in[0,1]}\frac{k_{\mathsf{i},2}}{\mu_{\nw}} k_{\nw,2}(s)\right\}. 
    \end{align}
  \item
  The capillary pressure saturation relationships $p_c^l(S_l) := \pnwln{} - \pwln{} $ are monotonically decreasing functions and therefore the saturations,  $S_l\bigl(p_c^l\bigr) = S_l\bigl(\pnwln{} - \pwln{}\bigr)$ are also monotonically decreasing as functions of $p_c^l$. Moreover, they are assumed to be Lipschitz continuous with Lipschitz constants $L_{S_l}$.
 \end{enumerate}
\end{assumption}

  Assumption \ref{LDD-TPR:env:assumptions:coercivity} is required
  to ensure the existence of a solution in each \LDDTPR solver step, and for the  convergence proof  of the \LDDTPR solver as well.   It excludes degeneracy and implicitly makes sure  that both phases are present on both sides of the interface avoiding  trapping effects.

Before stating our main result, we note the following lemma. In view of Assumption \ref{env:LDD-TPR:Assumptions:Data} it is a direct 
  consequence of the Lax-Milgram theorem  and guarantees that solving 
\autoref{env:LDD-TPR:LDD-scheme:weak} is always possible. %

\begin{lemma}
    \label{lemma:LDD-TPR:LDDTPRSolverHasSolution}
    Let Assumptions \ref{env:LDD-TPR:Assumptions:Data} hold true.  Given $f_{\alpha,l}^n \in \Fs_l'$, \autoref{env:LDD-TPR:LDD-scheme:weak} has a unique solution $\bigl(\pwlni[1]{i}, \pwlni[2]{i}, \pnwlni[2]{i}\bigr) \in \Fs_1\times\Fs_2\times\Fs_2$.
\end{lemma}

Considering a family of subsequent solutions to \autoref{env:LDD-TPR:LDD-scheme:weak}, we can now prove the following convergence result for the \LDDTPR solver.
\begin{theorem}[Convergence of the \LDDTPR solver]
\label{LDD-TPR:convergence:theorem}  Let Assumption \ref{env:LDD-TPR:Assumptions:Data} hold true and suppose that
   there exists a pair $\bigl((\pwln[1]{n},\pwln[2]{n}), \pnwln[2]{n}\bigr)\in \Fs\times  
  H^1_0(\dom{2}) $ that uniquely solves \autoref{env:LDD-TPR:problemFormulation:semiDiscrete:weak} 
   satisfying  for some $M>0$ the bound $\sup_{l,\alpha}\norm{\gradPalnPlusGravity{n} }_{L^\infty} \leq M $. %
  For $l\in \{1,2\}, \, \alpha \in \{w,\nw\}$ let $\lambda_\alpha > 0$
  and $L_{\alpha,l} > 0 $ be such that we have 
   \begin{align}
	\frac{1}{L_{S_2}\porosity{2}}-\sum_{\alpha}\frac{1}{2L_{\alpha,2}} > 0    \quad \mbox{ and } \quad
	\frac{1}{L_{S_1}\porosity{1}} - \frac{1}{2L_{w,1}} >0 \label{eqn:LDDrestriction.on.L}.
   \end{align}
  For arbitrary initial pressures $\pwlni{0} := \nu_{w,l} \in \Fs_l$, $l=1,2$, and $\pnwlni[2]{0} := \nu_{\nw,2} \in H^1_0(\dom{2})$, let $\bigl\{\bigl((\pwlni[1]{i},\pwlni[2]{i}), \pnwlni[2]{i} \bigr) \bigr\}_{i \in \NN_0} \in \bigl (\Fs_1\times\Fs_2 \times \Fs_2\bigr)^{\NN_0}$
  be a sequence of 
  solutions to \autoref{env:LDD-TPR:LDD-scheme:weak}, $\bigl\{\gali{i}\bigr\}_{i \in \NN_0}$ 
  be defined by
  (\ref{eqn:LDD-TPR:LDD-scheme:gliupdate:wetting:weak}) and (\ref{eqn:LDD-TPR:LDD-scheme:gliupdate:nonwetting:weak:1})–(\ref{eqn:LDD-TPR:LDD-scheme:gliupdate:nonwetting:weak:2}) and $\gal$ by 
  (\ref{LDD-TPR:eqn:gwln:wetting})–(\ref{LDD-TPR:eqn:gwln:nonwetting2}) for $\alpha \in\{w, \nw\}$ and $l\in\{1,2\}$.
  Assume, that the time step size $\tau$ has been chosen  to satisfy the conditions
  \begin{equation}
      C_2 := \frac{1}{L_{S_2}\porosity{2}} - \sum_{\alpha}\frac{1}{2L_{\alpha,2}} - \tau\sum_{\alpha} \frac{L_{k_{\alpha,2}}^2 M^2}{2m_2\porosity{2}^2}   > 0, \qquad 
      C_1 := \frac{1}{L_{S_1}\porosity{1}} - \frac{1}{2L_{w,1}} - \tau\frac{L_{k_{w,1}}^2 M^2}{2m_1\porosity{1}^2}   > 0.
      \label{eqn:timesteprestriction.2}
    \end{equation}
  Then,
  \begin{align}
    \begin{aligned}
        \pwlni{i} &\rightarrow \pwln{n} \quad \mbox{in } \Fs_l \nonumber\\ 
        \pnwlni[2]{i} &\rightarrow \pnwln[2]{n} \quad \mbox{in } \Fs_2
  \end{aligned},
  \quad \mbox{ and }
  \quad
  \begin{aligned}
    \gwli{i} &\rightarrow \gwl \quad \mbox{in } \Fs_l' \nonumber \\
    \gnwli{i} &\rightarrow \gnwl \quad \mbox{in } \Fs_l'
  \end{aligned},
  \end{align}
  for $l=1,2$ as $i\rightarrow \infty$. 
  Notably, $\onGamma{\pnwlni[2]{i}} \rightarrow 0$ in $\Tracespace$ as $i\rightarrow \infty$.
\end{theorem}
\begin{remark}[Explicit time step restriction]
 The implicit restrictions (\ref{eqn:timesteprestriction.2}) on the time step size translate to the explicit form 
    \begin{align}
    \tau < \min\left\{ \biggl(\frac{1}{L_{S_1}\porosity{1}} - \frac{1}{2L_{w,1}} \biggr)\frac{2m_1\porosity{1}^2}{L_{k_{w,1}}^2 M^2}, \,\,
      \frac{\frac{1}{L_{S_2}\porosity{2}} - \sum_{\alpha}\frac{1}{2L_{\alpha,2}}}{\sum_{\alpha}\bigr[ (L_{k_{\alpha,2}} M)^2/(2m_2\porosity{2}^2)\bigr]} 
      \label{eqn:timesteprestriction}
      \right\}.
  \end{align}
\end{remark}

\textit{Proof (of \autoref{LDD-TPR:convergence:theorem}):}
Note that for the proof, we will actually denote by $L_{k_{\alpha,l}}$ the Lipschitz constant of the function $\frac{\norm{k_{\mathsf{i},l}}_{\infty} k_{\alpha,l} }{\mu_\alpha}$  
by slight abuse of notation.
Define the iteration errors $\epwli{i}	:= \pwln{n}- \pwlni{i}$ and $\epnwli[2]{i}	:= \pnwln[2]{n}- \pnwlni[2]{i}$ as well as $\egali{i}	:= \gal - \gali{i}$ for $l=1,2$ and $\alpha \in \{w,\nw\}$. Add
$L_{w,l}\tspl \pwln{n},\testfunc{w}\tspr - L_{w,l}\tspl \pwln{n},\testfunc{w}\tspr$ to 
(\ref{eqn:LDD-TPR:limitSystem:wetting}) and  respectively $L_{\nw,2}\tspl
\pnwln[2]{n},\testfunc[2]{\nw}\tspr - L_{\nw,2}\tspl
\pnwln[2]{n},\testfunc[2]{\nw}\tspr$ to
(\ref{eqn:LDD-TPR:limitSystem:nonwetting})
 and subtract the corresponding equations (\ref{eqn:LDD-TPR:LDD-scheme:wetting}) as well as 
 (\ref{eqn:LDD-TPR:LDD-scheme:nonwetting}) to get
\begin{align} 
L_{\alpha,l}&\spl\epali{i} - \epali{i-1},\testfunc{\alpha} \spr %
    + \tau \spl \lambda_\alpha \epali{i}
    + \egali{i},\testfunc{\alpha} \spr_\Gamma 		%
    \nonumber\\[1ex]%
  &+ \tau \bspl-\fluxb[n]{\alpha}%
 		 {\color{\addnullcolor}-\kalni{i-1}\gradPalnPlusGravity{n} 
 		 +\kalni{i-1}\gradPalnPlusGravity{n}}%
 		 +\fluxai{i},\bm{\nabla}\testfunc{\alpha}\bspr 
 	   \label{LDD-TPR:mainproofstep1a} \\[1ex]
  &\quad= \underbrace{(-1)^{\krond{w}} \spl \Sln{n} -\Sln{n-1}, \testfunc{\alpha} \spr %
		-(-1)^{\krond{w}} \spl \Slni{i-1} - \Sln{n-1}, \testfunc{\alpha} \spr}_{%
		(-1)^{\krond{w}} \spl \Sln{n}- \Slni{i-1}, \testfunc{\alpha} \spr},
 		\nonumber 
\end{align}
where (\ref{LDD-TPR:mainproofstep1a}) is meaningful for the index combinations $(\alpha, l) \in \{(w,1), (w,2), (\nw,2)\}$.
Note the use of the Kronecker delta $\krond{w}$ to account for the minus sign of the time discretisation for the nonwetting phase.%

Inserting for all admissible index combinations  $\testfunc{\alpha} :=
\epali{i}$ in (\ref{LDD-TPR:mainproofstep1a})
all the while making use of the identity 
\begin{align*}
 L_{\alpha,l}\bspl\epali{i}-\epali{i-1},\epali{i}\bspr = \frac{L_{\alpha,l}}{2}\Bigl( \bigl\|\epali{i}\bigr\|^2 - 
\bigl\|\epali{i-1}\bigr\|^2 
+ 
\bigl\|\epali{i}-\epali{i-1}\bigr\|^2\Bigr), %
\end{align*}
leads to 
\begin{align}
  \frac{L_{\alpha,l}}{2}&\Bigl( \bigl\|\epali{i}\bigr\|^2 %
			  - \bigl\|\epali{i-1}\bigr\|^2 
			  + \bigl\|\epali{i}-\epali{i-1}\bigr\|^2\Bigr) 
    + \tau \lambda_\alpha \spl \epali{i},\epali{i} \spr_\Gamma
    = \spl \Sln{n}- \Slni{i-1}, (-1)^{\krond{w}}\epali{i} \spr%
    \nonumber \\[1ex]
  &- \tau  \spl \egali{i},\epali{i} \spr_\Gamma%
      -\tau \bspl\Bigl(\kaln{n} -\kalni{i-1}\!\Bigr)\gradPalnPlusGravity{n},\bm{\nabla}\epali{i}\bspr
    -\tau \bspl\kalni{i-1}\bm{\nabla} \epali{i},\bm{\nabla}\epali{i}\bspr.
  \label{LDD-TPR:mainproofstep2a}
\end{align}
Summing over phases $\alpha = w,\nw$ in (\ref{LDD-TPR:mainproofstep2a}) for $l=2$  and adding the term
$
    \spl \Sln[2]{n}- \Slni[2]{i-1}, \epwli[2]{i-1} -\epnwli[2]{i-1} \spr 
    $ 
on both sides of the equation, yields
\begin{align}
  \sum_{\alpha}\frac{L_{\alpha,2}}{2}&\Bigl( \bigl\|\epali[2]{i}\bigr\|^2 %
        - \bigl\|\epali[2]{i-1}\bigr\|^2 
        + \bigl\|\epali[2]{i}-\epali[2]{i-1}\bigr\|^2\Bigr) 
    {\color{\addnullcolor}+  \underbrace{\spl \Sln[2]{n}- \Slni[2]{i-1}, \epwli[2]{i-1} -\epnwli[2]{i-1} \spr}_{T\!P_1}} \nonumber\\
  &=\underbrace{\spl \Sln[2]{n}- \Slni[2]{i-1}, {\color{\addnullcolor} \epwli[2]{i-1}} - \epwli[2]{i} - {\color{\addnullcolor} \bigl(\epnwli[2]{i-1}} - \epnwli[2]{i}{\color{\addnullcolor}\bigr)} \spr}_{=: T\!P_2} - \tau \underbrace{\sum_{\alpha}\spl \lambda_{\alpha}\epali[2]{i} + \egali[2]{i},\epali[2]{i} \spr_\Gamma}_{=:I_{T\!P}} \nonumber \\  
  &\quad
  -\underbrace{\tau \sum_{\alpha} \bspl\kalni[2]{i-1}\bm{\nabla} \epali[2]{i},\bm{\nabla}\epali[2]{i}\bspr}_{=:T\!P_4}  
  -\underbrace{\tau \sum_{\alpha} \bspl\Bigl(\kaln[2]{n} -\kalni[2]{i-1}\!\Bigr)\gradPalnPlusGravity[2]{n},\bm{\nabla}\epali[2]{i}\bspr.}_{=: T\!P_3}
  \label{eqn:mainproofstep_summing.TP}
\end{align}
Similarly, adding 
$
    \spl \Sln[1]{n}- \Slni[1]{i-1}, \epwli[1]{i-1} \spr 
    $ 
to both sides of (\ref{LDD-TPR:mainproofstep2a}) for $l=1$, one gets
\begin{align}
  \frac{L_{w,1}}{2}&\Bigl( \bigl\|\epwli[1]{i}\bigr\|^2 %
        - \bigl\|\epwli[1]{i-1}\bigr\|^2 
        + \bigl\|\epwli[1]{i}-\epwli[1]{i-1}\bigr\|^2\Bigr) 
    {\color{\addnullcolor}+  \underbrace{\spl \Sln[1]{n}- \Slni[1]{i-1}, \epwli[1]{i-1} \spr}_{R_1}} \nonumber\\
  &= \underbrace{\spl \Sln[1]{n}- \Slni[1]{i-1}, {\color{\addnullcolor} \epwli[1]{i-1}} - \epwli[1]{i} \spr}_{=: R_2} %
  - \tau \underbrace{\spl \lambda_{w}\epwli[1]{i} + \egwli[1]{i},\epwli[1]{i} \spr_\Gamma}_{=:I_R} \nonumber \\  
  &\quad-\underbrace{\tau \bspl\Bigl(\kwln[1]{n} -\kwlni[1]{i-1}\!\Bigr)\gradPwlnPlusGravity[1]{n},\bm{\nabla}\epwli[1]{i}\bspr}_{=: R_3}
    -\underbrace{\tau \bspl\kwlni[1]{i-1}\bm{\nabla} \epwli[1]{i},\bm{\nabla}\epwli[1]{i}\bspr.}_{=:R_4} \label{eqn:mainproofstep_summing.R}
\end{align}
We proceed to estimate the assigned terms  $T\!P_1$--$T\!P_4$ and $R_1$--$R_4$ from (\ref{eqn:mainproofstep_summing.TP}) and (\ref{eqn:mainproofstep_summing.R}). 

\paragraph{$T\!P_1$, $R_1$} 
Recall that  $S_l\bigl(\pwln{},\pnwln{}\bigr) = (p_c^l)^{-1}\bigl(\pnwln{} - \pwln{}\bigr)$ for both $l$ and that ${p_c^l}' < 0 $ so that we actually have  the dependence $S_l\bigl(\pwln{},\pnwln{}\bigr) = S_l\bigl(\pnwln{} - \pwln{}\bigr)$ where $S_l$ 
as a function of $p_c^l$ is monotonically decreasing. 
Even though $\pnwln[1]{} = p_a = 0$ and there is no equation for $\pnwln[1]{}$, we estimate both $T\!P_1$ and $R_1$ with the same reasoning by setting $\pnwln[1]{n} = 0$ and $\pnwlni[1]{i} = 0$ for all $i\in \NN$.   Thereby, we have \emph{for $l=1,2$}
\begin{align}
 &\Bigl| \porosity{l}S_l\bigl(\pnwln{n} - \pwln{n}\bigr) - \porosity{l}S_l\bigl(\pnwlni{i-1} - \pwlni{i-1}\bigr) \Bigr|^2 \nonumber \\[1ex]
 &\quad\leq L_{S_l}\porosity{l} \Bigl| \porosity{l}S_l\bigl(\pnwln{n} - \pwln{n}\bigr) - \porosity{l}S_l\bigl(\pnwlni{i-1} - \pwlni{i-1}\bigr) \Bigr| 
 \Bigl| \epnwli{i-1} - \epwli{i-1} \Bigr| \nonumber \\[1ex]
 &\quad= L_{S_l}\porosity{l} \Bigl(\porosity{l}S_l\bigl(\pnwln{n} - \pwln{n}\bigr) - \porosity{l}S_l\bigl(\pnwlni{i-1} - \pwlni{i-1}\bigr) \Bigr)
 \Bigl( \epwli{i-1} -\epnwli{i-1} \Bigr) \label{LDD-TPR:eqn:estimate.TP_1R_1.intermediate}
\end{align}
as a result of the Lipschitz continuity of $S_l$. The monotonicity of $S_l$ allowed dropping the absolute value. 
Therefore, by integrating  (\ref{LDD-TPR:eqn:estimate.TP_1R_1.intermediate}), we estimate $T\!P_1$ and $R_1$ by 
\begin{align}
\frac{1}{L_{S_l}\porosity{l}} \bigl\| \Sln{n}  - \Slni{i-1} \bigr\|^2 %
  \leq \spl \Sln{n} - \Slni{i-1}, \epwli{i-1} - \epnwli{i-1} \spr.\label{LDD-TPR:eqn:estimate.TP_1R_1}
\end{align}
For $l=1$, (\ref{LDD-TPR:eqn:estimate.TP_1R_1}) is an estimate for $R_1$ since we had set $\epnwli[1]{i-1}=0$, and for $l=2$ it is an estimate for $T\!P_1$. In this manner, (\ref{LDD-TPR:eqn:estimate.TP_1R_1}) is a condensed notation of both estimates into one. 

\paragraph{$T\!P_2$, $R_2$} 
Young's inequality $\abs{xy}\leq \epsilon \abs{x}^2 + \frac{1}{4\epsilon}\abs{y}^2$ with $\epsilon > 0$ applied to the term $T\!P_2$, gives
\begin{align*}
  \abs{T\!P_2} = \Bigl| \spl \Sln[2]{n} - \Slni[2]{i-1}&, \epwli[2]{i-1} - \epwli[2]{i} -\bigl(\epnwli[2]{i-1} - \epnwli[2]{i}\bigr) \spr \Bigr| \leq \frac{L_{w,2}}{2}\bigl\|\epwli[2]{i-1} - \epwli[2]{i}\bigl\|^2 \nonumber \\[1ex]
  &+ \frac{L_{\nw,2}}{2}\bigl\|\bigl(\epnwli[2]{i-1} - \epnwli[2]{i}\bigr) \bigr\|^2   
  + \left(\frac{1}{2L_{w,2}} + \frac{1}{2L_{\nw,2}}\right)  \bigl\| \Sln[2]{n}- \Slni[2]{i-1}\bigr\|^2, %
\end{align*}
where we chose $\epsilon_{\alpha}^{2} = \frac{L_{\alpha,2}}{2}$ for $\alpha = w, \nw$. 
The analogous choice of $\epsilon_{w}^1 = \frac{L_{w,1}}{2}$ for $l=1$ yields 
\begin{align*}
  \abs{R_2} &= \Bigl| \spl \Sln[1]{n} - \Slni[1]{i-1}, \epwli[1]{i-1} - \epwli[1]{i}\spr \Bigr|
  \leq \frac{L_{w,1}}{2}\bigl\|\epwli[1]{i-1} - \epwli[1]{i}\bigl\|^2  
  + \frac{1}{2L_{w,1}}\bigl\| \Sln[1]{n}- \Slni[1]{i-1}\bigr\|^2. %
\end{align*}

\paragraph{$T\!P_3$, $R_3$} 
$T\!P_3$ and $R_3$ can be estimated together as well. We have 
\begin{align}
  \Bigl|\bspl\Bigl(\kaln{n} &-\kalni{i-1}\!\Bigr)\gradPalnPlusGravity{n},\bm{\nabla}\epali{i}\bspr \Bigr|
  \leq \bigl\|\bigl(\kaln{n} -\kalni{i-1}\bigr)\gradPalnPlusGravity{n}\bigr\| 
  \bigl\|\bm{\nabla}\epali{i}\bigr\| 				\nonumber \\[1ex]
  &\leq \frac{L_{\kaln{}}M}{\porosity{l}}\bigl\|\Sln{n} - \Slni{i}\bigr\| 
  \bigl\|\bm{\nabla}\epali{i}\bigr\| 
  \leq \frac{L_{\kaln{}}M}{\porosity{l}} \epsilon_{\alpha,l}\bigl\|\Sln{n} - \Slni{i}\bigr\|^2 + 
  \frac{L_{\kaln{}}M}{4\epsilon_{\alpha,l}\porosity{l}} 
\bigl\|\bm{\nabla}\epali{i}\bigr\|^2   \label{LDD-TPR:eqn:estimate.R_3}
\end{align}
for $(\alpha, l) \in \{(w,1), (w,2), (\nw,2)\}$.
  To derive (\ref{LDD-TPR:eqn:estimate.R_3}), the Lipschitz-continuity of $k_{\alpha,l}$ and the assumption 
  $\bigl\|\gradPalnPlusGravity{n}\bigr\|_\infty  < M$ was used.
  $\epsilon_{\alpha,l}$ will be chosen later. Equation (\ref{LDD-TPR:eqn:estimate.R_3}) is the estimate for $R_3$ for $l=1, \alpha =w$ and  the estimate for $T\!P_3$ is obtained by summing over the phase index $\alpha$,
\begin{align}
  \abs{T\!P_3} \leq \tau \sum_{\alpha}L_{\kaln[2]{}}\frac{M \epsilon_{\alpha,2}}{\porosity{2}}\bigl\|\Sln[2]{n} - \Slni[2]{i}\bigr\|^2 
  +\tau \sum_{\alpha} \frac{L_{\kaln[2]{}}M}{4\epsilon_{\alpha,2}\porosity{2}} 
\bigl\|\bm{\nabla}\epali[2]{i}\bigr\|^2. \label{LDD-TPR:eqn:estimate.TP_3}
\end{align}

\paragraph{$T\!P_4$, $R_4$}
Finally, by \autoref{env:LDD-TPR:Assumptions:Data}\ref{LDD-TPR:env:assumptions:coercivity}, we estimate $R_4$ by
\begin{align}
  R_4 = \tau \bspl\kwlni[1]{i-1}\bm{\nabla} \epwli[1]{i},\bm{\nabla}\epwli[1]{i}\bspr > \tau m_1 \bigl\|\bm{\nabla} \epwli[1]{i}  \bigr\|^2
 \label{LDD-TPR:eqn:estimate.R_4}
\end{align}
Analogously on the two-phase domain for $T\!P_4$,
we have the estimate 
\begin{align}
 T\!P_4 = \tau \sum_\alpha\bspl\kalni[2]{i-1}\bm{\nabla} \epali[2]{i},\bm{\nabla}\epali[2]{i}\bspr > \tau m_2 \sum_\alpha \bigl\|\bm{\nabla} \epali[2]{i}  \bigr\|^2.
 \label{LDD-TPR:eqn:estimate.TP_4}
\end{align}
Combining the just derived estimates in  (\ref{LDD-TPR:eqn:estimate.TP_1R_1})–(\ref{LDD-TPR:eqn:estimate.R_4}) and (\ref{LDD-TPR:eqn:estimate.TP_4})  with equations (\ref{eqn:mainproofstep_summing.TP}) and (\ref{eqn:mainproofstep_summing.R}), one arrives at
\begin{align}
    \sum_{\alpha}&\frac{L_{\alpha,2}}{2}\Bigl( \bigl\|\epali[2]{i}\bigr\|^2 %
        - \bigl\|\epali[2]{i-1}\bigr\|^2\Bigr) 
    + \frac{1}{L_{S_2}\porosity{2}} \bigl\| \Sln[2]{n} - \Slni[2]{i-1} \bigr\|^2 
    \nonumber \\
    & + \tau \sum_\alpha \spl \lambda_{\alpha}\epali[2]{i} + \egali[2]{i},\epali[2]{i} \spr_\Gamma 
    + \tau m_2 \sum_{\alpha}\bigl\|\bm{\nabla}\epali[2]{i}\bigr\|^2 \leq \sum_{\alpha} \frac{1}{2L_{\alpha,2}}\bigl\| \Sln[2]{n}-\Slni[2]{i-1}\bigr\|^2\nonumber\\
    &\hspace{4.5cm} 
  + \tau \sum_{\alpha} \frac{L_{\kaln[2]{}}M}{4\epsilon_{\alpha,2}\porosity{2}} 
  \bigl\|\bm{\nabla}\epali[2]{i}\bigr\|^2 
  + \tau \sum_{\alpha}L_{\kaln[2]{}}\frac{M\epsilon_{\alpha,2}}{\porosity{2}}\bigl\|\Sln[2]{n} - \Slni[2]{i}\bigr\|^2 %
  \label{eqn:mainproofstep3.TP}
\end{align}
on $\dom{2}$, and on $\dom{1}$, at
\begin{align}
  \frac{L_{w,1}}{2}&\Bigl( \bigl\|\epwli[1]{i}\bigr\|^2 %
        - \bigl\|\epwli[1]{i-1}\bigr\|^2\Bigr) 
    + \frac{1}{L_{S_1}\porosity{1}} \bigl\| \Sln[1]{n} - \Slni[1]{i-1} \bigr\|^2 
    + \tau \spl \lambda_{w}\epwli[1]{i} + \egwli[1]{i},\epwli[1]{i} \spr_\Gamma 
    + \tau m_1\bigl\|\bm{\nabla}\epwli[1]{i}\bigr\|^2 \nonumber\\
  &\qquad \leq \frac{1}{2L_{w,1}} \bigl\| \Sln[1]{n}-\Slni[1]{i-1}\bigr\|^2  
  + \tau \frac{L_{\kaln[1]{}}M}{4\epsilon_{w,1}\porosity{1}}
  \bigl\|\bm{\nabla}\epwli[1]{i}\bigr\|^2 
  + \tau L_{\kwln[1]{}}\frac{M\epsilon_{w,1}}{\porosity{1}}\bigl\|\Sln[1]{n} - \Slni[1]{i}\bigr\|^2. %
  \label{eqn:mainproofstep3.R}   
\end{align}

In order to handle the interface terms $I_R$ in (\ref{eqn:mainproofstep_summing.R}) and $I_{T\!P}$ in (\ref{eqn:mainproofstep_summing.TP}), recall 
that we defined the pressure functionals according to \autoref{note:LDD-TPR:pressureFunctionalsDefiniton}. This allows us to treat the interface terms in the following way: 
Subtracting
(\ref{eqn:LDD-TPR:LDD-scheme:gliupdate:wetting:weak}) from (\ref{eqn:LDD-TPR:limitSystem:glupdate}) for $\alpha=w$ and (\ref{eqn:LDD-TPR:LDD-scheme:gliupdate:nonwetting:weak:1}), (\ref{eqn:LDD-TPR:LDD-scheme:gliupdate:nonwetting:weak:2}) from (\ref{eqn:LDD-TPR:limitSystem:gli:nonwetting:1}) and (\ref{eqn:LDD-TPR:limitSystem:gli:nonwetting:2}) respectively for $\alpha=\nw$, all the while using 
the representations $\overline{\egali{i}} \in \Tracespace$ of the functionals $\egali{i}\in \Tracespace'$%
(cf.\ \autoref{rem_represent}), %
we obtain
\begin{align*}
    \overline{\egali{i}} = -2\lambda_{\alpha} \epaliOnGamma[3-l]{i-1} -  \overline{\egali[3-l]{i-1}}
\end{align*}
for $(\alpha,l)\in\{(w,1),(w,2),(\nw,1)\}$ as well as
\begin{align*} 
    \overline{\egnwli[2]{i}} = -\overline{\egnwli[1]{i-1}}.
\end{align*}
This leads to the relations
  \begin{align}
    \bigl\| \epwliOnGamma[l]{i} \bigr\|^2_{1/2}  &= \frac{1}{  4\lambda_{w}^2 } %
    \left(\bigl\| \overline{\egwli[3-l]{i+1}} \bigr\|^2_{1/2} %
    -\bigl\| \overline{\egwliOnGamma[l]{i}} \bigr\|^2_{1/2} %
    - 4\lambda_{w}  \spl \egwli[l]{i}, \epwliOnGamma[l]{i} \spr_\tGamma \right),
    \label{LDD-TPR:boundarytrickpre:wetting}\\
    \bigl\| \epnwliOnGamma[2]{i} \bigr\|^2_{1/2}  
    &= \frac{1}{  4\lambda_{\nw}^2 } %
    \left(\bigl\| \overline{\egnwli[1]{i+1}} \bigr\|^2_{1/2} %
    -\bigl\| \overline{\egnwliOnGamma[2]{i}} \bigr\|^2_{1/2} %
    - 4\lambda_{\nw}  \spl \egnwli[2]{i}, \epnwliOnGamma[2]{i}\spr_\tGamma \right),
    \label{LDD-TPR:boundarytrickpre:nonwetting2} \\
    0 &= \frac{1}{  4\lambda_{\nw}^2 } %
    \left(\bigl\| \overline{\egnwli[2]{i+1}} \bigr\|^2_{1/2} %
    -\bigl\| \overline{\egnwliOnGamma[1]{i}} \bigr\|^2_{1/2} \right).
    \label{LDD-TPR:boundarytrickpre:nonwetting1} 
  \end{align}
Inserting for $l=1$  equation (\ref{LDD-TPR:boundarytrickpre:wetting}) into (\ref{eqn:mainproofstep3.R}) and for $l=2$ (\ref{LDD-TPR:boundarytrickpre:wetting}) as well as (\ref{LDD-TPR:boundarytrickpre:nonwetting2}) in (\ref{eqn:mainproofstep3.TP}), 
yields
\begin{align}
    \biggl(\frac{1}{L_{S_2}\porosity{2}} - \sum_{\alpha}\Bigl(\frac{1}{2L_{\alpha,2}}  &+ \tau L_{\kaln[2]{}}\frac{M\epsilon_{\alpha,2}}{\porosity{2}}\Bigr)\!\biggr) \bigl\| \Sln[2]{n} - \Slni[2]{i-1} \bigr\|^2 
    +\tau \sum_{\alpha} \left(m_2 - \frac{L_{\kaln[2]{}}M}{4\epsilon_{\alpha,2}\porosity{2}}  \right)  
    \bigl\|\bm{\nabla}\epali[2]{i}\bigr\|^2 
    \nonumber\\
  &\leq \sum_{\alpha}\frac{L_{\alpha,2}}{2}\Bigl(\bigl\|\epali[2]{i-1}\bigr\|^2 - \bigl\|\epali[2]{i}\bigr\|^2 \Bigr) 
  + \tau\sum_{\alpha}\frac{1}{  4\lambda_{\alpha} } %
  \left(\bigl\| \egaliOnGamma[2]{i} \bigr\|^2_{1/2} - \bigl\| \egali[1]{i+1} \bigr\|^2_{1/2} \right) 
  \label{eqn:mainproofstep4.TP} 
\end{align}
for $\dom{2}$ and 
\begin{align}
\biggl(\frac{1}{L_{S_1}\porosity{1}} - \frac{1}{2L_{w,1}}  &- \tau L_{\kwln[1]{}}\frac{M\epsilon_{w,1}}{\porosity{1}}\!\biggr) \bigl\| \Sln[1]{n} - \Slni[1]{i-1} \bigr\|^2
    \tau \left(m_1 - \frac{L_{\kwln[1]{}}M}{4\epsilon_{w,1}\porosity{1}}  \right)  
    \bigl\|\bm{\nabla}\epwli[1]{i}\bigr\|^2 
    \nonumber\\
  &\leq \frac{L_{w,1}}{2}\Bigl(\bigl\|\epwli[1]{i-1}\bigr\|^2 - \bigl\|\epwli[1]{i}\bigr\|^2 \Bigr) 
  + \frac{\tau }{  4\lambda_{w} } %
  \left(\bigl\| \egwliOnGamma[1]{i} \bigr\|^2_{1/2} - \bigl\| \egwli[2]{i+1} \bigr\|^2_{1/2} \right),
  \label{eqn:mainproofstep4.R}
\end{align}
for $\dom{1}$, 
where we dropped denoting the representing element by $\bar{\cdot}$.

Now choose $\epsilon_{\alpha,l} = L_{\kaln{}}M/(2m_l\porosity{l})$ such that $m_l -L_{\kaln{}}M/(4\epsilon_{\alpha,l}\porosity{l}) = \frac{m_l}{2} > 0$. 
Recall that by assumption the numbers $L_{\alpha,l}$ have been chosen large enough that 
$\frac{1}{L_{S_2}\porosity{2}} - \frac12\sum_{\alpha}\frac{1}{L_{\alpha,2}} 
> 0$ as well as $\frac{1}{L_{S_1}\porosity{1}} - \frac{1}{2L_{w,1}} >0$,
and in addition the time step restriction (\ref{eqn:timesteprestriction.2}) is satisfied for a sufficiently small $\tau$.
Summing up the equations (\ref{eqn:mainproofstep4.TP}) and (\ref{eqn:mainproofstep4.R}), then adding zero in the form of (\ref{LDD-TPR:boundarytrickpre:nonwetting1}) to the result and thereafter summing with respect to iterations $i=1,\dots,r$  leads to
\begin{align}
 \sum_{i=1}^r &\sum_{l=1,2} C_l \bigl\| \Sln{n} - \Slni{i-1} \bigr\|^2 
  +\tau \sum_{i=1}^r\left(\frac{m_1}{2}\bigl\|\bm{\nabla}\epwli[1]{i}\bigr\|^2 +\frac{m_2}{2}\sum_{\alpha}   
  \bigl\|\bm{\nabla}\epali[2]{i}\bigr\|^2 \right)
  \nonumber\\
  &\hspace*{-.3cm} \leq \sum_{\alpha}\frac{L_{\alpha,2}}{2}\Bigl(\bigl\|\epali[2]{0}\bigr\|^2 - \bigl\|\epali[2]{r}\bigr\|^2 \Bigr)  
  + \frac{L_{w,1}}{2}\Bigl(\bigl\|\epwli[1]{0}\bigr\|^2 - \bigl\|\epwli[1]{r}\bigr\|^2 \Bigr)  %
  + \tau \sum_{l=1}^2\sum_\alpha\frac{1}{  4\lambda_{\alpha} } %
  \left(\bigl\| \egaliOnGamma{1} \bigr\|^2_{1/2} - \bigl\| \egali{r+1} \bigr\|^2_{1/2} \right),
  \label{eqn:mainproofstep5.TPR}
\end{align}
where $C_l$ is defined in (\ref{eqn:timesteprestriction.2}) and the telescopic nature of the sums on the right hand side have been exploited. 
Equation (\ref{eqn:mainproofstep5.TPR}) implies the estimates
\begin{align}
  \sum_{i=1}^r\sum_{l=1}^2C_l
  \bigl\| \Sln{n} - \Slni{i-1} \bigr\|^2  
  &\leq C, \label{eqn:mainproofstep6a} \\%
  \tau \sum_{i=1}^r\left(\frac{m_1}{2}\bigl\|\bm{\nabla}\epwli[1]{i}\bigr\|^2 +\frac{m_2}{2}\sum_{\alpha}   
  \bigl\|\bm{\nabla}\epali[2]{i}\bigr\|^2 \right)%
  &\leq C,
  \label{eqn:mainproofstep6b} \\
  \sum_{\alpha}\frac{L_{\alpha,2}}{2}\bigl\|\epali[2]{r}\bigr\|^2  
  + \frac{L_{w,1}}{2}\bigl\|\epwli[1]{r}\bigr\|^2
  + \tau \sum_{l=1}^2\sum_\alpha\frac{1}{  4\lambda_{\alpha} } %
  \bigl\| \egali{r+1} \bigr\|^2_{1/2} 
  &\leq C, \label{eqn:mainproofstep6c}
    \end{align}
  {with}
  \[
  C := \sum_{\alpha}\frac{L_{\alpha,2}}{2}\bigl\|\epali[2]{0}\bigr\|^2 
  + \frac{L_{w,1}}{2}\bigl\|\epwli[1]{0}\bigr\|^2 
  + \tau \sum_{l=1}^2\sum_\alpha\frac{1}{  4\lambda_{\alpha} } %
  \bigl\| \egaliOnGamma{1} \bigr\|^2_{1/2}.
\]
Since $C$ is independent of $r$, we thereby conclude that %
\begin{align}
    \bigl\| \Sln{n} - \Slni{i-1} \bigr\|, \,\bigl\|\gradepali{i} \bigr\|\longrightarrow 0 \quad \mbox{as } i\rightarrow \infty, 
    \label{eqn:LDD-TPR:SlAndGradPconverge}
\end{align}
for all appearing combinations of $l\in\{1, 2\}$ and $\alpha\in \{w, \nw\}$.
Due to the partial homogeneous Dirichlet boundary, the Poincaré inequality is applicable %
  cf. \cite[Theorem A.2.5, p. 252]{Berninger2009}
  for functions in $\Fs_l$. 
  Thus, equation (\ref{eqn:LDD-TPR:SlAndGradPconverge}) further 
  implies $\bigl\|\epali{i} \bigr\|$ $\rightarrow 0$ 	 as $i\rightarrow \infty$ for all admissible index combinations.

  In order to show that $\egali{i} \rightarrow 0$ in $\Fs_l'$ for all appearing indices, we subtract  
  (\ref{eqn:LDD-TPR:LDD-scheme:wetting}) from (\ref{eqn:LDD-TPR:limitSystem:wetting}) for $\alpha=w$,
  $l=1,2$, and (\ref{eqn:LDD-TPR:LDD-scheme:nonwetting}) from (\ref{eqn:LDD-TPR:limitSystem:nonwetting}) 
  for $\alpha = \nw$, $l=2$
  and consider only test functions in $\testfunc{\alpha} \in C_0^\infty(\dom{l})$, i.e.
  \begin{align}
    -\tau \bspl\fluxn{\alpha} - \fluxai{i},\bm{\nabla}\testfunc{\alpha}\bspr
      &= L_{\alpha,l}\bspl\epali{i-1} -\epali{i},\testfunc{\alpha} \bspr %
    + (-1)^{\krond{w}}\bspl \Sln{n}- \Slni{i-1}, \testfunc{\alpha} \bspr. \label{eqn:giconvergencestep1}
  \end{align}
  Thus, $\dv\bigl(\fluxn{\alpha} - \fluxai{i}\bigr)$ exists in $L^2(\dom{l})$ and
   \begin{align}
      -\tau \dv\bigl(\fluxn{\alpha} - \fluxai{i}\bigr)
      = L_{\alpha,l} \bigl(\epali{i} - \epali{i-1}\bigr)  
      +(-1)^{\krond{\nw}}\bigl( \Sln{n}- \Slni{i-1}\bigr) \label{eqn:giconvergencestep2}
   \end{align}
  almost everywhere, from which we deduce for $\testfunc{\alpha}$ now taken to be in $\Fs_l$
  \begin{align}
      \Bigl|\bspl\dv\bigl(\fluxn{\alpha} - \fluxai{i}\bigr),\testfunc{\alpha} \bspr\Bigr| 
	&\leq \frac{L_{\alpha,l}}{\tau}\bigl\|\epali{i} - \epali{i-1} \bigr\|\bigl\|\testfunc{\alpha}\bigr\| 
    + \frac{1}{\tau}\bigl\| \Sln{n}- \Slni{i-1}\bigr\| \bigl\| \testfunc{\alpha} \bigr\|.
      \label{eqn:giconvergencestep3}
  \end{align}
  Introducing the abbreviation $\bigl|\Psi_{\alpha,l}^{n,i}\bigl(\testfunc{\alpha}\bigr)\bigr|$ for the left hand side of (\ref{eqn:giconvergencestep3}), the limit 
  \begin{align*}
    \sup_{\stackrel{\testfunc{\alpha} \in \Fs_l}{\testfunc{\alpha} \neq 0}} 
    &\frac{\bigl|\Psi_{\alpha,l}^{n,i}\bigl(\testfunc{\alpha}\bigr)\bigr|}{\norm{\testfunc{\alpha}}_{\Fs_l}} 
    \leq \frac{L_{\alpha,l}}{\tau}\bigl\|\epali{i} - \epali{i-1} \bigr\|
    + \frac{1}{\tau}\bigl\| \Sln{n}- \Slni{i-1}\bigr\|
    \longrightarrow 0 %
  \end{align*}
  as $i \rightarrow \infty$ follows as a consequence of (\ref{eqn:mainproofstep6a}) and (\ref{eqn:mainproofstep6b}). In other words $\bigl\|\Psi_{\alpha,l}^{n,i} \bigr\|_{\Fs_l'} \rightarrow 0$ as 
$i 
  \rightarrow \infty$.
  On the other hand, starting again from (\ref{LDD-TPR:mainproofstep1a}) (without the added zero term), this time however 
inserting   $\testfunc{\alpha} \in \Fs_l$ and integrating by parts, keeping in mind (\ref{eqn:giconvergencestep2}), 
one deduces that
  \begin{align}
    \bspl \egali{i},\testfunc{\alpha} \bspr_\Gamma   
      = - \lambda_{\alpha} &\bspl\epaliOnGamma{i},\testfunc{\alpha} \bspr_\Gamma 
      +\bspl \bigl( \fluxn{\alpha}- \fluxai{i} \bigr)\cdot\vt{n_l},\testfunc{\alpha} 
  \bspr_\Gamma . \label{eqn:giconvergencestep5}
  \end{align}
    We already know, that $\bigl\|\epali{i}\bigr\|_{\Fs_l} \rightarrow 0$ as $i \rightarrow 0$ and we will use the 
continuity   of the trace operator to deal with the term $\spl\epaliOnGamma{i},\testfunc{\alpha} \spr_\Gamma $. 
For the last summand in (\ref{eqn:giconvergencestep5}) we have by the 
integration by parts formula
  \begin{align}
    \bspl \bigl( \fluxn{\alpha}- \fluxai{i}  \bigr)\cdot\vt{n_l},\testfunc{\alpha} 
  \bspr_\Gamma 
    &= \Psi_{\alpha,l}^{n,i}(\testfunc{\alpha}) 
      + \bspl \fluxn{\alpha}- \fluxai{i}, \bm{\nabla}\testfunc{\alpha}\bspr,
  \label{eqn:giconvergencestep6}
  \end{align}
  and the second term can be estimated by
  \begin{align*}
    \Bigl|\bspl  \kaln{n} \gradPalnPlusGravity{n}  
    &- \kalni{i-1}\gradPalniPlusGravity{i} 
    ,\bm{\nabla}\testfunc{\alpha}\bspr\Bigr| 
    \nonumber \\
    &
    \leq
    \Bigl|\bspl \bigl( \kaln{n} -\kalni{i-1}\bigr) \gradPalnPlusGravity{n} 
    - \kalni{i-1}\gradepali{i} ,\bm{\nabla}\testfunc{\alpha}\bspr\Bigr|
\nonumber   \\%
    &\leq  \frac{L_{\kaln{}} M}{\porosity{l}} \bigl\|\Sln{n} - \Slni{i-1}\bigr\| \bigl\|\testfunc{\alpha}\bigr\|_{\Fs_l} + 
    M_{\kaln{}}\bigl\|   \gradepali{i}\bigr\| \bigl\|\testfunc{\alpha} \bigr\|_{\Fs_l},
  \end{align*}
  where we used the same reasoning as in (\ref{LDD-TPR:eqn:estimate.R_3}) and $\max|\kaln{}| \leq M_{\kaln{}}$. 
  With this, we get 
  \begin{align*}
    \sup_{\stackrel{\testfunc{\alpha} \in \Fs_l}{\norm{\testfunc{\alpha}}_{\Fs_l} = 1}}
    \Bigl|\bspl \bigl( \fluxn{\alpha}&- \fluxai{i} \bigr)\cdot\vt{n_l},\testfunc{\alpha}\bspr_\Gamma\Bigr| 
    \leq \bigl\|\Psi_{\alpha,l}^{n,i}\bigr\|_{\Fs_l'}
    + \frac{L_{\kaln{}} M}{\porosity{l}} \bigl\|\Sln{n} - \Slni{i-1}\bigr\|  + M_{\kaln{}}\bigl\| \gradepali{i}\bigr\| 
  \longrightarrow 0,  %
  \end{align*}
  as $i \rightarrow \infty $ from (\ref{eqn:giconvergencestep6}). 
  Finally, we deduce from (\ref{eqn:giconvergencestep5}) and the continuity of the trace operator (with constant $\tilde{C}$) 
  on Lipschitz domains
   \begin{align*}
    \sup_{\stackrel{\testfunc{\alpha} \in \Fs_l}{\testfunc{\alpha} \neq 0}}
    \frac{\bigl|\spl \egali{i},\testfunc{\alpha} \spr_\Gamma\bigr|}{\norm{\testfunc{a}}_{\Fs_l}} 
    &\leq \lambda_{\alpha}\tilde{C}\bigl\| \epali{i} \bigr\|_{\Fs_l} 
      + \bigl\|\Psi_{\alpha,l}^{n,i}\bigr\|_{\Fs_l'} 
      + \frac{L_{\kaln{}} M}{\porosity{l}} \bigl\|\Sln{n} - \Slni{i-1}\bigr\|  
      + M_{\kaln{}}\bigl\| \gradepali{i} \bigr\| \longrightarrow 0 , 
  \end{align*}
  as $i \rightarrow \infty$.
  This shows $\egali{i} \rightarrow 0$ in $\Fs_l'$ for all valid index combinations and concludes the proof. \qed

\section{The LDD-TPR solver for the  multi-domain case} 

\label{LDD-TPR:section:problemFormulation:multiDomain}
In this section we provide a generalisation of \autoref{env:LDD-TPR:problemFormulation:semiDiscrete:weak} 
to a multi-domain setting.  An in-depth presentation  with a multi-domain convergence result can be found in \cite[Section 4.4]{SeusThesis2021}.

We start with  a generalisation of our geometric notations.
   \label{LDD-TPR:Notation:multidomain}
The domain  $\dom{} \subset \RR^d$ is partitioned into  a finite number of non-overlapping Lipschitz subdomains $\dom{l} \subset 
	\dom{}$  such that
	$\overline{\dom{}} = \bigcup_{l=1}^{\subdomNum} \overline{\dom{l}}$. 
	The interior of intersections of the boundary of neighbouring domains, 
	$\intf{kl} := \overline{\dom{k}} \cap \overline{\dom{l}}\setminus\partial(\overline{\dom{k}}\cup \overline{\dom{l}})$ that in addition have non-zero $(d-1)$-dimensional Hausdorff measure are called interfaces and are submanifolds 
	of dimension $d-1$. 
	As a consequence, the outer normal $\vt{n_{kl}} \in S^{d-1}$ 
	pointing from $\dom{k}$ to $\dom{l}$ is defined almost everywhere on $\intf{kl}$, cf. \cite[p. 97 ff.]{MacLean2000} and \cite{Agranovich2015} for more details on definitions.
	Figure \ref{fig:LDD-TPR:NumericalDomain}(\subref{fig:LDD-TPR:NumericalDomain:FiveDomainInnerPatch}), p.\ \pageref{fig:LDD-TPR:NumericalDomain}, %
    illustrates the notation.
	 Given $\dom{l}$, let $\ind_l \subset \ind = \{1, \dots, \subdomNum\}$ be the set of indices 
	denoting those neighbouring subdomains $\dom{k}$ for which $\intf{lk}$ is an interface.
	Furthermore, let $\intf{l} := \text{\normalfont int}\bigl(\partial\dom{} \cap \overline {\dom{l}}\bigr)$ denote that particular part of the 
	boundary of $\dom{l}$ intersecting with $\partial\dom{}$ for all $l\in \ind$ for which $\intf{l}$ 
	is $(d-1)$-dimensional. Then $\partial \dom{} = \bigcup_{l\in\ind}\overline{\intf{l}}$ and 
	$\partial \dom{l} =  \bigcup_{k\in\ind_{l}}\bigl(\overline{\intf{kl}}  \cup \overline{\intf{l}} \bigr)$. 

Let $\indR, \indTP \subset \ind$ with $ \ind = \union{\indR}{\indTP}$ be  the (possibly empty) 
	sets of indices denoting the subdomains $\dom{l}$ on which the \Richards equation $(l\in \indR)$ or 
	the full two-phase system $(l \in \indTP)$ is imposed.  
	If neither $\indR \neq \emptyset$ nor $\indTP \neq \emptyset$, denote by $\indTPR \subset 
	\indTP $ the indices of domain patches, that model the full two-phase flow but have at
	least one
	neighbouring subdomain that assumes the \Richards model sharing an interface of dimension 
	$(d-1)$. Similarly, define the subset $\indRTP \subset \indR$ of \Richards subdomains with 
	a two-phase neighbour. Given a subdomain $\dom{k}$ with $k \in \indTPR$, let 
	$\indTPRl{k} \subset \ind_k$ be the set of indices of neighbouring \Richards subdomains. 
	Analogously, for $\dom{k}$ with $k \in \indRTP$, let $\indRTPl{k} \subset \ind_k$ be the set of 
	indices of neighbouring two-phase subdomains.

Turning to function spaces, we decompose $H^1_0(\dom{})$  %
   into spaces
  \begin{align*}
   \Fs_l &:= %
      \left\{ u \in H^1(\dom{l})\, \bigl|\bigr. \, \restrictTo{u}{\oB{l}}
      = 0 \right\} 
    \end{align*}
    and set 
    \begin{align*}
  \prodFs[]{w}   &:= \prod_{j=1}^\subdomNum\Fs_j, 
  \quad \mbox{and }  \quad
  \Fs^w   = \bigl\{ (u_1,\dots, u_\subdomNum) \in \prodFs[]{w} \, \bigl|\bigr. \, \restrictTo{u_l}{\intf{lk}} = \restrictTo{u_k}{\intf{lk}}, l \in \ind, k\in \ind_l \bigr\},
  \end{align*}
  the latter being $H^1_0(\dom{})$.
  In case $\dom{k}$ is an internal subdomain, i.e. $\oB{k} = \emptyset$, we have $\Fs_k = H^1(\dom{k})$. \\
  If a purely \Richards or purely two-phase domain decomposition is considered, i.e. no
  \TPR coupling occurs, the space $\Fs^w$ can be used for the wetting phases, (Richards) and also the nonwetting phases (two-phase). 
  If both models are present, however, we refine the notion of $\Fs^{\nw}$, the space for the nonwetting phase.
  For all $k\in \indTP$ define $\Fs^{\nw}_l := \Fs_l$, $\Fs^{\nw}_k = \{0\}$ for $k\in \indR$ and set
  the general space for the nonwetting phase 
  \begin{align*}
    \prodFs[]{\nw}   &:= \prod_{j=1}^\subdomNum\Fs_j %
    \quad \mbox{and } \quad
    \Fs^{\nw}   := \bigl\{ (u_1,\dots, u_\subdomNum) \in \prodFs[]{\nw} \, \bigl|\bigr. \, \restrictTo{u_l}{\intf{lk}} = \restrictTo{u_k}{\intf{lk}}, \, l \in \indTP \bigr\}.
  \end{align*}
  Note, that since for $k\in \indR$ we set $\Fs^{\nw}_k = \{0\}$, 
  $\restrictTo{u_l}{\intf{lk}} = \restrictTo{u_k}{\intf{lk}}$ for $k\in \indTPRl{l}$
  actually means $\restrictTo{u_l}{\intf{lk}} = 0$.
  Lastly, in order to define Neumann traces, we need those subspaces of the spaces $\Fs_{k}^{\alpha}$ and $\Fs^{\alpha}$ for which traces on each interface to neighbours can be extended by zero. We set
  \begin{align*}
    \Fs_{k,00}^{\alpha} &:= \bigl\{ u\in \Fs_{k}^{\alpha} \, \bigl| \, \restrictTo{u}{\intf{kl}} \in \Ts{kl},\, l\in \ind_k \bigr\} %
    \quad \mbox{and}\quad
    \prodFs[00]{\alpha}:= \prod_{j=1}^\subdomNum \Fs^{\alpha}_{j,00},
  \end{align*}
  for $\alpha\in \{w, \nw\}$.

With the above notations a multi-domain semi-discrete formulation of \autoref{env:LDD-TPR:problemFormulation:semiDiscrete:weak} reads as 
\begin{problem}[Semi-discrete \TPR problem, multi-domain]~
\label{env:LDD-TPR:problemFormulation:semiDiscrete:multi-domain:weak} 
Given functions 
$(p_w^{n-1}, p_{\nw}^{n-1}) \in \Fs^w\times\Fs^{\nw}$,  
find 
$(p_w^{n}, p_{\nw}^{n}) \in \Fs^w\times\Fs^{\nw}$,
such that all fluxes fulfil $\fluxn[l]{\alpha}\cdot \vt{n_{lk}} \in \Ts{lk}'$ for $l \in \ind$, $k\in \ind_l$, $\alpha \in \{w, \nw\}$, and the equations
\begin{align*}
     \spl \Sln{n} - \Sln{n-1},\testfunc[l]{w}\spr  %
      &- \tau\spl \fluxn[l]{w},\bm{\nabla}\testfunc[l]{w} \spr 
      + \tau \sum_{k\in\ind_l}\spl \fluxn[k]{w}\cdot \vt{n_{lk}},\testfunc[l]{w}\spr_{\intf{lk}}
      =  \tau\spl  f_{w,l}^n ,\testfunc[l]{w} \spr
\end{align*}
as well as 
\begin{align*}
     \spl \Sln[j]{n} - \Sln[j]{n-1},\testfunc[j]{w}\spr  %
      &- \tau\spl \fluxn[j]{w},\bm{\nabla}\testfunc[j]{w} \spr 
      + \tau \sum_{k\in\ind_j} \spl \fluxn[k]{w}\cdot \vt{n_{jk}},\testfunc[j]{w}\spr_{\intf{jk}}
      =  \tau\spl  f_{w,j}^n ,\testfunc[j]{w}\spr,  
    \\[1ex]
    -\spl \Sln[j]{n} - \Sln[j]{n-1},\testfunc[j]{\nw}\spr  %
      &- \tau\spl \fluxn[j]{\nw},\bm{\nabla}\testfunc[j]{\nw} \spr
      + \tau \sum_{k\in \ind_j} \spl \fluxn[k]{\nw}\cdot \vt{n_{jk}},\testfunc[j]{\nw}\spr_{\intf{jk}}
      =  \tau\spl  f_{\nw,j}^n ,\testfunc[j]{\nw} \spr
  \end{align*}
  are satisfied for  $l \in \indR$, $j \in \indTP$ and for all $(\testfunc[1]{\alpha},\testfunc[2]{\alpha}, \dots, \testfunc[\subdomNum]{\alpha}) \in \prodFs[00]{\alpha}$.
\end{problem}

As we have seen in the previous section, cf.\ \autoref{lemma:LDD-TPR:iterationinterpretationlemma}, introducing a Robin type formulation allows to drop the pressure continuity that is implicitly contained in the definition of our spaces $\Fs^w$
and $\Fs^{\nw}$. Instead, the pressure continuity becomes part of the equations to solve and is thereby more accessible to implementation.
  Setting analogously to equations (\ref{LDD-TPR:eqn:gwln:wetting})–(\ref{LDD-TPR:eqn:gwln:nonwetting2}) for $\lambda_\alpha^{lk}=\lambda_\alpha^{kl} > 0$
  \begin{align*}
  \gali[lk]{} &:= -\lambda_{\alpha}^{lk}\restrictTo{\paln{n}}{\intf{lk}} + \fluxn{\alpha}\cdot \vt{n_{lk}}, \quad (l\in \ind), \, (k\in \ind_l) \mbox{ and }
   \alpha \in \{n, \nw\},%
\end{align*}
in $\dualTs{lk}$,  where as in the previous sections $\gnwli[lk]{} = \fluxn{\nw}\cdot \vt{n_{lk}}$ on \Richards subdomains $l\in\indR$ are only 
  nonzero if the neighbour $\dom{k}$ assumes the two-phase model, i.e., $k\in\indRTPl{l}$ and gravity is included,   \autoref{env:LDD-TPR:problemFormulation:semiDiscrete:multi-domain:weak}   can be equivalently reformulated into

\begin{problem}[Semi-discrete \TPR problem,  limit formulation, multi-domain]
\label{problem:LDD-TPR:semiDiscrete:multi-domain:limitSystem}
  Let functions 
$(p_w^{n-1}, p_{\nw}^{n-1}) \in \Fs^w\times\Fs^{\nw}$ as well as  real numbers $\lambda_\alpha^{lk}=\lambda_\alpha^{kl} > 0$ be given for all interfaces $\intf{lk}$, $l\in \ind$, $k\in\ind_l$, and appearing phases $\alpha$.  

Find 
$(p_w^{n}, p_{\nw}^{n}) \in \prodFs{w}\times\prodFs{\nw}$ and 
$\gali[lk]{} \in \dualTs{lk}$, $l \in \ind$, $k\in \ind_l$, $\alpha \in \{w, \nw\}$ such that on \Richards domains, i.e., $l\in \indR$,
the terms $\gali[lk]{}$ are given by
  \begin{align}
    \begin{aligned}
    \gnwli[lk]{} &=  \fluxn[l]{\nw}\cdot \vt{n_{lk}}, &&\mbox{for} \quad l\in\indRTP, \quad\mbox{and}\quad k\in \indRTPl{l}, \\[1ex]
    \gnwli[lk]{} &=  0, &&\mbox{for}\quad l\in\indRTP \quad\mbox{and}\quad k\in \ind_l\setminus\indRTPl{l},
    \\[1ex]
    \gnwli[lk]{} &=  0, &&\mbox{for} \quad l\in\indR\setminus\indRTP, \quad\mbox{and}\quad k\in \ind_j,
    \end{aligned}
    \label{eqn:LDD-TPR:multi-domain:gnwlRFluxes}
  \end{align}
  where $\fluxn[l]{\nw}\cdot \vt{n_{lk}}\in \dualTs{lk}$ are defined by the right hand side of either 
(\ref{eqn:LDD-TPR:ContinuousCouplingConditions:strong:alternative})
or (\ref{eqn:LDD-TPR:CouplingConditions:continuous:zero_flux}),
  and these functions $(p_w^{n}, p_{\nw}^{n})$ and 
$\gali[lk]{}$ fulfil the equations
  \begin{align}
     \spl \Sln{n} - \Sln{n-1},\testfunc[l]{w}\spr  %
      &- \tau\spl \fluxn[l]{w},\bm{\nabla}\testfunc[l]{w} \spr 
      + \sum_{k\in\ind_l}\tau \spl \lambda_{w}^{lk}\pwln{n} + \gwli[lk]{},\testfunc[l]{w}\spr_{\intf{lk}}
      =  \tau\spl  f_{w,l}^n ,\testfunc[l]{w} \spr,
      \label{eqn:LDD-TPR:multi-domain:limitsystem:Richards} 
\end{align}
for $l\in \indR$ with 
\begin{align}
   \spl \gwli[lk]{},\testfunc{w}\spr_{\intf{lk}} &= \spl -2\lambda_{w}^{lk}
  \pwln[k]{n} -   \gwli[kl]{},\testfunc{w}\spr_{\intf{lk}} 
  && (k\in \ind_l),
  \label{eqn:LDD-TPR:multi-domain:Rdomains:galupdate:wetting} \\[1ex]
  \spl \gnwli[lk]{},\testfunc[k]{\nw}\spr_{\intf{lk}} &= \spl -2\lambda_{\nw}^{lk}
  \pnwln[k]{n} -   \gnwli[kl]{},\testfunc[k]{\nw}\spr_{\intf{lk}} 
  && (k\in \indRTPl{l}),
  \label{eqn:LDD-TPR:multi-domain:Rdomains:galupdate:nonwetting} 
  \end{align}
and for $j\in\indTP$ the equations
\begin{align}
     \spl \Sln[j]{n} - \Sln[j]{n-1},\testfunc[j]{w}\spr  %
      &- \tau\spl \fluxn[j]{w},\bm{\nabla}\testfunc[j]{w}\spr + \tau \sum_{k\in\ind_j} \spl \lambda_{w}^{jk}\pwln{n} + \gwli[jk]{},\testfunc[j]{w}\spr_{\intf{jk}}
      =  \tau\spl  f_{w,j}^n ,\testfunc[j]{w}\spr,  
      \label{eqn:LDD-TPR:multi-domain:TPdomains:limitsystemTPR:TP:w}\\[1ex]
    -\spl \Sln[j]{n} - \Sln[j]{n-1},\testfunc[j]{\nw}\spr  %
      &- \tau\spl \fluxn[j]{\nw},\bm{\nabla}\testfunc[j]{\nw} \spr+ \tau \sum_{k\in \ind_j} \spl \lambda_{\nw}^{jk}\pnwln{n} + \gnwli[jk]{},\testfunc[j]{\nw}\spr_{\intf{jk}}
      =  \tau\spl  f_{\nw,j}^n ,\testfunc[j]{\nw} \spr,
     \label{eqn:LDD-TPR:multi-domain:TPdomains:limitsystemTPR:TP:nw}
\end{align}
together with
\begin{align}
  \spl \gwli[jk]{},\testfunc[j]{w}\spr_{\intf{jk}} &= \spl -2\lambda_{w}^{jk}
  \pwln[k]{n} -   \gwli[kj]{},\testfunc[j]{w}\spr_{\intf{jk}} 
  && (k\in \ind_j), 
  \label{eqn:LDD-TPR:multi-domain:TPdomains:galupdate:wetting} \\[1ex]
  \spl \gnwli[jk]{},\testfunc[j]{w}\spr_{\intf{jk}} &= \spl -2\lambda_{\nw}^{jk}
  \pnwln[k]{n} -   \gnwli[kj]{},\testfunc[j]{\nw}\spr_{\intf{jk}} 
  && (k\in \ind_j \cap \indTP),
  \label{eqn:LDD-TPR:multi-domain:TPdomains:galupdate:nonwetting:1}\\[1ex]
  \spl \gnwli[jk]{},\testfunc[j]{\nw}\spr_{\intf{jk}} 
       &= \spl -\gnwli[kj]{},\testfunc[j]{\nw}\spr_{\intf{jk}} 
  && (k\in \indTPRl{j}),
  \label{eqn:LDD-TPR:multi-domain:TPdomains:galupdate:nonwetting:2}
\end{align}
 for all $(\testfunc[1]{\alpha},\testfunc[2]{\alpha}, \dots, \testfunc[\subdomNum]{\alpha}) \in \prodFs[00]{\alpha}$.
\end{problem}%

\begin{remark}[Necessity of ${\prodFs[00]{\alpha}}$ as test function space]
  In order to define Neumann traces on parts of a boundary 
     $\interface \subset\partial\dom{}$ of a Lipschitz domain $\dom{}$, test functions 
     $\varphi\in H^{1/2}(\interface)$ need to be extendable by zero and these are precisely the functions in $\Tracespace$. If we tested the above problems with functions $\testfunc[l]{\alpha}\in\Fs_l^{\alpha}$ the traces $\restrictTo{\varphi_l}{\interface_{lk}}$ for $k\in \ind_l$ a priori would only lie in $H^{1/2}(\interface_{lk})$. 
     For the Neumann traces appearing in \autoref{env:LDD-TPR:problemFormulation:semiDiscrete:multi-domain:weak}
     and \autoref{problem:LDD-TPR:semiDiscrete:multi-domain:limitSystem} to be well-defined, we need $\restrictTo{\varphi_l}{\interface_{lk}} \in \Ts{lk}$ for $k\in \ind_l$, however. Testing with $\testfunc{\alpha}\in {\prodFs[00]{\alpha}}$ precisely alleviates that problem.
     
\end{remark}
As before, \autoref{problem:LDD-TPR:semiDiscrete:multi-domain:limitSystem} shows how to design the  multi-domain \LDDTPR solver step.
\begin{problem}[\LDDTPR solver step, multi-domain, version 1] 
\label{env:LDD-TPR:LDD-scheme:multi-domain:weak:ver1}
Given $(p_w^{n-1}, p_{\nw}^{n-1}) \in \Fs^w \times \Fs^{\nw}$,  
set  
on all subdomains $\dom{l}$, $l\in\ind$, for some $L_{\alpha,l} > 0$ as initial iterates 
\begin{align}
    \palni{0}:= \paln{n-1},
    \label{eqn:LDD-TPR:LDDsolver:intialIterate:pressure:ver1}
\end{align}
where $\alpha \in \{w\}$ for $l\in \indR$ and $\alpha \in \{w, \nw\}$ for $l\in \indTP$
as well as 
\begin{align}
    \gali[lk]{0} := \fluxb[n-1]{\alpha} \cdot \vt{n_{lk}} - \lambda_{\alpha}^{lk} \restrictTo{\paln[l]{n-1}}{\intf{lk}}
    \label{eqn:LDD-TPR:LDDsolver:intialIterate:gla:ver1}
\end{align}
in $\dualTs{lk}$ for $k\in \ind_l$ and $\lambda_{\alpha}^{lk}=\lambda_{\alpha}^{kl} > 0$.
As before, for \Richards domains, $\dom{l}$ with $l\in\indR$, and on interfaces $\interface_{lk}$ to a two-phase domain, i.e. $k\in\indRTPl{l}$, equation (\ref{eqn:LDD-TPR:LDDsolver:intialIterate:gla:ver1}) becomes 
\begin{align}
    \gali[lk]{0} := \fluxb[n-1]{\alpha} \cdot \vt{n_{lk}}
    \label{eqn:LDD-TPR:LDDsolver:intialIterate:gla2:ver1}
\end{align}
and the fluxes $\fluxb[n-1]{\alpha} \cdot \vt{n_{lk}}$ are defined by the right hand side of either 
(\ref{eqn:LDD-TPR:ContinuousCouplingConditions:strong:alternative})
or (\ref{eqn:LDD-TPR:CouplingConditions:continuous:zero_flux}). On interfaces $\interface_{lk}$ between \Richards domains, we take
\begin{align}
    \gnwli[lk]{0}=0.
    \label{eqn:LDD-TPR:LDDsolver:intialIterate:gla3:ver1}
\end{align}
Given the iterates
$(p_{w}^{n,i-1}, p_{\nw}^{n,i-1}) \in \prodFs{w}\times\prodFs{\nw}$, as well as $\gali[lk]{i-1} \in \Ts{lk}' $,  $(\NN\in i \geq 1)$, one step of the \LDDTPR solver consists of finding 
$(p_{w}^{n,i}, p_{\nw}^{n,i}) \in \prodFs{w}\times\prodFs{\nw}$
such that on \Richards subdomains, i.e. $l \in \indR$, the equations
\begin{align}
  L_{w,l}\spl \pwlni{i},\testfunc{w} \spr
    &- \tau\spl \fluxwi{i},\bm{\nabla}\testfunc{w} \spr
    + \tau \sum_{k\in\ind_l}\spl \lambda_w^{lk}\pwlni{i}
    + \gwli[lk]{i} ,\testfunc{w}\spr_{\intf{lk}}   \nonumber\\[1ex]
  &= L_{w,l}\spl \pwlni{i-1},\testfunc{w}\spr
      - \spl \Slni{i-1} - \Sln{n-1},\testfunc{w}\spr
    + \tau\spl  f_{w,l}^n ,\testfunc{w} \spr,   
  \label{eqn:LDD-TPR:LDD-scheme:multi-domain:Richards:ver1}	
\end{align}
with
\begin{align}
   \spl \gwli[lk]{i},\testfunc{w}\spr_{\intf{lk}}  &:= \spl -2\lambda_{w}^{lk}
  \pwlni[k]{i-1} -   \gwli[kl]{i-1},\testfunc{w}\spr_{\intf{lk}}
  && (k\in \ind_l),
  \label{eqn:LDD-TPR:LDD-scheme:multi-domain:gliupdate:Richards:wetting:weak:ver1} \\[1ex]
  \spl \gnwli[lk]{i},\testfunc[k]{\nw}\spr_{\intf{lk}}  &:= \spl -2\lambda_{\nw}^{lk}
  \pnwlni[k]{i-1} -   \gnwli[kl]{i-1},\testfunc[k]{\nw}\spr_{\intf{lk}} 
  && (k\in \indRTPl{l})
  \label{eqn:LDD-TPR:LDD-scheme:multi-domain:gliupdate:Richards:nonwetting:weak:ver1} 
  \end{align}
 are satisfied, and on two-phase domains, $(j\in\indTP)$, the equations
  \begin{align}
  L_{\alpha,j}\spl \palni[j]{i},\testfunc[j]{\alpha} \spr
    &- \tau\spl \fluxai[j]{i},\bm{\nabla}\testfunc[j]{\alpha} \spr
    + \tau \sum_{k\in\ind_j} \spl \lambda_{\alpha}^{jk}\palni[j]{i}
    + \gali[jk]{i} ,\testfunc[j]{\alpha}\spr_{\intf{jk}} \nonumber \\[1ex]
  &= L_{\alpha,j}\spl \palni[j]{i-1},\testfunc[j]{\alpha}\spr
      +  (-1)^{\krond{w}}\spl \Slni[j]{i-1} - \Sln[j]{n-1},\testfunc[j]{\alpha}\spr
    + \tau\spl  f_{\alpha,2}^n ,\testfunc[j]{\alpha} \spr,   
  \label{eqn:LDD-TPR:LDD-scheme:multi-domain:two-phase:ver1}
\end{align}
for $\alpha \in \{n, \nw \}$  along with 
\begin{align}
  \spl \gwli[jk]{i},\testfunc[j]{w}\spr_{\intf{jk}}  &:= \spl -2\lambda_{w}^{jk}
  \pwlni[k]{i-1} -   \gwli[kj]{i-1},\testfunc[j]{w}\spr_{\intf{jk}} 
  && (k\in \ind_j), \label{eqn:LDD-TPR:LDD-scheme:multi-domain:gliupdate:twophase:wetting:weak:ver1} \\[1ex]
  \spl \gnwli[jk]{i},\testfunc[j]{w}\spr_{\intf{jk}} &:= \spl -2\lambda_{\nw}^{jk}
  \pnwlni[k]{i-1} -   \gnwli[kj]{i-1},\testfunc[j]{\nw}\spr_{\intf{jk}} 
  && (k\in \ind_j \cap \indTP), 
  \label{eqn:LDD-TPR:LDD-scheme:multi-domain:gliupdate:twophase:nonwetting:TP:weak:ver1}\\[1ex]
  \spl \gnwli[jk]{i},\testfunc[j]{\nw}\spr_{\intf{jk}} 
       &:= \spl -\gnwli[kj]{i-1},\testfunc[j]{\nw}\spr_{\intf{jk}} 
  && (k\in \indTPRl{j})  
  \label{eqn:LDD-TPR:LDD-scheme:multi-domain:gliupdate:Richards:twophase:nonwetting:R:weak:ver1}
\end{align}
are fulfilled for all test functions $(\varphi_w,\varphi_{\nw})\in \prodFs[00]{w}\times\prodFs[00]{\nw}$.
\end{problem}

\begin{remark}
        We note  that the iterates $(p_{w}^{n,i}, p_{\nw}^{n,i})$ are only required to be in 
        $\prodFs{w}\times\prodFs{\nw}$ and need not to be in $\Fs^w\times\Fs^{\nw}$ (the latter meaning continuity over interfaces). If the family of subsequent solutions to the \LDDTPR solver step, \autoref{env:LDD-TPR:LDD-scheme:multi-domain:weak}, converge to a solution of 
        \autoref{problem:LDD-TPR:semiDiscrete:multi-domain:limitSystem}, 
        then the continuity of the pressures is guaranteed in the limit. 
\end{remark}
From the proof of \autoref{LDD-TPR:convergence:theorem} we expect  that the convergence of the solver holds also in the multi-domain case irrespectively of the choice of initial iterates.
Therefore, it is possible to chose other initial iterates than given in
(\ref{eqn:LDD-TPR:LDDsolver:intialIterate:pressure:ver1}) and (\ref{eqn:LDD-TPR:LDDsolver:intialIterate:gla:ver1}). 
In particular, $\gali[lk]{0}$ can instead be chosen to belong to $H^{1/2}(\interface_{lk})'$ providing $\gali[lk]{i}\in H^{1/2}(\interface_{lk})'$ as well, 
and \autoref{env:LDD-TPR:LDD-scheme:multi-domain:weak:ver1} can be tested with functions $(\varphi_w,\varphi_{\nw})\in \prodFs[]{w}\times\prodFs[]{\nw}$ instead of 
$(\varphi_w,\varphi_{\nw})\in \prodFs[00]{w}\times\prodFs[00]{\nw}$.
While in general it is not clear whether it is possible to  approximate the Neumann fluxes in \autoref{problem:LDD-TPR:semiDiscrete:multi-domain:limitSystem} by functionals in $\Ts[]{lk}'$, in situations, where $\prodFs[]{w}\times\prodFs[]{\nw}$ can be chosen as test function space, 
it is useful to do so for two reasons. 
First, this makes the Lax–Milgram arguments from \autoref{lemma:LDD-TPR:LDDTPRSolverHasSolution} carry over to the multi-domain situation here, so that a solution to each iteration of the solver can be guaranteed. 
Secondly, implementation is facilitated, as the requirement 
$(\varphi_w,\varphi_{\nw})\in \prodFs[00]{w}\times\prodFs[00]{\nw}$
is more difficult to achieve in an implementation than 
$(\varphi_w,\varphi_{\nw})\in \prodFs[]{w}\times\prodFs[]{\nw}$.
Thus,  a more practical formulation of the \LDDTPR solver  used in Section \ref{LDD-TPR:section:NumericalSection} below is given by
    
\begin{problem}[\LDDTPR solver step, multi-domain, version 2] 
\label{env:LDD-TPR:LDD-scheme:multi-domain:weak}
Let functions $(p_w^{n-1}, p_{\nw}^{n-1})$ $\in \Fs^w \times \Fs^{\nw}$ be given and  
define on all subdomains $\dom{l}$, $l\in\ind$, for arbitrary $\nu_l^\alpha \in \Fs_l^\alpha$ and 
$\zeta_{lk}^{\alpha} \in \Ts[]{lk}'$, $k\in\ind_l$, as initial iterates 
\begin{align*}
    \palni{0}:= \nu_l^\alpha,
\end{align*}
where  $\alpha \in \{w\}$ for $l\in \indR$ and $\alpha \in \{w, \nw\}$ for $l\in \indTP$
as well as 
\begin{align*}
    \gali[lk]{0} := \zeta_{lk}^{\alpha}
\end{align*}
in $\Ts[]{lk}'$. 
On interfaces $\interface_{lk}$ between \Richards domains, we take
\begin{align*}
    \gnwli[lk]{0}=0.
\end{align*}
In addition, choose on each domain $\dom{l}$, $l\in\ind$, some real number $L_{\alpha,l} > 0$  and on all interfaces $\intf{lk}$, $k\in \ind_l$, real numbers 
$\lambda_{\alpha}^{lk}=\lambda_{\alpha}^{kl} > 0$.
Given previously known iterates
$(p_{w}^{n,i-1}, p_{\nw}^{n,i-1}) \in \prodFs{w}\times\prodFs{\nw}$, as well as $\gali[lk]{i-1} \in \Ts[]{lk}' $,  $(\NN\in i \geq 1)$, one step of the \LDDTPR solver consists of finding 
$(p_{w}^{n,i}, p_{\nw}^{n,i}) \in \prodFs{w}\times\prodFs{\nw}$
such that on \Richards subdomains, $l \in \indR$, the equations \eqref{eqn:LDD-TPR:LDD-scheme:multi-domain:Richards:ver1}
together with 
\eqref{eqn:LDD-TPR:LDD-scheme:multi-domain:gliupdate:Richards:wetting:weak:ver1}, 
\eqref{eqn:LDD-TPR:LDD-scheme:multi-domain:gliupdate:Richards:nonwetting:weak:ver1}
 are satisfied, and on two-phase domains, $j\in\indTP$, the equations
 \eqref{eqn:LDD-TPR:LDD-scheme:multi-domain:two-phase:ver1}
for $\alpha \in \{n, \nw \}$  along with 
\eqref{eqn:LDD-TPR:LDD-scheme:multi-domain:gliupdate:twophase:wetting:weak:ver1},
\eqref{eqn:LDD-TPR:LDD-scheme:multi-domain:gliupdate:twophase:nonwetting:TP:weak:ver1} and
\eqref{eqn:LDD-TPR:LDD-scheme:multi-domain:gliupdate:Richards:twophase:nonwetting:R:weak:ver1}
are fulfilled for all test functions $(\varphi_w,\varphi_{\nw})\in \prodFs[]{w}\times\prodFs[]{\nw}$.
\end{problem}
\begin{remark}
  Note that the key difference between \autoref{env:LDD-TPR:LDD-scheme:multi-domain:weak:ver1} and \autoref{env:LDD-TPR:LDD-scheme:multi-domain:weak} is the test function space and consequently the space on which the functionals $\gali{i}$ act. 
\end{remark}


\section{Numerical validation of the \LDDTPR solver} 
\label{LDD-TPR:section:NumericalSection}
In this section, we turn to the numerical validation of the \LDDTPR solver for the case $d=2$. 
We provide examples for two different substructurings. For a two-domain case, we  compare the performance of the \LDDTPR solver to  the full two-phase flow model.  In addition, we discuss  the choice  of solver parameters. For a multi-domain example involving an  inner subdomain, we  illustrate the performance as well. 
Both domain partitions are displayed in \autoref{fig:LDD-TPR:NumericalDomain}.

All experiments were implemented using {\Python}  and \fenics' main library \dolfin, cf. 
\cite{Fenics:AlnaesBlechta2015a,
Fenics:LoggWellsEtAl2012a}. %
The code for all examples along with its documentation can be found at \cite{LDDcode}.
For a   detailed description of  the design principles  we refer to \cite{SeusThesis2021}.  Here we restrict ourselves to a listed summary
regarding the grid and the ansatz functions.

 {\scshape Substructuring and meshes.} All subdomains $\dom{l} \subset \dom{} \subset \RR^2$ and triangular  meshes are constructed by the \fenics {} mesh tool \mshr. 
 To ensure that the meshes are matching, submeshes $\submesh{l}$ on each subdomain 
 $\dom{l}$  are always extracted from a global conforming mesh $\mesh$ on  $\dom{}$. This means, mesh vertices and faces always lie on the polygons defining the interfaces, and no facets intersect interfaces.  In this way, neighbouring subdomains share vertices and facets over interfaces.
 \Dolfin was instructed to use {\scshape ParMETIS} as mesh partitioner.
 An example of such a mesh can be seen in Figure \ref{fig:fiveDomain:mesh}). %
 If $d_{\triangle}$ is the diameter (two times the circumradius) of a mesh cell (triangle) $\triangle \in \submesh{l}$, the mesh size $h_l$ on each domain is defined as $h_l = \max\{d_{\triangle}\,|\, \triangle\in\submesh{l}\}$. %
On each subdomain mesh, \fenics' first-order Lagrange finite elements, $\FEMp$,  were used as ansatz spaces 
 $\Fs_{h,l}\subset \Fs_l$.
 
   {\scshape Interfaces terms and communication.}  
 The calculation of the Robin-interface terms across interfaces and the data exchange across interfaces  requires manual assembly of the fluxes involving gradients  of $\FEMp$ functions.  
 The calculation of interface terms is done dof-wise and their communication over interfaces  needs to take into account 
 the different mesh and dof numberings on each subdomain adjacent to a given interface. 
 The calculation of the approximations $\galih[lk]{i}$ of the $\gali[lk]{i}$-terms uses discontinuous 
 Galerkin elements of degree $1$, $\FEMgl$. 
 The reason is twofold. On the one hand, the calculation of 
 $\galih[lk]{0}$ necessitates 
 the assembly of fluxes (since we use the initial iterates of 
 \autoref{env:LDD-TPR:LDD-scheme:multi-domain:weak:ver1}), involving thereby the 
 gradient of a $P_1$ function, hence the need for discontinuous 
 ansatz functions, and on the other hand it seemed desirable to have the same number of degrees of 
 freedom as the pressures that need to be added to these terms.\\ 
 The implementation of the $\galih[lk]{i}$ terms is done in the following way. 
 The \LDD solver, upon entering time step $n$, first assembles $\galih[lk]{0}$:
 On each domain $\dom{l}$, $l\in \ind$, the approximation $\fluxaih{n-1}$ of the flux 
 $\fluxb[n-1]{\alpha}$ is assembled in $\FEMgl \times \FEMgl$ and 
 $\fluxaih{n-1}\cdot \vt{n_{lk}} \in \FEMgl$ is added dof-wise 
 to $\palnh{n-1}|_{\Gamma_{lk}}$ for dofs that lie on facets belonging to the interface 
 $\Gamma_{lk}$, $k \in \ind_l$.
 The resulting dofs of the $\galih{0}$-term are then saved to interface dictionaries for communication. 
 During the $i$th iteration of the \LDD solver on $\dom{l}$, the $\galih[kl]{i-1}$ and $\palnih[k]{i-1}$ 
 dofs of the neighbour $k$ are read from these interface dictionaries and are added – again dof-wise –
 along $\Gamma_{lk}$ to get $\galih[lk]{i} \in \FEMgl$.
 Since the form assembly of \autoref{env:LDD-TPR:LDD-scheme:multi-domain:weak} is done in $\FEMp$,
 the $\galih[lk]{i}$ terms enter the form as projections $\Pi\galih[lk]{i} \in \FEMp$, where 
 $\Pi : \FEMgl \rightarrow \FEMp$ is the projection onto  $\FEMp$.

All appearing linear systems were solved using the Generalised Minimal Residual Method (GMRES) in conjunction with Incomplete LU preconditioning (ILU) as realised in  the \Fenics library.
We will use the following 
\begin{notation}
\label{notation:TPTP:numericalSection}
By $\paln{e}$ we denote manufactured solutions, and  $\paln{e,n}:= \paln{e}(\cdot,\cdot,t_n)$ are their evaluation at time step $t_n$. 
We  have  $\paln{e,n} = \paln{n}$ since manufactured solutions solve the semi-discrete \TPR problem.
     Numerical approximations are denoted with an additional $h$, i.e.\ 
    $\palnh{n}$ is the numerical approximation of $\paln{e,n}\bigl/\paln{n}$ and the symbols $\palnih{i}$ denote the numerical approximation of the iterates $\palni{i}$ of the \LDDTPR solver. Note, that the index $n$ is dropped in this case.  This means that $\palnih{i}$ \emph{always} denotes the iterates in the calculation of the $n$-th time step. 
\end{notation}

\begin{figure}[t]
  \centering
  \begin{subfigure}[b]{.4\textwidth}
  \centering
  \ifdevelopmentVersion
    \input{\LDDTPRfigSource/LDD-TPR-NumericalDomain2Patch.tikz}
  \else
    \includegraphics[width=\textwidth]{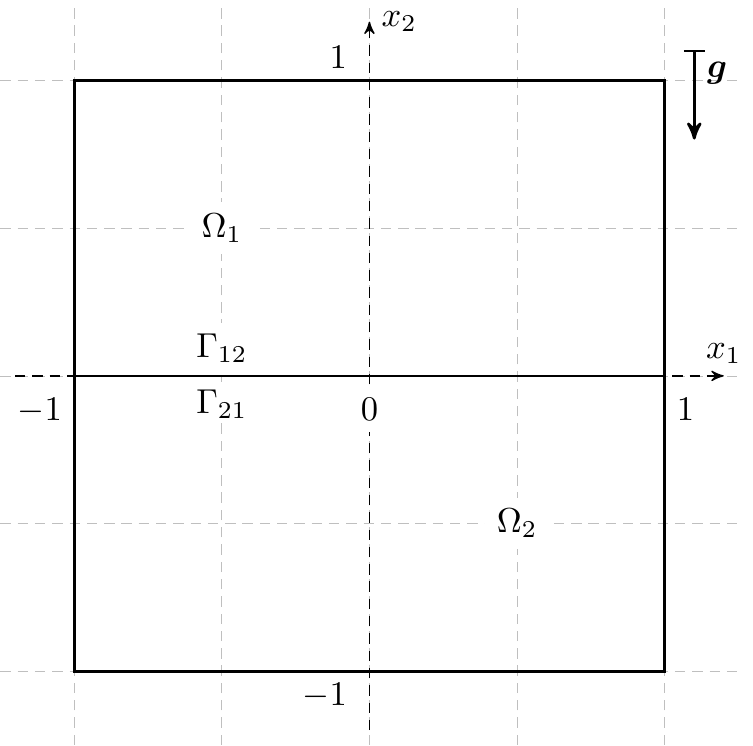}
  \fi
  \caption{Two-domain substructuring. \label{fig:LDD-TPR:NumericalDomain:2patch} }
  \end{subfigure}\hspace{1ex}
  \begin{subfigure}[b]{.4\textwidth}
  \centering
  \ifdevelopmentVersion
    \input{\LDDTPRfigSource/LDD-TPR-NumericalFiveDomain.tikz}
  \else
    \includegraphics[width=\textwidth]{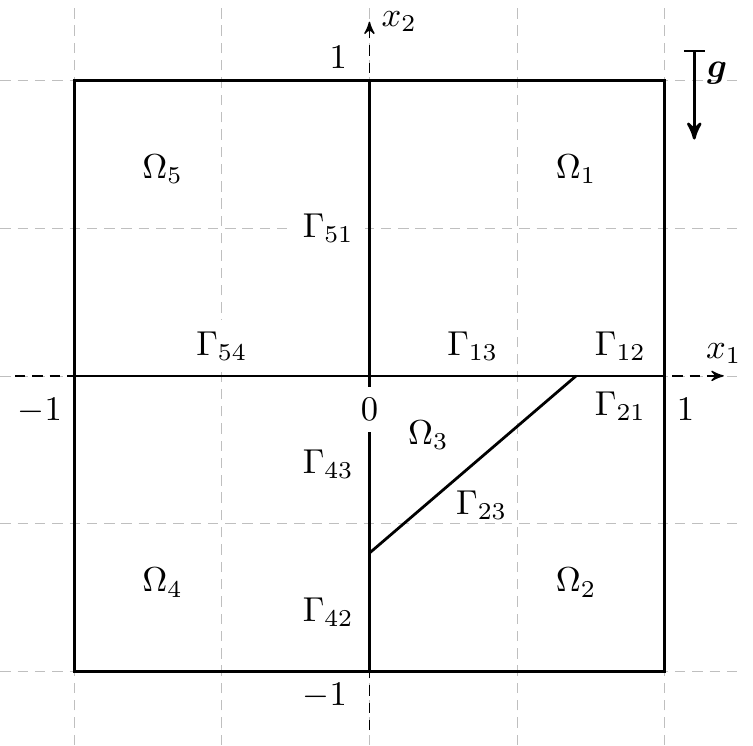}
  \fi
  \caption{Five domains with inner subdomain. \label{fig:LDD-TPR:NumericalDomain:FiveDomainInnerPatch}}
  \end{subfigure}%
  \caption[Numerical domains]{Domains used in the numerical experiments. 
           All domains are polygonal subdivisions of the 
           unit square $[0,1] \times [0,1]$. 
           The nomenclature of the interfaces follows the 
           conventions introduced in \autoref{LDD-TPR:Notation:multidomain}.%
           \label{fig:LDD-TPR:NumericalDomain}}
\end{figure}

\subsection{Two-domain computations}
\label{LDD-TPR:section:NumericalSection:TPRexamples}
We start the numerical validation of the \LDDTPR solver for the  two-domain case  shown in   
\autoref{fig:LDD-TPR:NumericalDomain}(\subref{fig:LDD-TPR:NumericalDomain:2patch}).%

\paragraph{Homogeneous intrinsic permeability and porosity}
We assume the  permeability and porosity in both domains to be
the same and   demonstrate the convergence of the scheme using a manufactured solution.
Modelling the flow of water and air, all soil parameters are listed in Table \ref{table:TPR:2patch:soilParameters}(\subref{table:TPR:2patch:soilParameters:sameIntrinsic}). For the relative permeabilities, $S$-$p_c$ relationships as well as the manufactured solution expressions we refer to
\autoref{LDD-TPR:table:TPR:2patch:sameIntrinsic:exact_solution}.
\begin{table}[t]
  \centering
  \begin{tabular}{c  c  c}
\toprule
    {\scshape Data}
        &  $\dom{1}$
            & $\dom{2}$	\\\midrule
    $\kwln{}(s)$
        & $s^2$
            & $s^3$	\\
    $\knwln{}(s)$
        &$(1-s)^2$
            &$(1-s)^3$	\\
    $S_l(p_c)$
        &
        $\begin{cases}
            \frac{1}{(1 + p_c)^{1/2}} & p_c \geq 0 \\
            \hspace{7.25mm}1           & p_c < 0
        \end{cases} $
            & $\begin{cases}
                \frac{1}{(1 + p_c)^{1/3}} & p_c \geq 0 \\
                \hspace{7.25mm}1          & p_c < 0
              \end{cases}$\\ \addlinespace[0.75ex]
\midrule
%
    $\pwln{e}(x, y, t)$
        &$-7 - \bigl(1+t^2\bigr)\bigl(1 + x^2 + y^2\bigr)$ 
            & $ -7 - \bigl(1+t^2\bigr)\bigl(1 + x^2\bigr)$ \\
    $\pnwln{e}(x,y, t)$
        &  -
            & $\bigl(-2-t(1.1+y + x^2)\bigr)y^2$\\
    \addlinespace[0.5ex]
    \bottomrule
\end{tabular}

  \caption[\LDDTPRcaptionPrefix two-domain, same soil parameters, coefficients and exact solution]{\TPR coupling on two domains: 
             coefficient functions and exact solutions.}
  \label{LDD-TPR:table:TPR:2patch:sameIntrinsic:exact_solution}
\end{table}
\begin{table}[t]
  \centering
  \begin{subtable}[t]{0.47\textwidth}
  \centering
  \begin{tabularx}{.79\textwidth}{ c c c c c} 
\toprule
\multicolumn{1}{r}{\scshape Parameter} & \multicolumn{2}{c}{$\dom{1}$}
    & \multicolumn{2}{c}{$\dom{2}$}	\\\midrule
$\Phi_l$& \multicolumn{2}{c}{\num{0,22}}
    & \multicolumn{2}{c}{\num{0,22}}	\\
$k_{\mathsf{i},l}$& \multicolumn{2}{c}{\num{0,01}}
    & \multicolumn{2}{c}{\num{0,01}}	\\
%
%
   \cmidrule(lr){2-3}\cmidrule(lr){4-5} 
   $\alpha$ & $w$ & $\nw$ & $w$ & $\nw$ \\\cmidrule(lr){2-3}\cmidrule(lr){4-5}
$\mu_{\alpha}[\si{\kilo\gram\per(\meter\second)}]$
    & \num{1} & - & \num{1} & $\frac{1}{50}$
    \\
$\rho_{\alpha}[\si{\kilo\gram\per\cubic\meter}]$
    & \num{997} & - & \num{997} & \num{1,225} 
    \\
\bottomrule
\end{tabularx}

  \caption[\LDDTPRcaptionPrefix two-domain, same soil parameters]{\TPR coupling on two domains:  soil 
          parameters for the case with same intrinsic permeabilities and porosities.}
  \label{table:TPR:2patch:soilParameters:sameIntrinsic}
\end{subtable}\hspace{2ex}
\begin{subtable}[t]{0.47\textwidth}
    \centering
    \begin{tabularx}{.83\textwidth}{ c c c c c} 
\toprule
\multicolumn{1}{r}{\scshape Parameter} & \multicolumn{2}{c}{$\dom{1}$}
    & \multicolumn{2}{c}{$\dom{2}$}	\\\midrule
$\Phi_l$& \multicolumn{2}{c}{\num{0,22}}
    & \multicolumn{2}{c}{\num{0,022}}	\\
$k_{\mathsf{i},l}$& \multicolumn{2}{c}{\num{0,01}}
    & \multicolumn{2}{c}{\num{0,0001}}	\\
%
%
    \cmidrule(lr){2-3}\cmidrule(lr){4-5}
    $\alpha$& $w$ & $\nw$ & $w$ & $\nw$ \\\cmidrule(lr){2-3}\cmidrule(lr){4-5}
$\mu_{\alpha}[\si{\kilo\gram\per(\meter\second)}]$
    & \num{1} & $\frac{1}{50}$ & \num{1} & $\frac{1}{50}$
    \\
$\rho_{\alpha}[\si{\kilo\gram\per\cubic\meter}]$
    & \num{997} & \num{1,225} & \num{997} & \num{1,225} 
    \\
\bottomrule
\end{tabularx}

    \caption[\LDDTPRcaptionPrefix two-domain, varying soil parameters]{%
      \TPR coupling on two domains: case with varying intrinsic
      permeabilities and porosities.%
      }
    \label{table:TPR:2patch:soilParameters:differentIntrinsic}
    \end{subtable}
    \caption[two-domain, soil parameters]{\TPR coupling on two domains:  soil 
          parameters for same (\subref{table:TPR:2patch:soilParameters:sameIntrinsic}) and varying (\subref{table:TPR:2patch:soilParameters:differentIntrinsic}) intrinsic permeabilities and porosities.}
  \label{table:TPR:2patch:soilParameters}
\end{table}

\autoref{LDD-TPR:fig:TPR:2patch:sameIntrinsic} shows the results for a simulation over 
$1500$ time steps of size $\tau = 0.001$ on the time interval $[0,T] = [0,1.5]$ using a mesh size 
$h\approx \num{0.071}$. 
The algorithm  was set to terminate 
after the stopping criterion 
\begin{align*}
  \|\palnih{i} - \palnih{i-1}\|_{{}_2} < \epsilon_s:= 2\cdot\num{e-6}
\end{align*}
had been reached for all appearing $l$ and $\alpha$.
Parameters of the \TPR solver were chosen as $L_{w,1} = \num{0.007}$ and 
$L_{\alpha,2} = \num{0.005}$ for all phases $\alpha\in\{w,\nw\}$ and 
$\lambda_w^{12} = \lambda_{\nw}^{12}= \num{0,75}$. 
\autoref{LDD-TPR:fig:TPR:2patch:sameIntrinsic}(\subref{LDD-TPR:fig:TPR:2patch:sameIntrinsic:errornorm})
shows the relative error norms with respect to the exact solution over time, demonstrating that the accuracy remains invariant 
over time. 
The relative error of the nonwetting phase remains steadily around $\num{0.009}$,  
that of the wetting phases below $\num{0.01}$\%.   %
The nonwetting phase shows a greater approximation error, which is not unexpected, since no nonwetting phase equation is assumed in $\dom{1}$. 
\autoref{LDD-TPR:fig:TPR:2patch:sameIntrinsic}(\subref{LDD-TPR:fig:TPR:2patch:sameIntrinsic:subsequent_errors})
displays the  errors of the solver for the time step $\num{1500}$ at time $T = \num{1.5}$.

\begin{figure}[hp]
  \centering
  \begin{subfigure}[b]{.45\textwidth}
    \centering
    \ifdevelopmentVersion
      \input{\LDDTPRfigSource/2020-07-08-TPR-2P-realistic-same-intrinsic-perm0.01_error_norms.tikz}    
    \else
      \includegraphics[width=\textwidth]{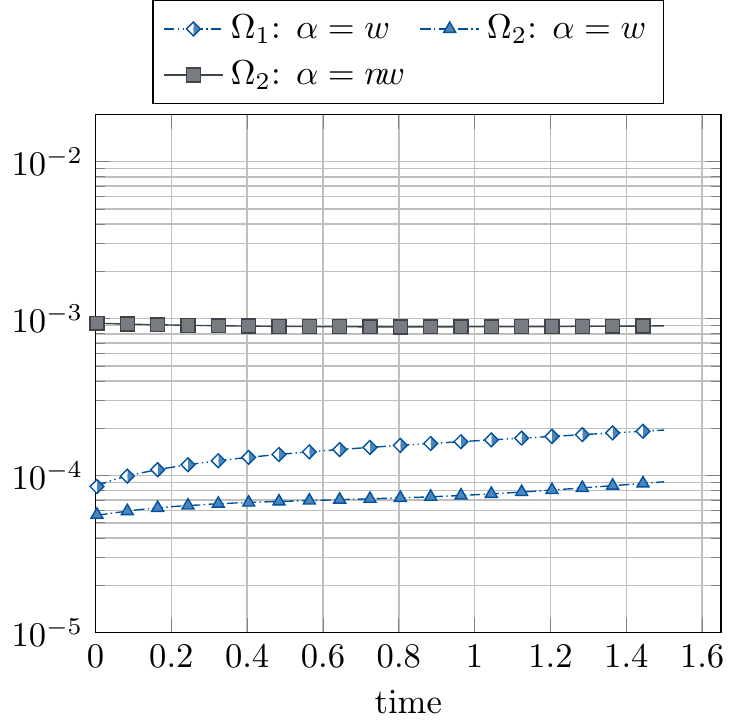}
    \fi
    \caption{Relative error norms $\|\paln{e,n} - \palnh{n}\|_{{}_2}\bigl/\|\paln{e,n}\|_{{}_2}$ over time 
             $t$ for $(\alpha,l) \in \{(w,1), (w,2), (\nw,2)\}$.}
    \label{LDD-TPR:fig:TPR:2patch:sameIntrinsic:errornorm}
  \end{subfigure}\hspace{2ex}
  \begin{subfigure}[b]{.45\textwidth}
    \centering
    \ifdevelopmentVersion
      \input{\LDDTPRfigSource/2020-07-08-TPR-2P-realistic-same-intrinsic-perm0.01_subsequent_iterations.tikz}
    \else
      \includegraphics[width=\textwidth]{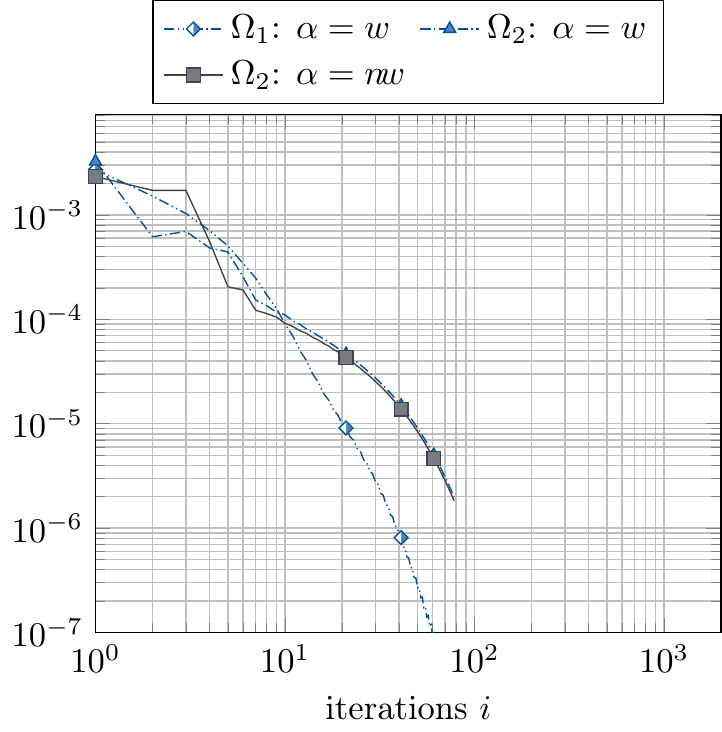}
    \fi
    \caption{Subsequent errors $\|\palnih{i} - \palnih{i-1}\|_{{}_2}$ over iteration $i$ at time step 
             $n=1500$, $t_n = 1.5$. \label{LDD-TPR:fig:TPR:2patch:sameIntrinsic:subsequent_errors} }
  \end{subfigure} \\[4ex]
  \begin{subfigure}[b]{.45\textwidth}
    \centering
    \ifdevelopmentVersion
      \input{\LDDTPTPfigSource/2020-07-21-TP-TP-2patch-realistic-same-intrinsic-nw-zero_errornorm.tikz}    
    \else
      \includegraphics[width=\textwidth]{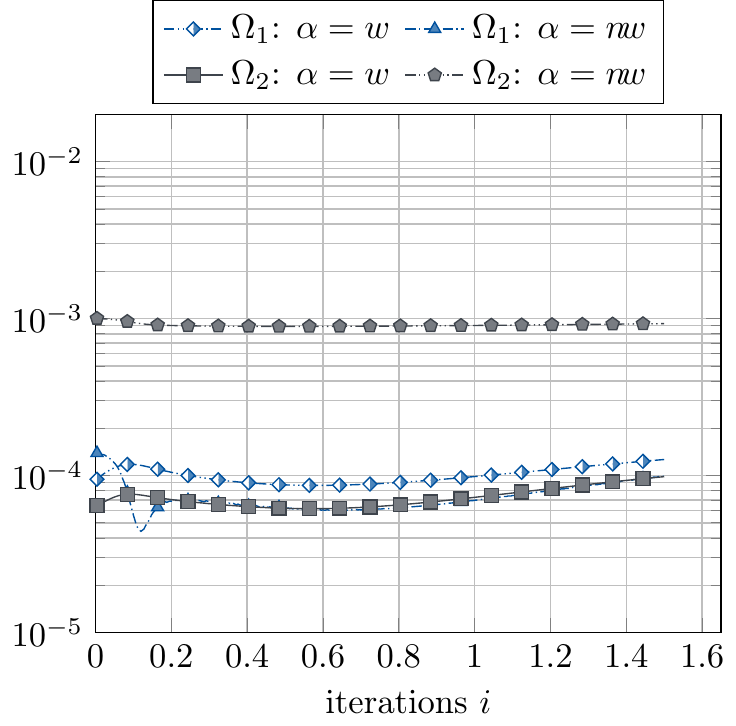}
    \fi
    \caption{Relative error norms $\|\paln{e,n} - \palnh{n}\|_{{}_2}\bigl/\|\paln{e,n}\|_{{}_2}$ over time 
             $t$ for $\alpha\in\{w,\nw\}$ and $l\in\{1,2\}$.}
    \label{LDD-TPR:fig:TPTP:2patch:sameIntrinsic:ZeroNonwetting:errornorm}
  \end{subfigure}\hspace{2ex}
  \begin{subfigure}[b]{.45\textwidth}
    \centering
    \ifdevelopmentVersion
      \input{\LDDTPTPfigSource/2020-07-21-TP-TP-2patch-realistic-same-intrinsic-nw-zero_subsequent_errors.tikz}
    \else
      \includegraphics[width=\textwidth]{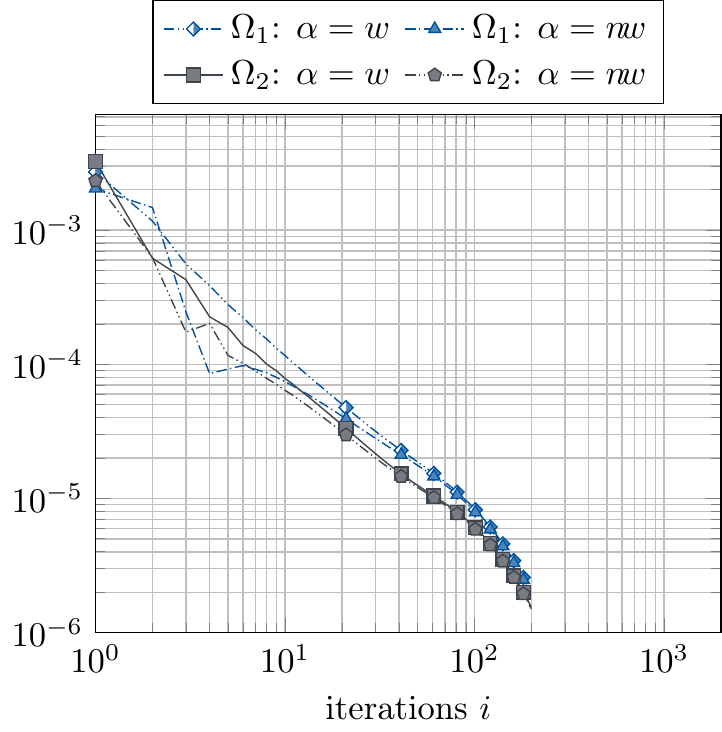}
    \fi
    \caption{Subsequent errors $\|\palnih{i} - \palnih{i-1}\|_{{}_2}$ over iteration $i$ at time step 
             $n=1500$, $t_n = 1.5$. \label{LDD-TPR:fig:TPTP:2patch:sameIntrinsic:ZeroNonwetting:subsequent_errors} }
  \end{subfigure}  
  \caption[\LDDTPRcaptionPrefix two-domain, same soil parameters]{\TPR coupling on two domains: 
          Relative error norms (\subref{LDD-TPR:fig:TPR:2patch:sameIntrinsic:errornorm})
          and subsequent errors at a fixed time step 
          (\subref{LDD-TPR:fig:TPR:2patch:sameIntrinsic:subsequent_errors})
          for a simulation over $1500$ time steps with same intrinsic 
          permeabilities and porosities,  and parameters 
          $h\approx 0.071$, $\tau = 1\cdot 10^{-3}$, $L_{w,1} = \num{0.007}$ and 
          $L_{\alpha,2} = \num{0.005}$, $\alpha\in\{w,\nw\}$ and 
          $\lambda_w^{12} = \lambda_{\nw}^{12}= \num{0,75}$.
          Relative error norms 
          (\subref{LDD-TPR:fig:TPTP:2patch:sameIntrinsic:ZeroNonwetting:errornorm})
          and subsequent errors at a fixed time step 
          (\subref{LDD-TPR:fig:TPTP:2patch:sameIntrinsic:ZeroNonwetting:subsequent_errors}) 
          for a simulation of the same situation using the  same parameters, but assuming a \TPTP coupling and assumed zero  nonwetting phase. 
          }
  \label{LDD-TPR:fig:TPR:2patch:sameIntrinsic}
\end{figure}

To determine how the use of the \TPR coupling in this situation affects both accuracy and performance, 
we compare  with a simulation of the same setting, assuming constant nonwetting pressure, 
$\pnwln[1]{} \equiv 0$ on $\dom{1}$ and use 
the LDD solver  for two-phase flow
equations in  $\dom{1}$ and $\dom{2}$ (\LDDTPTP solver, see \cite{SeusEnumath2019} for details). 
As 
\autoref{LDD-TPR:fig:TPR:2patch:sameIntrinsic}(\subref{LDD-TPR:fig:TPTP:2patch:sameIntrinsic:ZeroNonwetting:errornorm})
shows, the same precision is achieved in both cases, using either the \LDDTPR or the \LDDTPTP solver. 
The worst relative error can be observed for the nonwetting phase on $\dom{2}$, 
similarly to the case of the \TPR coupling shown in  
\autoref{LDD-TPR:fig:TPR:2patch:sameIntrinsic}(\subref{LDD-TPR:fig:TPR:2patch:sameIntrinsic:errornorm}). 
This suggests that the error is not dominated 
by the use of the \TPR coupling in place of the a complete \TPTP coupling.

Naturally, the \LDDTPTP solver is slower, having to solve an additional system. 
The subsequent errors at a fixed time step, 
\autoref{LDD-TPR:fig:TPR:2patch:sameIntrinsic}(\subref{LDD-TPR:fig:TPTP:2patch:sameIntrinsic:ZeroNonwetting:subsequent_errors}), show
in addition, 
that the \LDDTPTP solver needs $199$ iterations in the $1500$th time step to achieve the same stopping criterion. $29$ iterations were needed in the first time step. 

These results show that in situations in which the assumptions for the validity of the \Richards equation  hold, the hybrid  \LDDTPR solver excels  over the \LDDTPTP solver as there is a noticeable performance gain at virtually no loss of approximation accuracy. 

\paragraph{Heterogeneous intrinsic permeabilities and porosities}
We investigate numerically the influence of heterogeneneous soil parameters  running a test case 
with varying intrinsic permeabilities and porosities. 
The values used are listed in 
\ref{table:TPR:2patch:soilParameters}(\subref{table:TPR:2patch:soilParameters:differentIntrinsic}).
Relative permeabilities, $p_c$–$S$ relationships and the exact solutions are the same as before, cf. 
\autoref{LDD-TPR:table:TPR:2patch:sameIntrinsic:exact_solution}.
Grid parameters remain the same, namely $h \approx 0.071$,
and $\tau = \num{0,001}$ for the time step. 
\autoref{LDD-TPR:fig:TPR:2patch:differentIntrinsic} shows results for a simulation comprising 
 $1500$ time steps using \LDDTPR parameters
$\lambda_\alpha^{12} = \num{0,5}$, $L_{w,1}= \num{0,007}$ and  
$L_{\alpha,2}= \num{0,0005}$, for $\alpha\in\{w, \nw\}$.
The stopping criterion was set to $\epsilon_s = 2\cdot\num{e-6}$.
\begin{figure}[t]
  \centering
  \begin{subfigure}[b]{.45\textwidth}
    \centering
    \ifdevelopmentVersion
      \input{\LDDTPRfigSource/2020-07-08-TPR-2-patch-realistic-different-intrinsic_errornorms.tikz}    
    \else
      \includegraphics[width=\textwidth]{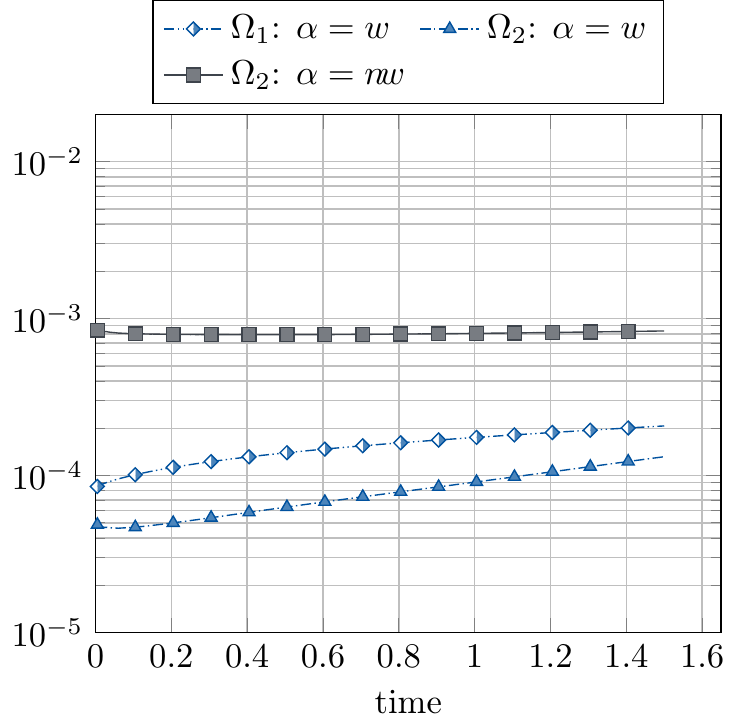}
    \fi
    \caption{Relative error norms $\|\paln{e,n} - \palnh{n}\|_{{}_2}\bigl/\|\paln{e,n}\|_{{}_2}$ over time 
             $t$ for $(\alpha,l) \in \{(w,1), (w,2), (\nw,2)\}$.}
    \label{LDD-TPR:fig:TPR:2patch:differentIntrinsic:errornorm}
  \end{subfigure}\hspace{2ex}
  \begin{subfigure}[b]{.45\textwidth}
    \centering
    \ifdevelopmentVersion
      \input{\LDDTPRfigSource/2020-07-08-TPR-2-patch-realistic-different-intrinsic_subsequent_iterations.tikz}
    \else
      \includegraphics[width=\textwidth]{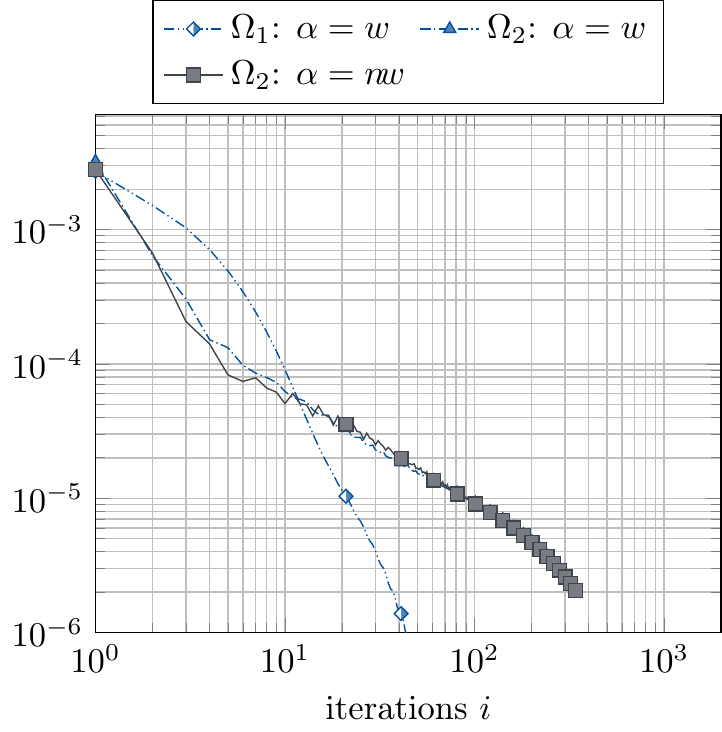}
    \fi
    \caption{Subsequent errors $\|\palnih{i} - \palnih{i-1}\|_{{}_2}$ over iteration $i$ at time step 
             $n=1500$, $t_n = 1.5$. \label{LDD-TPR:fig:TPR:2patch:differentIntrinsic:subsequent_errors} }
  \end{subfigure}
  \caption[\LDDTPRcaptionPrefix two-domain, varying soil parameters]{\TPR coupling on two domains: 
          Relative error norms (\subref{LDD-TPR:fig:TPR:2patch:differentIntrinsic:errornorm})
          and subsequent errors at a fixed time step 
          (\subref{LDD-TPR:fig:TPR:2patch:differentIntrinsic:subsequent_errors})
          for a simulation over $1500$ time steps with varying intrinsic 
          permeabilities and porosities,  \emph{excluding gravity} and parameters 
          $h\approx 0.071$, $\tau = 1\cdot 10^{-3}$, $\lambda_\alpha^{12} = \num{0,5}$, 
          $L_{w,1}= \num{0,007}$ and $L_{\alpha,2}= \num{0,0005}$, for $\alpha\in\{w, \nw\}$.
          }
  \label{LDD-TPR:fig:TPR:2patch:differentIntrinsic}
\end{figure}
As can be seen from 
\autoref{LDD-TPR:fig:TPR:2patch:differentIntrinsic}(\subref{LDD-TPR:fig:TPR:2patch:differentIntrinsic:errornorm})
the final approximation precision is unaffected by the more challenging soil parameters compared to
the case with same intrinsic permeabilities. 
Mind the adjusted \LDDTPR parameters $L_{\alpha,l}$ and
$\lambda_\alpha^{12}$, however. 
In contrast to the previously shown case more iterations are needed to achieve the stopping criterion precision as is visible in \autoref{LDD-TPR:fig:TPR:2patch:differentIntrinsic}(\subref{LDD-TPR:fig:TPR:2patch:differentIntrinsic:subsequent_errors}). 
Required iterations ranged from $63$ iterations in the first time step to $348$ iterations in time step 
$1500$ shown in 
\autoref{LDD-TPR:fig:TPR:2patch:differentIntrinsic}(\subref{LDD-TPR:fig:TPR:2patch:differentIntrinsic:subsequent_errors}).

\paragraph{Comparison to coarser time step size.} 
\autoref{LDD-TPR:fig:TPR:2patch:differentIntrinsic:dt0.01}(\subref{LDD-TPR:fig:TPR:2patch:differentIntrinsic:dt0.01:errornorm})
and
\autoref{LDD-TPR:fig:TPR:2patch:differentIntrinsic:dt0.01}(\subref{LDD-TPR:fig:TPR:2patch:differentIntrinsic:dt0.01:subsequent_errors})
show the same situation but simulated with a coarser time step $\tau = 1\cdot\num{e-2}$ reaching 
$T = 15$ after $\num{1500}$ iterations.
Interestingly, the error norms of all phases are in the same range as for the simulation with the 
finer time step, cf. \autoref{LDD-TPR:fig:TPR:2patch:differentIntrinsic:dt0.01}(\subref{LDD-TPR:fig:TPR:2patch:differentIntrinsic:dt0.01:errornorm}), 
albeit the errors of the wetting phases behaving noticeably worse. The increase of error that starts 
taking place around $T=10$ is due to solver maxing out the maximal iterations number, $1000$. 
\begin{figure}[hp]
    \begin{subfigure}[b]{.45\textwidth}
    \centering
    \ifdevelopmentVersion
      \input{\LDDTPRfigSource/2020-07-08-TPR-2-patch-realistic-different-intrinsic_dt0.01_errornorms.tikz}    
    \else
      \includegraphics[width=\textwidth]{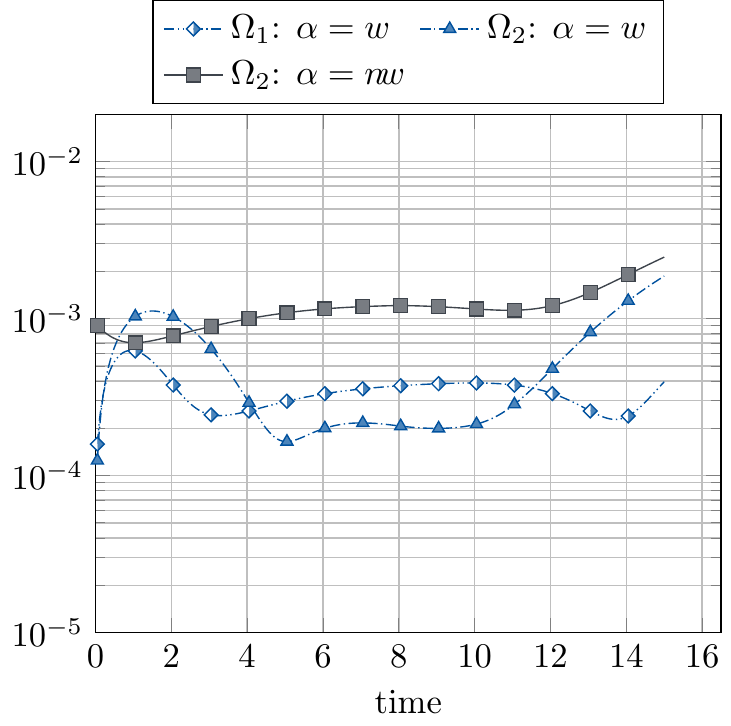}
    \fi
    \caption{Relative error norms $\|\paln{e,n} - \palnh{n}\|_{{}_2}\bigl/\|\paln{e,n}\|_{{}_2}$ over time 
             $t$ for $(\alpha,l) \in \{(w,1), (w,2), (\nw,2)\}$.}
    \label{LDD-TPR:fig:TPR:2patch:differentIntrinsic:dt0.01:errornorm}
  \end{subfigure}\hspace{2ex}
  \begin{subfigure}[b]{.45\textwidth}
    \centering
    \ifdevelopmentVersion
      \input{\LDDTPRfigSource/2020-07-08-TPR-2-patch-realistic-different-intrinsic_dt0.01_subsequent_iterations.tikz}
    \else
      \includegraphics[width=\textwidth]{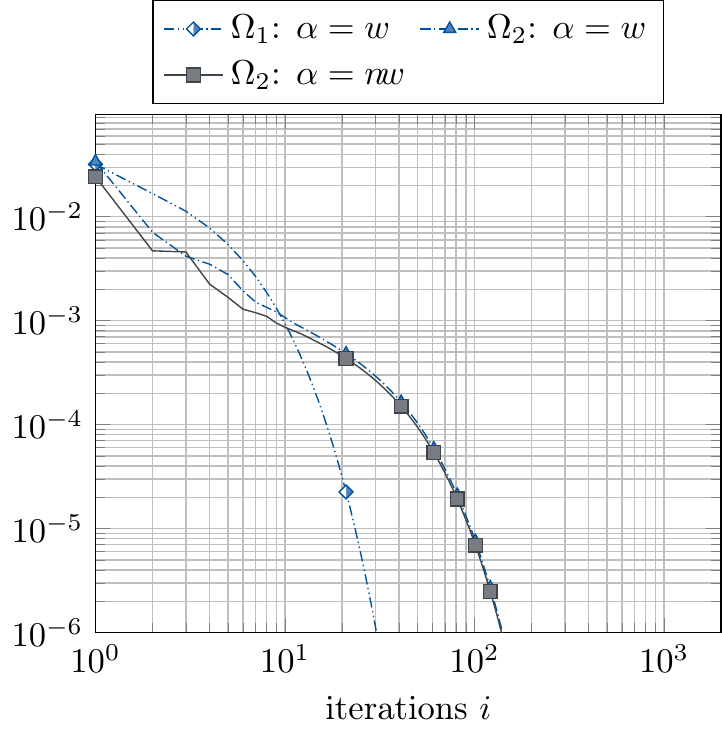}
    \fi
    \caption{Subsequent errors $\|\palnih{i} - \palnih{i-1}\|_{{}_2}$ over iteration $i$ at time step 
             $n=167$, $t_n = 1.67$. \label{LDD-TPR:fig:TPR:2patch:differentIntrinsic:dt0.01:subsequent_errors} }
  \end{subfigure}\\[2ex]
  \centering
  \begin{subfigure}[b]{.45\textwidth}
    \centering
    \ifdevelopmentVersion
      \input{\LDDTPRfigSource/2020-07-08-TPR-2-patch-realistic-different-intrinsic_dt0.01_subsequent_iteration1500.tikz}
    \else
      \includegraphics[width=\textwidth]{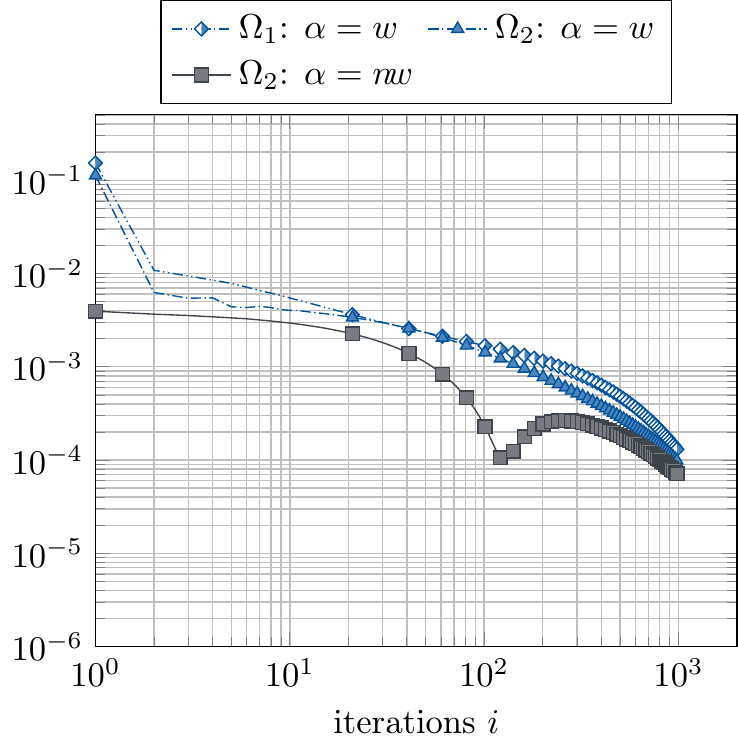}
    \fi
    \caption{Subsequent errors $\|\palnih{i} - \palnih{i-1}\|_{{}_2}$ over iteration $i$ at time step $n=1500$, $t_n = 15$. 
    \label{LDD-TPR:fig:TPR:2patch:differentIntrinsic:dt0.01:subsequent_errors:1500} }
  \end{subfigure}
  \caption[\LDDTPRcaptionPrefix two-domain, varying soil parameters, coarse $\tau$]{\TPR coupling on two domains: 
      Relative error norms (\subref{LDD-TPR:fig:TPR:2patch:differentIntrinsic:dt0.01:errornorm})
          and subsequent errors at a fixed time step 
          (\subref{LDD-TPR:fig:TPR:2patch:differentIntrinsic:dt0.01:subsequent_errors})
          and (\subref{LDD-TPR:fig:TPR:2patch:differentIntrinsic:dt0.01:subsequent_errors:1500})
          for the same situation as in \autoref{LDD-TPR:fig:TPR:2patch:differentIntrinsic} but 
          simulated using a \emph{coarser time step} 
          $\tau = 1\cdot 10^{-2}$.                                                                          
      }
  \label{LDD-TPR:fig:TPR:2patch:differentIntrinsic:dt0.01}
\end{figure}
After $t=10$ the solver always iterates $1000$ times but fails to reach the stopping criterion. 
The effect is shown in \autoref{LDD-TPR:fig:TPR:2patch:differentIntrinsic:dt0.01:subsequent_errors}(\subref{LDD-TPR:fig:TPR:2patch:differentIntrinsic:dt0.01:subsequent_errors:1500})
depicting the $1000$ iterations in time step $t_{1500}$ that the solver uses without reducing the subsequent errors sufficiently.
Notably, convergence is very slow.

To compare the behaviour of the solver at a similar time than is depicted in 
\autoref{LDD-TPR:fig:TPR:2patch:differentIntrinsic}(\subref{LDD-TPR:fig:TPR:2patch:differentIntrinsic:subsequent_errors}), 
\autoref{LDD-TPR:fig:TPR:2patch:differentIntrinsic:dt0.01}(\subref{LDD-TPR:fig:TPR:2patch:differentIntrinsic:dt0.01:subsequent_errors})
shows the behaviour of the solver at $t=\num{1,67}$. %
used to achieve the error tolerance $\epsilon_s = 1\cdot\num{e-6}$.
This means that up to this point in time, the \LDDTPR solver needs less iterations in each time step of the simulation using $\tau=\num{0,01}$ than in the example using $\tau=\num{0,001}$, all the while achieving the same level of approximation error!

\paragraph{Influence of \LDDTPR parameters}
\begin{figure}[t]
  \centering
  \begin{subfigure}[b]{.45\textwidth}
    \centering
    \ifdevelopmentVersion
      \input{\LDDTPRfigSource/2020-06-17-TP-R-2-patch-realistic-no-gravity-but-varying-intrinsic-perm_bad_params_errornorms.tikz}    
    \else
      \includegraphics[width=\textwidth]{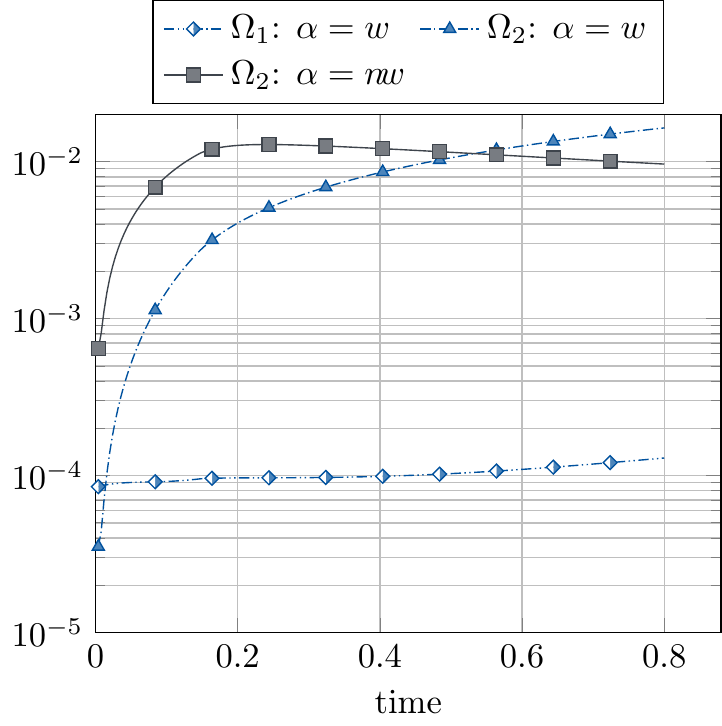}
    \fi
    \caption{Relative error norms $\|\paln{e,n} - \palnh{n}\|_{{}_2}\bigl/\|\paln{e,n}\|_{{}_2}$ over time 
             $t$ for $(\alpha,l) \in \{(w,1), (w,2), (\nw,2)\}$.}
    \label{LDD-TPR:fig:TPR:2patch:differentIntrinsic:badParams:errornorm}
  \end{subfigure}\hspace{2ex}
  \begin{subfigure}[b]{.45\textwidth}
    \centering
    \ifdevelopmentVersion
      \input{\LDDTPRfigSource/2020-06-17-TP-R-2-patch-realistic-no-gravity-but-varying-intrinsic-perm_bad_params_subsequent_iterations.tikz}
    \else
      \includegraphics[width=\textwidth]{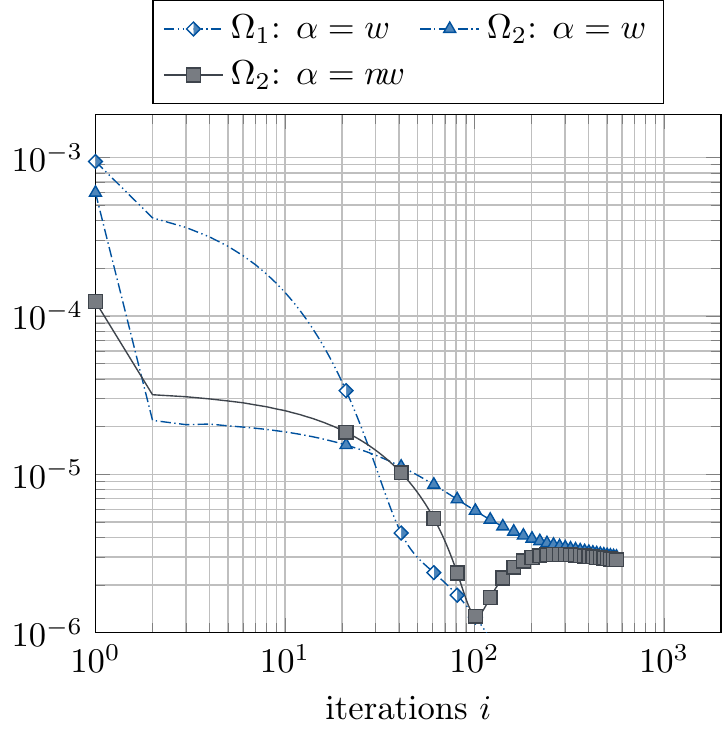}
    \fi
    \caption{Subsequent errors $\|\palnih{i} - \palnih{i-1}\|_{{}_2}$ over iteration $i$ at time step 
             $n=800$, $t_n = 0.8$. \label{LDD-TPR:fig:TPR:2patch:differentIntrinsic:badParams:subsequent_errors} }
  \end{subfigure}
  \caption[\LDDTPRcaptionPrefix two-domain, varying soil parameters, bad parameter choice]{\TPR coupling on two domains: 
          Relative error norms (\subref{LDD-TPR:fig:TPR:2patch:differentIntrinsic:badParams:errornorm})
          and subsequent errors at a fixed time step 
          (\subref{LDD-TPR:fig:TPR:2patch:differentIntrinsic:badParams:subsequent_errors})
          for a simulation over $800$ time steps with varying intrinsic 
          permeabilities and porosities, and parameters 
          $h\approx 0.071$, $\tau = 1\cdot 10^{-3}$, $\lambda_\alpha^{12} = \num{4}$, 
          $L_{w,1}= \num{0,025}$ and $L_{w,2}= \num{0,05}$, $L_{\nw,2}= \num{0,025}$.
          }
  \label{LDD-TPR:fig:TPR:2patch:differentIntrinsic:badParams}
\end{figure}
The \LDD solver is sensitive to the numerical parameters. 
To illustrate this dependence, we revisit the previous example with varying permeabilities and 
porosities keeping the grid parameters 
the same, %
but varying the \LDDTPR parameters:
\autoref{LDD-TPR:fig:TPR:2patch:differentIntrinsic:badParams} shows results of a simulation 
of $800$ time steps using 
$\lambda_\alpha^{12} = \num{4}$, $L_{w,1}= \num{0,025}$ and  
$L_{w,2}= \num{0,05}$, $L_{\nw,2}= \num{0,025}$  as well as  $\epsilon_s = 3\cdot\num{e-6}$
as the stopping criterion. 
While the wetting phase error on $\dom{1}$ at around 
$\num{0,01}$\%  compares  to the one in  
\autoref{LDD-TPR:fig:TPR:2patch:differentIntrinsic}(\subref{LDD-TPR:fig:TPR:2patch:differentIntrinsic:errornorm}),  
\autoref{LDD-TPR:fig:TPR:2patch:differentIntrinsic:badParams}(\subref{LDD-TPR:fig:TPR:2patch:differentIntrinsic:badParams:errornorm})
shows that the errors of both phases on $\dom{2}$ are an order of magnitude worse than 
what has been shown in
\autoref{LDD-TPR:fig:TPR:2patch:differentIntrinsic}(\subref{LDD-TPR:fig:TPR:2patch:differentIntrinsic:errornorm}). 
Accordingly, 
\autoref{LDD-TPR:fig:TPR:2patch:differentIntrinsic:badParams}(\subref{LDD-TPR:fig:TPR:2patch:differentIntrinsic:badParams:subsequent_errors})
indicates by the high number of required iterations as well as the tilts observable in 
the subsequent error plots of the phases on $\dom{2}$, 
that the solver is struggles to find the solution. 
The required iterations to achieve the stopping criterion in this case ranged from $40$ in the 
first time step to $577$ in the $800$th time step depicted in 
\autoref{LDD-TPR:fig:TPR:2patch:differentIntrinsic:badParams}(\subref{LDD-TPR:fig:TPR:2patch:differentIntrinsic:badParams:subsequent_errors}).

\subsection{Multi-domain computations}
\label{section:LDD-TPR:numericalSection:multi-domain}
\begin{figure}[t]
  \centering
  \begin{subfigure}[t]{.45\textwidth}
  \centering
  \ifdevelopmentVersion
    \input{\LDDTPRfigSource/LDD-TPR-NumericalFiveDomainmodels.tikz}
  \else
    \includegraphics[width=\textwidth]{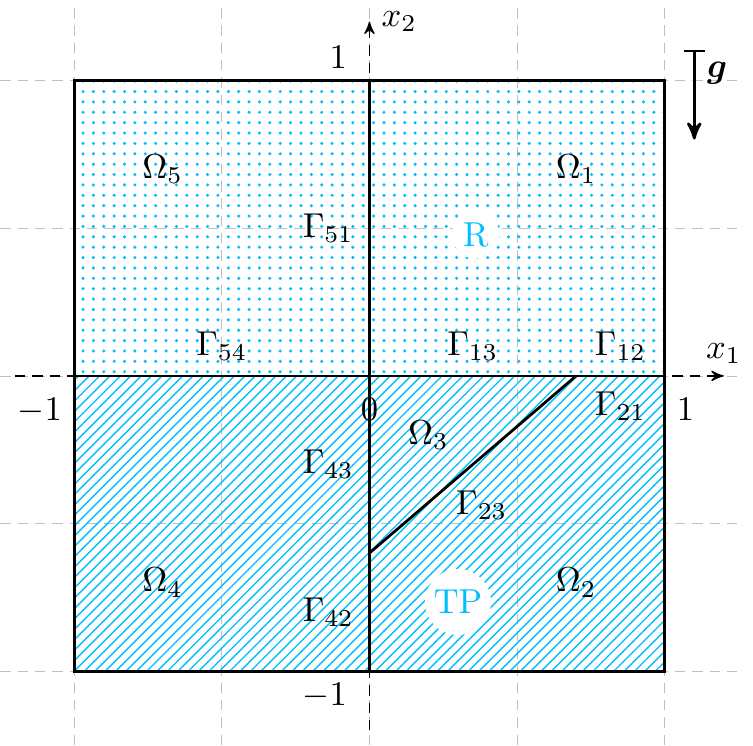}
  \fi
  \caption{Five domains with inner subdomain and highlighted model change. \label{figure:LDD-TPR:TPR:numericalDomains:modelHighlighted:FiveDomainInnerPatch} }
  \end{subfigure}\hspace{2ex}
  \begin{subfigure}[t]{0.45\textwidth}
    \centering
    \ifdevelopmentVersion
      \input{\LDDTPRfigSource/LDD-TPR-NumericalFiveDomainMesh.tikz}    
    \else
      \includegraphics[width=\textwidth]{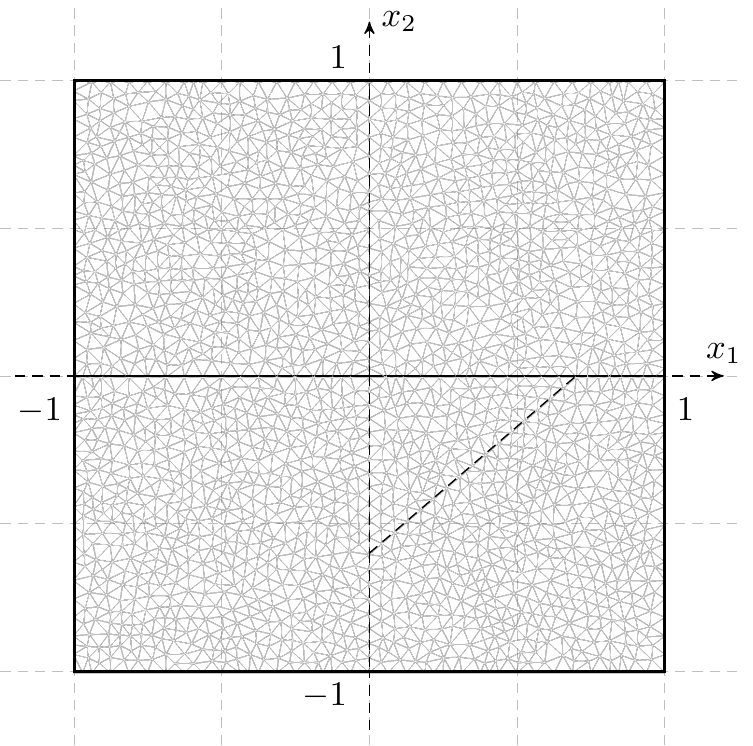}
    \fi
    \caption[Unstructured mesh for five domain substructuring]{Unstructured mesh for the five domain substructuring at \meshrez\ $= 32$.
    \label{fig:fiveDomain:mesh}}
  \end{subfigure}
  \caption[Numerical domain with highlighted models]{\TPR coupling on a five domain substructuring. Highlighted are the areas where different models are being used (\subref{figure:LDD-TPR:TPR:numericalDomains:modelHighlighted:FiveDomainInnerPatch}). The dotted 
          areas features \Richards equations and the striped areas the two-phase flow model.
          Unstructured triangular mesh (\subref{fig:fiveDomain:mesh}) for \meshrez\ $= 32$.
          \label{figure:LDD-TPR:TPR:numericalDomains:modelHighlighted}}
\end{figure}
\begin{table}[t]
  \centering
  \begin{tabular}{c  c  c}
\toprule
        & \Richards 
            & two-phase \\
    {\scshape Data}
        &   $\dom{1}$, $\dom{5}$
            &  $\dom{2}$-$\dom{4}$	\\\midrule
    $\kwln{}(s)$
        & $s^2$
            & $s^3$	\\
    $\knwln{}(s)$
        &$(1-s)^2$
            &$(1-s)^3$	\\
    $S_l(p_c)$
        &
        $\begin{cases}
            \frac{1}{(1 + p_c)^{1/2}} & p_c \geq 0 \\
            \hspace{7.25mm}1           & p_c < 0
        \end{cases} $
            & $\begin{cases}
                \frac{1}{(1 + p_c)^{1/3}} & p_c \geq 0 \\
                \hspace{7.25mm}1          & p_c < 0
              \end{cases}$\\ \addlinespace[0.75ex]
\midrule
%
    $\pwln{e}(x, y, t)$
        &$-7 - \bigl(1+t^2\bigr)\bigl(1 + x^2 + y^2\bigr)$ 
            & $ -7 - \bigl(1+t^2\bigr)\bigl(1 + x^2\bigr)$ \\
    $\pnwln{e}(x,y, t)$
        &  -
            & $\bigl(-3-t(1+y + x^2) - t^2\bigr)y^2$\\
    \addlinespace[0.5ex]
    \bottomrule
\end{tabular}

  \caption[\LDDTPRcaptionPrefix five-domain with inner subdomain: coefficient functions and manufactured solutions]{\TPR coupling on a five-domain with inner subdomain: 
            assumed coefficient functions and manufactured solutions.}
  \label{LDD-TPR:table:TPR:fiveDomainIP:sameIntrinsic:exact_solution}
\end{table}
We advance from  the two-domain examples to  a multi-domain example featuring an inner subdomain, see 
Figure 
\autoref{fig:LDD-TPR:NumericalDomain}(\subref{fig:LDD-TPR:NumericalDomain:FiveDomainInnerPatch}).
We assume the \Richards equation on subdomains $1$, $5$ and the full two-phase flow model on subdomains
$2$–$4$. The inner subdomain is $\dom{3}$. 
According to  \autoref{LDD-TPR:Notation:multidomain}, 
this means $\ind = \{1,2,3,4,5\}$, $\indR = \{1,5\}$ and $\indTP = \{2,3,4\}$, cf. illustration in \autoref{figure:LDD-TPR:TPR:numericalDomains:modelHighlighted}(\subref{figure:LDD-TPR:TPR:numericalDomains:modelHighlighted:FiveDomainInnerPatch}).
\paragraph{Excluding gravity}

We first use an example {excluding gravity} featuring the manufactured solutions, relative permeabilities and $p_c$–$S$-relations 
given in \autoref{LDD-TPR:table:TPR:fiveDomainIP:sameIntrinsic:exact_solution}. Soil parameters are listed in \autoref{LDD-TPR:table:TPR:fiveDomainIP:sameIntrinsic:soilParameters}.
\begin{table}[t]
  \centering
  \begin{tabular}{ c c c c c c c} 
\toprule
        & \multicolumn{2}{c}{\Richards} 
            & \multicolumn{4}{c}{two-phase} \\
\multicolumn{1}{r}{\scshape Parameter} 
    & \multicolumn{2}{c}{$\dom{1}$, $\dom{5}$}
        & \multicolumn{2}{c}{$\dom{3}$} 
            & \multicolumn{2}{c}{$\dom{2}$, $\dom{4}$}	\\\midrule
$\Phi_l$
    & \multicolumn{2}{c}{\num{0,2}}
        & \multicolumn{2}{c}{\num{0,2}} 
            & \multicolumn{2}{c}{\num{0,2}}	\\
$k_{\mathsf{i},l}$
    & \multicolumn{2}{c}{\num{0,01}}
        & \multicolumn{2}{c}{\num{0,01}} 
            & \multicolumn{2}{c}{\num{0,01}}	\\
%
$h_{l}$
    & \multicolumn{2}{c}{\num{0.071}} 
        & \multicolumn{2}{c}{\num{0.070}} 
            & \multicolumn{2}{c}{\num{0.071}} \\\cmidrule(lr){2-3} \cmidrule(lr){4-5} \cmidrule(lr){6-7}
  $\alpha$  & $w$ & $\nw$ 
        & $w$ & $\nw$ 
            & $w$ & $\nw$ \\\cmidrule(lr){2-3} \cmidrule(lr){4-5} \cmidrule(lr){6-7}
$\mu_{\alpha}[\si{\kilo\gram\per(\meter\second)}]$
    & \num{1} & - 
        & \num{1} & $\frac{1}{50}$ 
            & \num{1} & $\frac{1}{50}$
    \\
$\rho_{\alpha}[\si{\kilo\gram\per\cubic\meter}]$
    & \num{997} & - 
        & \num{997} & \num{1,225} 
            & \num{997} & \num{1,225}  
    \\
\bottomrule
\end{tabular}

  \caption[\LDDTPRcaptionPrefix five-domain with inner subdomain, same soil parameters]{\TPR coupling on five-domain with inner subdomain: assumed soil 
          parameters for case with same intrinsic permeabilities and porosities.}
  \label{LDD-TPR:table:TPR:fiveDomainIP:sameIntrinsic:soilParameters}
\end{table}
Notice that the same porosity $\porosity{l} = \num{0,2}$ and intrinsic permeability $k_{i,l} = \num{0,01}$ is assumed on all subdomains, $l\in\ind$.

\autoref{LDD-TPR:fig:TPR:fiveDomainIP:sameIntrinsic}(\subref{LDD-TPR:fig:TPR:fiveDomainIP:sameIntrinsic:errornorm})
shows the relative error norms over time for a simulation of $1000$ time steps of size 
$\tau=1\cdot\num{e-3}$. 
The \LDDTPR parameters were set to 
$L_{w,l} = \num{0,01}$, $L_{\nw,l} = \num{0,004}$ and $\lambda_w^{lk} = 1$ as well as 
$\lambda_{\nw}^{lk} = 0.25$ for all $l\in \ind$ and $k\in \ind_k$.
The error of the nonwetting phase of the inner subdomain $\dom{3}$, which is the worst of all phases 
and subdomains stays consistently below $5\cdot\num{e-3}$. 
The nonwetting phases of $\dom{2}$ and 
$\dom{4}$ don't surpass $1\cdot\num{e-3}$ and all wetting phases stay below $2\cdot\num{e-4}$. 
\begin{figure}[t]
  \centering
  \begin{subfigure}[b]{.45\textwidth}
    \centering
    \ifdevelopmentVersion
      \input{\LDDTPRfigSource/2020-07-16-TPR-FiveDomainIP-same-intrinsic_error_norms.tikz}    
    \else
      \includegraphics[width=\textwidth]{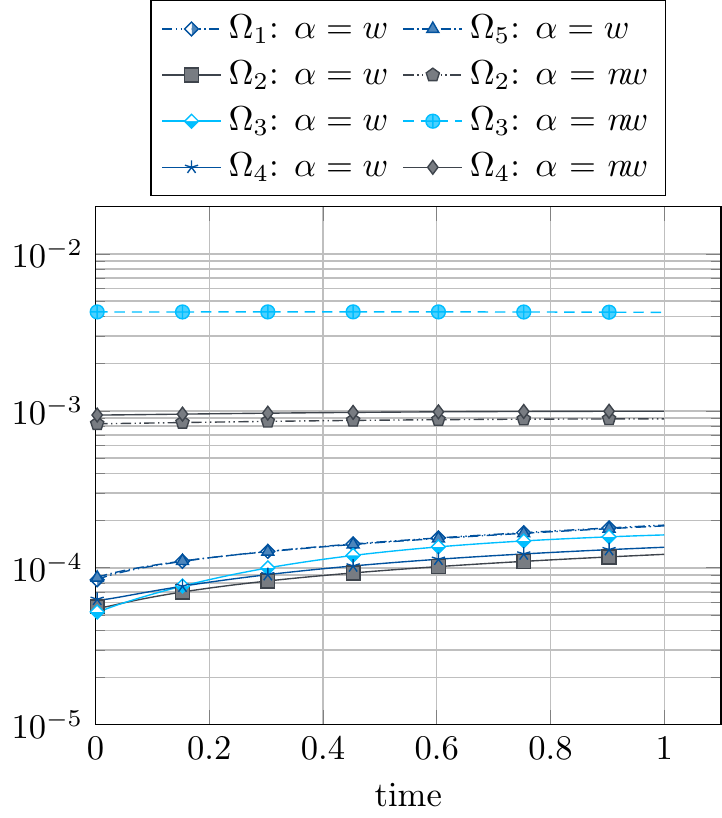}
    \fi
    \caption{Relative error norms $\|\paln{e,n} - \palnh{n}\|_{{}_2}\bigl/\|\paln{e,n}\|_{{}_2}$ over time 
             $t$ for all occurring $(\alpha,l)$.}
    \label{LDD-TPR:fig:TPR:fiveDomainIP:sameIntrinsic:errornorm}
  \end{subfigure}\hspace{4ex}
  \begin{subfigure}[b]{.45\textwidth}
    \centering
    \ifdevelopmentVersion
      \input{\LDDTPRfigSource/2020-07-16-TPR-FiveDomainIP-same-intrinsic_subsequent_iterations.tikz}
    \else
      \includegraphics[width=\textwidth]{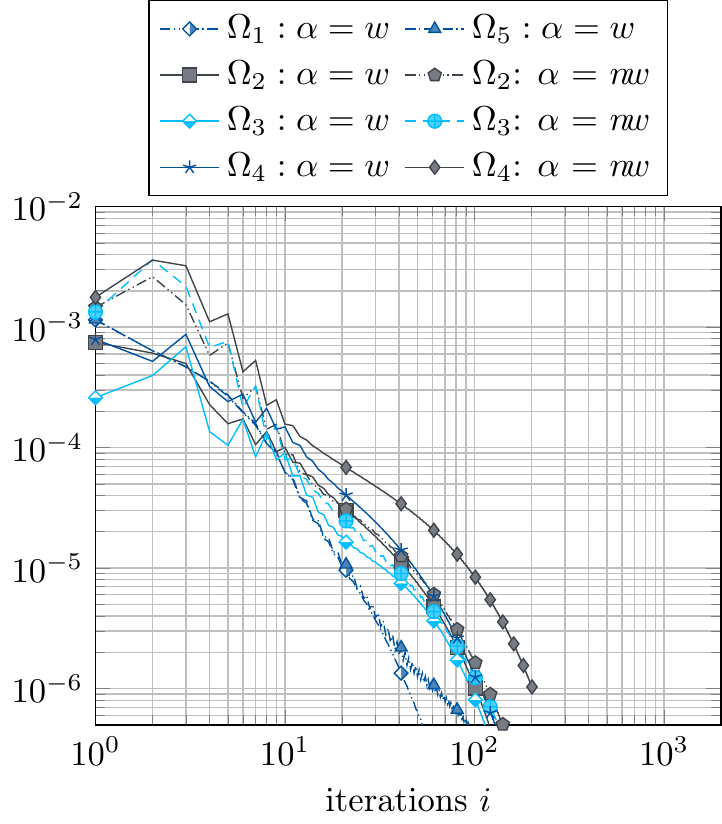}
    \fi
    \caption{Subsequent errors $\|\palnih{i} - \palnih{i-1}\|_{{}_2}$ over iteration $i$ at time step 
             $n=1000$, $t_n = \num{1,0}$. 
             \label{LDD-TPR:fig:TPR:fiveDomainIP:sameIntrinsic:subsequent_errors} }
  \end{subfigure}  
  \caption[\LDDTPRcaptionPrefix five-domain with inner subdomain, same soil parameters]{\TPR coupling on 
          five-domain substructuring with inner subdomain: 
          Relative error norms (\subref{LDD-TPR:fig:TPR:fiveDomainIP:sameIntrinsic:errornorm})
          and subsequent errors at a fixed time step 
          (\subref{LDD-TPR:fig:TPR:fiveDomainIP:sameIntrinsic:subsequent_errors})
          for a simulation over $1000$ time steps with same intrinsic 
          permeabilities and porosities,  \emph{excluding gravity} and parameters 
          $h\approx 0.070$–$0.071$, $\tau = 1\cdot 10^{-3}$, $L_{w,l} = \num{0.01}$, 
          $L_{\nw,l} = \num{0.004}$
          and $\lambda_{w}^{lk} =\num{1}$, $\lambda_{\nw}^{lk} =\num{0,25}$, $l\in\ind$, $k\in\ind_l$.
          \label{LDD-TPR:fig:TPR:fiveDomainIP:sameIntrinsic}
          }
\end{figure}
The nonwetting phase on the inner subdomain shows a degradation in  accuracy whereas the other phase errors are in line with the two-domain examples. 

\autoref{LDD-TPR:fig:TPR:fiveDomainIP:sameIntrinsic}(\subref{LDD-TPR:fig:TPR:fiveDomainIP:sameIntrinsic:subsequent_errors}) shows the subsequent errors for the time step $t_{1000}$. 
The stopping criterion was $\epsilon_s = 1\cdot\num{e-6}$ and was reached after $202$ iterations. 
The first time step required $106$ iterations.

\paragraph{Including gravity}
The behaviour of the solver  when gravity is taken into account
is shown in \autoref{LDD-TPR:fig:TPR:fiveDomainIP:sameIntrinsic:withGravity}(\subref{LDD-TPR:fig:TPR:fiveDomainIP:sameIntrinsic:withGravity:subsequent_errors}). 
To stabilise the solver, it was necessary to adjust the \LDDTPR parameters to 
$L_{\alpha,l} = \num{0.5}$ and $\lambda_{\alpha}^{lk} = 4$ for all phases $\alpha = w, \nw$, $l\in\ind$ and $k\in\ind_l$.
A tilting behaviour in the subsequent error curves can be seen.
This leads to plateaus in the curves and consequently, 
$370$ iterations were required for time step $t_{1000} = \num{1,0}$ to achieve the stopping criterion with 
$\epsilon_s = 5\cdot\num{e-6}$. During the calculation for the first time step,  $193$ iterations were required.
Albeit the solver exhibiting more struggle, the overall approximation quality, despite oscillating a bit seems unaffected as 
\autoref{LDD-TPR:fig:TPR:fiveDomainIP:sameIntrinsic:withGravity}(\subref{LDD-TPR:fig:TPR:fiveDomainIP:sameIntrinsic:withGravity:errornorm}) shows. 
\begin{figure}[t]
  \centering
  \begin{subfigure}[b]{.45\textwidth}
    \centering
    \ifdevelopmentVersion
      \input{\LDDTPRfigSource/2020-06-29-TPR-FiveDomainIP-same-intrinsic-with-g_error_norms.tikz}    
    \else
      \includegraphics[width=\textwidth]{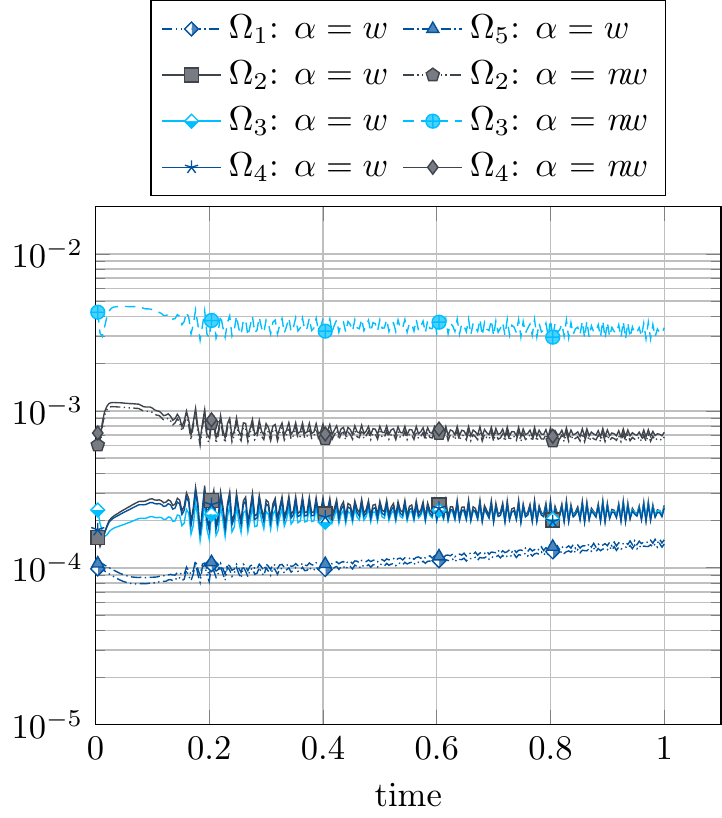}
    \fi
    \caption{Relative error norms $\|\paln{e,n} - \palnh{n}\|_{{}_2}\bigl/\|\paln{e,n}\|_{{}_2}$ over time 
             $t$ for all occurring $(\alpha,l)$.}
    \label{LDD-TPR:fig:TPR:fiveDomainIP:sameIntrinsic:withGravity:errornorm}
  \end{subfigure}\hspace{2ex}
  \begin{subfigure}[b]{.45\textwidth}
    \centering
    \ifdevelopmentVersion
      \input{\LDDTPRfigSource/2020-06-29-TPR-FiveDomainIP-same-intrinsic-with-g_subsequent_iterations.tikz}
    \else
      \includegraphics[width=\textwidth]{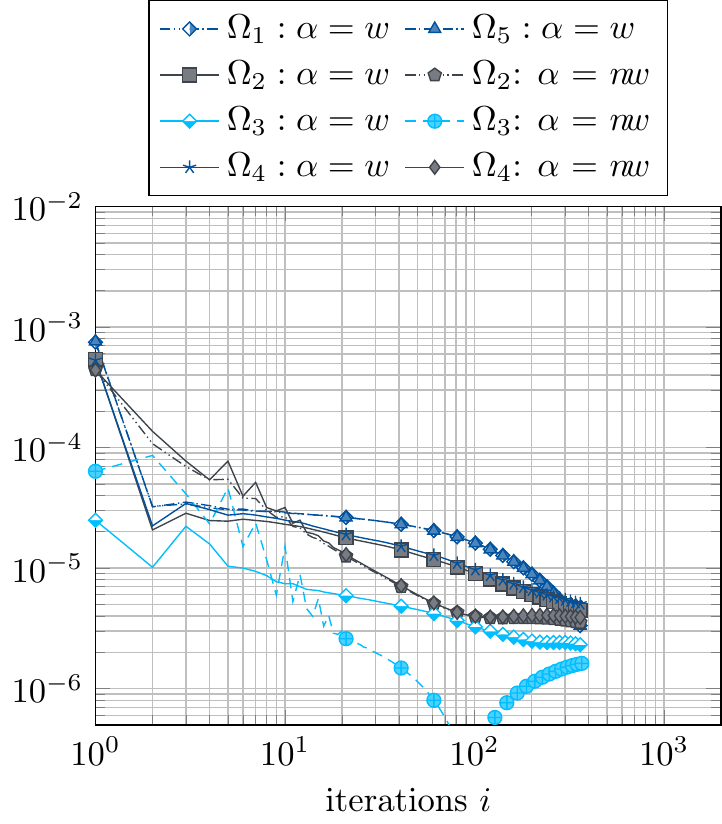}
    \fi
    \caption{Subsequent errors $\|\palnih{i} - \palnih{i-1}\|_{{}_2}$ over iteration $i$ at time step 
             $n=1000$, $t_n = \num{1,0}$. 
             \label{LDD-TPR:fig:TPR:fiveDomainIP:sameIntrinsic:withGravity:subsequent_errors} }
  \end{subfigure}  
  \caption[\LDDTPRcaptionPrefix five-domain with inner subdomain, same soil parameters with gravity]{\TPR coupling on 
          five-domain substructuring with inner subdomain: 
          Relative error norms (\subref{LDD-TPR:fig:TPR:fiveDomainIP:sameIntrinsic:withGravity:errornorm})
          and subsequent errors at a fixed time step 
          (\subref{LDD-TPR:fig:TPR:fiveDomainIP:sameIntrinsic:withGravity:subsequent_errors})
          for a simulation over $1000$ time steps with same intrinsic 
          permeabilities and porosities,  \emph{including gravity} and parameters 
          $h\approx 0.070$–$0.071$, $\tau = 1\cdot 10^{-3}$, $L_{\alpha,l} = \num{0.5}$, 
           and $\lambda_{\alpha}^{lk} =\num{4}$, for all appearing $\alpha \in\{w, \nw\}$  and 
          $l\in\ind$, $k\in\ind_l$.
          \label{LDD-TPR:fig:TPR:fiveDomainIP:sameIntrinsic:withGravity}
          }
\end{figure}
The error of the nonwetting phase of the inner subdomain remains under $5\cdot \num{e-3}$ for all times, and all other errors are lower. 
The tilting behaviour and occurence of plateaus was observed for all examples featuring the inclusion of gravity and is most probably due to the inherent instability of standard finite element methods  for 
advection-dominated regimes.

\section{Conclusions}
In this work we proposed a new domain decomposition approach for hybrid two-phase flow systems.  
For new coupling conditions between domains with different  two-phase flow models  we developped an 
approach combining an $L$-type linearisation of the nonlinearities with a generalised nonoverlapping 
alternating Schwarz method,  the \LDDTPR solver. 
This formulation unifies the work of both, 
\cite{SeusMitra2018} and \cite{SeusEnumath2019} on 
homogeneous two-phase models
and allows for the treatment of complex modelling situations involving very heterogeneous soil
parameters. 
The \LDDTPR solver has been analysed rigorously on the 
time-discrete level. Numerical experiments for two- and multi-domain settings confirm the theoretical 
findings. In particular, they show the possible gain 
of computing time when using the hybrid model instead
of employing  an expensive full two-phase model on the entire domain.

As the \LDDTPR solver linearises and decouples the substructured problem, it  can either be used as a pure domain decomposition method,
as a basis for effective parallel computation,  
or in a model-adaptive domain decomposition setting, in which  an envisioned model change (two-phase/\richards) dictate the substructuring. Future work will be directed to design
such an algorithm that  might also include an adaptive choice of models based on our error
analysis. We envisage that our approach 
is not only effective for the basic two-phase flow 
models encountered here but can be also extended 
to more complex model hierarchies for  multi-phase flow and/or multi-component transport.

\section*{Acknowledgements}

The authors thank the German Research Foundation (DFG) for funding this work (Project Number 327154368 -- SFB 1313).   
In parts, this work was supported by E.ON Stipendienfonds (Project Number T0087/30890/17) which funded a research stay at the University of Bergen (UIB) for which the authors are grateful.

\bibliography{\bibDir/Bibliography.bib}

\begin{thebibliography}{10}
\expandafter\ifx\csname url\endcsname\relax
  \def\url#1{\texttt{#1}}\fi
\expandafter\ifx\csname urlprefix\endcsname\relax\def\urlprefix{URL }\fi
\expandafter\ifx\csname href\endcsname\relax
  \def\href#1#2{#2} \def\path#1{#1}\fi

\bibitem{SeusMitra2018}
D.~Seus, K.~Mitra, I.~S. Pop, F.~A. Radu, C.~Rohde,
  \href{https://doi.org/10.1016/j.cma.2018.01.029}{A linear domain
  decomposition method for partially saturated flow in porous media}, Computer
  Methods in Applied Mechanics and Engineering 333 (2018) 331--355.

\bibitem{SeusEnumath2019}
D.~Seus, F.~A. Radu, C.~Rohde,
  \href{https://doi.org/10.1007/978-3-319-96415-7_55}{A linear domain
  decomposition method for two-phase flow in porous media}, in: F.~A. Radu,
  K.~Kumar, I.~Berre, J.~M. Nordbotten, I.~S. Pop (Eds.), Numerical Mathematics
  and Advanced Applications ENUMATH 2017, Lecture Notes in Computational
  Science and Engineering, Springer International Publishing, 2019, pp.
  603--614.

\bibitem{Lions1988}
P.-L. Lions, {On the {S}chwarz alternating method}, in: R.~Glowinski, G.~H.
  Golub, G.~A. Meurant, J.~Periaux (Eds.), Proceedings of the 1st International
  Symposium on Domain Decomposition Methods for Partial Differential Equations,
  SIAM, Philadelphia, 1988, pp. 1--42.

\bibitem{Pop2004}
I.~S. Pop, F.~A. Radu, P.~Knabner,
  \href{http://www.sciencedirect.com/science/article/pii/S037704270301001X}{Mixed
  finite elements for the {R}ichards’ equation: linearization procedure},
  Journal of Computational and Applied Mathematics 168~(1–2) (2004) 365--373.

\bibitem{List2016}
F.~List, F.~A. Radu, \href{http://dx.doi.org/10.1007/s10596-016-9566-3}{{A
  study on iterative methods for solving Richards' equation}}, Computational
  Geosciences 20~(2) (2016) 341--353.

\bibitem{Agranovich2015}
M.~S. Agranovich, \href{https://www.springer.com/de/book/9783319146478}{Sobolev
  Spaces, Their Generalizations and Elliptic Problems in Smooth and Lipschitz
  Domains}, Springer Monographs in Mathematics, Springer International
  Publishing Switzerland, 2015.

\bibitem{QuarteroniValli2005}
A.~Quarteroni, A.~Valli, Domain decomposition methods for partial differential
  equations, repr. Edition, Numerical mathematics and scientific computation,
  Clarendon Press, Oxford [u.a.], 2005.

\bibitem{Dolean2015}
V.~Dolean, P.~Jolivet, F.~Nataf,
  \href{http://dx.doi.org/10.1137/1.9781611974065.ch1}{An introduction to
  domain decomposition methods}, Society for Industrial and Applied Mathematics
  (SIAM), Philadelphia, PA, 2015, algorithms, theory, and parallel
  implementation.

\bibitem{Benne2016}
D.~Bennequin, M.~J. Gander, L.~Gouarin, L.~Halpern, Optimized {S}chwarz
  waveform relaxation for advection reaction diffusion equations in two
  dimensions, Numerische Mathematik 134 (2016) 513--567.

\bibitem{Skogestad2013}
J.~O. Skogestad, E.~Keilegavlen, J.~M. Nordbotten,
  \href{http://www.sciencedirect.com/science/article/pii/S002199911200589X}{Domain
  decomposition strategies for nonlinear flow problems in porous media},
  Journal of Computational Physics 234 (2013) 439--451.

\bibitem{Berninger2015}
H.~Berninger, R.~Kornhuber, O.~Sander,
  \href{http://dx.doi.org/10.1007/s10596-014-9461-8}{{A multidomain
  discretization of the Richards equation in layered soil}}, Computational
  Geosciences 19~(1) (2015) 213--232.

\bibitem{DIPIETRO2014163}
D.~A. {Di Pietro}, E.~Flauraud, M.~Vohralík, S.~Yousef, A posteriori error
  estimates, stopping criteria, and adaptivity for multiphase compositional
  {D}arcy flows in porous media, Journal of Computational Physics 276 (2014)
  163--187.

\bibitem{Ahmed2019}
E.~Ahmed, S.~A. Hassan, C.~Japhet, M.~Kern, M.~Vohral\'{i}k,
  \href{https://smai-jcm.centre-mersenne.org/articles/10.5802/smai-jcm.47/}{A
  posteriori error estimates and stopping criteria for space-time domain
  decomposition for two-phase flow between different rock types}, The SMAI
  Journal of Computational Mathematics 5 (2019) 195--227.

\bibitem{KURAZ20142}
M.~Kuraz, P.~Mayer, P.~Pech,
  \href{https://www.sciencedirect.com/science/article/pii/S0377042714001502}{Solving
  the nonlinear {R}ichards equation model with adaptive domain decomposition},
  Journal of Computational and Applied Mathematics 270 (2014) 2--11, fourth
  International Conference on Finite Element Methods in Engineering and
  Sciences (FEMTEC 2013).

\bibitem{gander2021nonoverlapping}
M.~Gander, S.~Lunowa, C.~Rohde, Non-overlapping {S}chwarz waveform-relaxation
  for nonlinear advection-diffusion equations, preprint:
  \texttt{http://www.uhasselt.be/Documents/CMAT/Preprints/2021/UP2103.pdf}
  (2021).

\bibitem{Lunowa2020}
S.~B. Lunowa, I.~S. Pop, B.~Koren,
  \href{https://www.sciencedirect.com/science/article/pii/S0045782520305491}{Linearized
  domain decomposition methods for two-phase porous media flow models involving
  dynamic capillarity and hysteresis}, Computer Methods in Applied Mechanics
  and Engineering 372 (2020) 113364.

\bibitem{AHMED2020113294}
E.~Ahmed, C.~Japhet, M.~Kern, Space–time domain decomposition for two-phase
  flow between different rock types, Computer Methods in Applied Mechanics and
  Engineering 371 (2020) 113294.

\bibitem{BorregalesReveron2021}
M.~A. Borregales~Rever\'{o}n, K.~Kumar, J.~M. Nordbotten, F.~A. Radu,
  \href{https://doi.org/10.1007/s10596-020-09983-0}{Iterative solvers for biot
  model under small and large deformations}, Computational Geosciences 25~(2)
  (2021) 687--699.

\bibitem{Illiano2021}
D.~Illiano, I.~S. Pop, F.~A. Radu,
  \href{https://doi.org/10.1007/s10596-020-09949-2}{Iterative schemes for
  surfactant transport in porous media}, Computational Geosciences 25~(2)
  (2021) 805--822.

\bibitem{MitraPop2019}
K.~Mitra, I.~Pop,
  \href{http://www.sciencedirect.com/science/article/pii/S0898122118305546}{A
  modified {L}-scheme to solve nonlinear diffusion problems}, Computers \&
  Mathematics with Applications 77~(6) (2019) 1722 -- 1738, 7th International
  Conference on Advanced Computational Methods in Engineering (ACOMEN 2017).

\bibitem{Helmig97}
R.~Helmig, \href{http://swbplus.bsz-bw.de/bsz061703095err.htm}{Multiphase flow
  and transport processes in the subsurface: a contribution to the modeling of
  hydrosystems}, Springer, Berlin, 1997.

\bibitem{Bear2007}
J.~Bear, Hydraulics of Groundwater, Dover Books on Engineering, Dover
  Publications, 2007.

\bibitem{Richards1931}
L.~A. Richards, \href{http://dx.doi.org/10.1063/1.1745010}{Capillary conduction
  of liquids through porous mediums}, Physics 1~(5) (1931) 318--333.

\bibitem{Richardson1922}
L.~F. Richardson,
  \href{https://archive.org/stream/weatherpredictio00richrich#page/n0/mode/2up/search/Darcy}{Weather
  prediction by numerical process}, Camebridge University Press, 1922.

\bibitem{HenryHilhorstEymard2012}
M.~Henry, D.~Hilhorst, R.~Eymard,
  \href{https://doi.org/10.3934/dcdss.2012.5.93}{Singular limit of a two-phase
  flow problem in porous medium as the air viscosity tends to zero}, Discrete
  and Continous Dynamical Systems - Series S (DCDS - S) 5~(1) (2012) 93--113.

\bibitem{Cances2012}
C.~Cancès, M.~Pierre, \href{https://doi.org/10.1137/11082943X}{An existence
  result for multidimensional immiscible two-phase flows with discontinuous
  capillary pressure field}, SIAM Journal on Mathematical Analysis 44~(2)
  (2012) 966--992.

\bibitem{Berninger2009}
H.~Berninger,
  \href{http://www.diss.fu-berlin.de/diss/receive/FUDISS_thesis_000000008972}{{Domain
  decomposition methods for elliptic problems with jumping nonlinearities and
  application to the Richards equation}}, Ph.D. thesis, FB Mathematik und
  Informatik, Freie Universit\"{a}t Berlin (2009).

\bibitem{Brezzi1991}
F.~Brezzi, M.~Fortin, {Mixed and Hybrid Finite Element Methods}, Vol.~15 of
  Springer Series in Computational Mathematics, Springer, 1991.

\bibitem{MacLean2000}
W.~MacLean,
  \href{http://bvbr.bib-bvb.de:8991/F?func=service&doc_library=BVB01&doc_number=009028341&line_number=0002&func_code=DB_RECORDS&service_type=MEDIA}{Strongly
  elliptic systems and boundary integral equations}, 1st Edition, Cambridge
  University Press, Cambridge [u.a.], 2000.

\bibitem{Lions1990III}
P.~L. Lions, On the {Schwarz} alternating method {III}: A variant for
  nonoverlapping subdomains, in: T.~F. Chan, R.~Glowinski, J.~P\'{e}riaux,
  O.~B. Widlund (Eds.), Third International Symposium on Domain Decomposition
  Methods for Partial Differential Equations, Society for Industrial and
  Applied Mathematics, 1990, pp. 202--223.

\bibitem{SeusThesis2021}
D.~Seus, {LDD} schemes for two-phase flow systems, Ph.D. thesis, University of
  Stuttgart (2021).

\bibitem{Radu2017}
F.~A. Radu, K.~Kumar, J.~M. Nordbotten, I.~S. Pop,
  \href{https://doi.org/10.1093/imanum/drx032}{A robust, mass conservative
  scheme for two-phase flow in porous media including {H}\"{o}lder continuous
  nonlinearities}, IMA Journal of Numerical Analysis 38 (2017) 884–920.

\bibitem{Fenics:AlnaesBlechta2015a}
M.~S. Aln{\ae}s, J.~Blechta, J.~Hake, A.~Johansson, B.~Kehlet, A.~Logg,
  C.~Richardson, J.~Ring, M.~E. Rognes, G.~N. Wells, The {FE}ni{CS} project
  version 1.5, Archive of Numerical Software 3~(100).

\bibitem{Fenics:LoggWellsEtAl2012a}
A.~Logg, G.~N. Wells, J.~Hake, DOLFIN: a C++/Python Finite Element Library,
  Springer, 2012, Ch.~10.

\bibitem{LDDcode}
D.~Seus, \texttt{https://gitlab.com/davidseus/ldd-for-two-phase-flow-systems}
  (2021).

\end{thebibliography}
\end{document}